\documentclass[10pt,a4]{article}
\usepackage[english]{babel}
\usepackage[dvips,ps2pdf,colorlinks,linkcolor=black,citecolor=black]{hyperref}
\usepackage[hmargin=1.5cm, vmargin=2cm]{geometry} 
\usepackage{amsmath,amssymb}
\usepackage{array}
\usepackage{pst-all}
\usepackage{ntheorem}
\usepackage{stmaryrd}
\usepackage{pgf}

\title{Admissible Diagrams in $ $ $U_{q}^{w}(\mathfrak{g})$ $ $ and Combinatoric Properties \\ of Weyl Groups.}
\author{Antoine M\'eriaux, G\'erard Cauchon\\
\\
\it Laboratoire d'\'equations aux d\'eriv\'ees partielles et physique
math\'ematique,\\ \it U.F.R.
Sciences, B.P. 1039, 51687 Reims Cedex 2, France.}
\date{\ }
\newenvironment{itemizedot}{\begin{itemize} }{\end{itemize}}

\begin{document}
\maketitle

\bigskip
\def\Le{\hbox{\rotatedown{$\Gamma$}}}
\newcommand{\bqwg}{$B_{q}^{w}(\mathfrak{g}) $}
\newcommand{\rqg}{$R_{q}[\bf G] \rm $}
\newcommand{\G}{$ \bf G \rm $}
\newcommand{\uqtg}{$ \check{U}_{q}(\mathfrak{g}) $}
\newcommand{\uqbp}{$U_{q}(\mathfrak{b}^{+}) $}
\newcommand{\uqbm}{$U_{q}(\mathfrak{b}^{-}) $}
\newcommand{\uqtbp}{$\check{U}_{q}(\mathfrak{b}^{+}) $}
\newcommand{\uqtbm}{$\check{U}_{q}(\mathfrak{b}^{-}) $}
\newcommand{\uqnp}{$U_{q}(\mathfrak{n}^{+}) $}
\newcommand{\uqnm}{$U_{q}(\mathfrak{n}^{-}) $}
\newcommand{\bqtwg}{$\check{B}_{q}^{w}(\mathfrak{g})$}

\newcommand{\g}{\mathfrak{g}}
\newcommand{\uqg}{U_{q}(\mathfrak{g})}
\newcommand{\uqpg}{U_{q}^{+}(\mathfrak{g})}
\newcommand{\upl}{U^{+}}
\newcommand{\upw}{U^{+}[w]}

\font\bbfnt=msbm10
\newcommand{\bbN}{\hbox{\bbfnt{\char'116}}}
\newcommand{\bbZ}{\hbox{\bbfnt{\char'132}}}
\newcommand{\bbC}{\hbox{\bbfnt{\char'103}}}
\newcommand{\bbQ}{\hbox{\bbfnt{\char'121}}}

\newcommand{\spn}{\hbox{\bbfnt{\char'045}}}
\newcommand{\sbn}{\hbox{\bbfnt{\char'044}}}
\newcommand{\sr}{$Spec_{w}(R)$}
\newcommand{\srb}{$Spec_{w}(\overline{R})$}
\newcommand{\ubqwg}{$\overline{U}_{q}^{w}(\mathfrak{g})$}
\newcommand{\ucqwg}{$\widehat{U}_{q}^{w}(\mathfrak{g})$}
\newcommand{\eoig}{$k[X_{1}][X_{2};\sigma _{2},\delta_{2}]\cdots
[X_{N};\sigma _{N},\delta_{N}]$}
\newcommand{\Rb}{$\overline{R}$}
\newcommand{\jzj}{$\mathcal{P} \it _{j}^{0} (R_{j})$}
\newcommand{\jzjp}{$\mathcal{P} \it _{j}^{0}(R_{j+1})$}
\newcommand{\juj}{$\mathcal{P} \it _{j}^{1}(R_{j})$}
\newcommand{\jujp}{$\mathcal{P} \it _{j}^{1}(R_{j+1})$}
\newcommand{\xb}{$\overline{x}$}
\newcommand{\Ab}{$\overline{A}$}
\newcommand{\vf}{$\varphi$}
\newcommand{\ds}{$\displaystyle$}
\newcommand{\wb}{$\overline{w}\mbox{ }$}
\newcommand{\bx}{\begin{flushright}
$\square$
\end{flushright}}

$ $ \\ $ $ \\ $ $ \\ $ $ \\ $ $ \\

\begin{center}

\bf Abstract

\end{center}

\begin{flushleft}
Consider a complex simple Lie algebra $ $ $\g$ $ $ of rank $ $ $n$. Denote by $ $ $\Pi$ $ $ a system of simple roots, by $ $ $W$ $ $ the corresponding Weyl group, consider a reduced expression $ $ $ $ $w \mbox{ } = \mbox{ } s_{\alpha_{1}} \circ \mbox{ } ... \mbox{ } \circ s_{\alpha_{t}}$ $ $ $ $ (each $ $ $\alpha_{i}$ $\in$ $\Pi$) $ $ of some $ $ $w$ $\in$ $W$ $ $ and call \it diagram \rm any subset of $\llbracket  1,$ ... ,$t \rrbracket $. $ $ We denote by $ $ $\upw$ $ $ (or $ $ $U_{q}^{w}(\mathfrak{g})$) $ $ the "quantum nilpotent" algebra defined as in \cite{J}. \\ $ $ \\

We prove (theorem 5.3. 1) that the \it positive diagrams \rm naturally associated with the positive subexpressions (of the reduced expression of $ $ $w$) in the sense of R. Marsh and K. Rietsch \cite{MR}, coincide with the \it admissible diagrams \rm constructed by G. Cauchon \cite{C} which describe the natural stratification of $ $ $Spec(\upw).$ \\ $ $ 

This theorem implies in particular (corollaries 5.3. 1 and 5.3. 2):

\begin{enumerate}

\item The map $ $ $\zeta$ : $ $ $ $ $\Delta$ = $\{j_{1}$ $<$ $...$ $<$ $j_{s}\}$ $ $ $\mapsto$ $ $ $u$ =  $s_{\alpha_{j_{1}}}$ $\circ$ $...$ $\circ$ $s_{\alpha_{j_{s}}}$ $ $ $ $ is a bijection from the set of admissible diagrams onto the set $ $ $ $ $\{u \in W \mbox{ } | \mbox{ } u \mbox{ } \leq \mbox{ } w \}$.

\item For each admissible diagram $ $ $\Delta$ = $\{j_{1}$ $<$ $...$ $<$ $j_{s}\}$, $ $ $ $ $s_{\alpha_{j_{1}}}$ $\circ$ $...$ $\circ$ $s_{\alpha_{j_{s}}}$ $ $ is a reduced expression of $ $ $u$ = $\zeta (\Delta)$. 

\item The map $ $ $\zeta^{\prime}$ : $ $ $ $ $\Delta$ = $\{j_{1}$ $<$ $...$ $<$ $j_{s}\}$ $ $ $\mapsto$ $ $ $u^{\prime}$ =  $s_{\alpha_{j_{s}}}$ $\circ$ $...$ $\circ$ $s_{\alpha_{j_{1}}}$ $ $ $ $ is a bijection from the set of admissible diagrams onto the set $ $ $ $ $\{u \in W \mbox{ } | \mbox{ } u \mbox{ } \leq \mbox{ } v \mbox{ } = \mbox{ } w^{-1} \}$.

\item For each admissible diagram $ $ $\Delta$ = $\{j_{1}$ $<$ $...$ $<$ $j_{s}\}$, $ $ $ $ $s_{\alpha_{j_{s}}}$ $\circ$ $...$ $\circ$ $s_{\alpha_{j_{1}}}$ $ $ is a reduced expression of $ $ $u^{\prime}$ = $\zeta^{\prime} (\Delta)$. 

\end{enumerate}

If the Lie algebra $ $ $\g$ $ $ is of type $ $ $A_{n}$ $ $ and $ $ $w$ $ $ is choosen in order that $ $ $\upw$ $ $ is the quantum matrices algebra $ $ $O_{q}(M_{p,m}(k))$ $ $ with $ $ $m$ = $n-p+1$ (see section 2.1), then, by \cite{CC}, the admissible diagrams are known to be the $ $ \Le $ $ - diagrams in the sense of A. Postnikov \cite{P}. In this particular case,  the equality of Le - diagrams and positive subexpressions (of the reduced expression of $w$) have also been proved (with quite different methods) by A. Postnikov (\cite{P}, theorem 19.1.) and by T. Lam and L. Williams (\cite{LW}, theorem 5.3.).

\newpage

\tableofcontents

$ $ \\
$ $ \\
\section{Introduction.} $ $\\

Given a complex simple Lie algebra $ $ $\g$ $ $ of rank $ $ $n$, $ $ denote by $ $ $\Pi$ $ $ a system of simple roots, by $ $ $W$ $ $ the corresponding Weyl group, and consider a reduced expression $ $ $w \mbox{ } = \mbox{ } s_{\alpha_{1}} \circ \mbox{ } ... \mbox{ } \circ s_{\alpha_{t}}$ $ $ of some $ $ $w$ $\in$ $W$.\\
Denote by $ $ $\upw$ $ $ (or $ $ $U_{q}^{w}(\mathfrak{g})$) $ $ the "quantum nilpotent" algebra associated to $ $ $w$ $ $ as in \cite{J} and by $ $ $X_{1}$, ... , $X_{t}$ $ $ the canonical generators of $ $ $\upw$ $ $ associated with this choosen reduced expression of $ $ $w$ $ $ (see section 2.1 for more details). \\
The natural action of the torus on the prime spectrum $ $ $Spec(\upw)$ $ $ induces a finite stratification which is completely described by the \it admissible diagrams \rm (also called Cauchon diagrams) which are some particular subsets of $ $ $\llbracket  1,$ ... ,$ $ $t \rrbracket $ $ $ (see \cite{C}). Those admissible diagrams are very closely connected with the quantum commutation rules satisfied by the canonical generators. For example, in the $ $ $A_{n}$ $ $ type, one can choose $ $ $w$ $ $ such that $ $ $\upw$ $ $ is a quantum matrices algebra
$$O_{q}(M_{p,m}(k)) \mbox{ } = \mbox{ } k<X_{i,j}> \mbox{ } \mbox{ } \mbox{ } (1 \mbox{ } \leq \mbox{ } i \leq \mbox{ } p, \mbox{ } 1 \mbox{ } \leq \mbox{ } j \leq \mbox{ } m)$$
and we know \cite{CC} that the commutation relations
$$X_{u,v}X_{i,j} \mbox{ } = \mbox{ } X_{i,j}X_{u,v} \mbox{ } - \mbox{ } (q-q^{-1})X_{i,v}X_{u,j} \mbox{ } \mbox{ } \mbox{ } (i \mbox{ } < \mbox{ } u, \mbox{ } j \mbox{ } < \mbox{ } v) \mbox{ } \mbox{ } \mbox{ } \mbox{ } \mbox{ } \mbox{ } \mbox{ } \mbox{ } (*)$$
determine the shape of admissible diagrams by the following observation: \\ If $ $ $\mathcal{P}$ $ $ is any completely prime ideal of $ $ $O_{q}(M_{p,m}(k))$, $ $ we have immediately \\
$$X_{u,v} \mbox{ } \in  \mbox{ } \mathcal{P} \it \mbox{ } \Rightarrow  \mbox{ } (X_{i,v} \mbox{ } \in  \mbox{ } \mathcal{P} \it \mbox{ } \mbox{ }  \rm or \it \mbox{ } \mbox{ } X_{u,j} \mbox{ } \in  \mbox{ } \mathcal{P}) \mbox{ } \mbox{ } \mbox{ } \mbox{ } \mbox{ } \mbox{ } \mbox{ } \mbox{ } (**)$$
This implies \cite{CC} that $ $ $\Delta$ $\subset$ $\llbracket 1,  \mbox{ } ...  \mbox{ } ,  \mbox{ } p \rrbracket $ $\times$ $\llbracket 1,  \mbox{ } ...  \mbox{ } ,  \mbox{ } m \rrbracket $ $ $ (the set of indexes) is an admissible diagram if and only if, for each $ $ $i, \mbox{ } j, \mbox{ } u, \mbox{ } v$ $ $ as over, we have:
$$(u,v) \mbox{ } \in  \mbox{ } \Delta \it \mbox{ } \Rightarrow  \mbox{ } ((i,v) \mbox{ } \in  \mbox{ } \Delta \mbox{ } \mbox{ } \mbox{ } \mbox{ } \rm or \mbox{ } \mbox{ } \mbox{ }\it \mbox{ } (u,j) \mbox{ } \in  \mbox{ } \Delta)$$
This means that the admissible diagrams are the unions of truncated rows and columns, namely the $ $ \Le $ $ - diagrams in the sense of A. Postnikov (see \cite{P}). \\

So, we see on this example that admissible diagrams are quantum objects since, in the non quantum case (when $ $ $q$ = $1$), $ $ the formulas $ $ $(*)$ $ $ become $ $ $X_{u,v}X_{i,j} \mbox{ } = \mbox{ } X_{i,j}X_{u,v}$, $ $ so that observation $ $ $(**)$ $ $ is not valid any longer and the admissible diagrams become invisible. \\ $ $ \\

On the other hand, R. Marsh and K. Rietsch \cite{MR} defined the notion of positive subexpression of the reduced $ $ $w$'s $ $ expression considered over. These positive subexpressions are defined by particular subsets of $ $ $\llbracket  1,$ ... ,$ $ $t \rrbracket $ $ $ that we call the \it positive diagrams. \rm R. Marsh and K. Rietsch proved in \cite{MR} that they are in one to one natural correspondence with the elements of the Weyl group which are smaller or equal to $ $ $w$ $ $ (for the Bruhat order). In this paper, we prove (theorem 5.3. 1) that \it the positive diagrams coincide with the admissible diagrams, \rm which can be interpreted saying that \it R. Marsh and K. Rietsch positive subexpressions are quantum objects. \\ \rm $ $ \\
In particular, this implies (corollary 5.3. 1):

\begin{enumerate}

\item The map $ $ $\zeta$ : $ $ $ $ $\Delta$ = $\{j_{1}$ $<$ $...$ $<$ $j_{s}\}$ $ $ $\mapsto$ $ $ $u$ = $s_{\alpha_{j_{1}}}$ $\circ$ $...$ $\circ$ $s_{\alpha_{j_{s}}}$ $ $ $ $ is a bijection from the set of admissible diagrams onto the set $ $ $ $ $\{u \in W \mbox{ } | \mbox{ } u \mbox{ } \leq \mbox{ } w \}$,

\item Consider an admissible diagram $ $ $\Delta$ = $\{j_{1}$ $<$ $...$ $<$ $j_{s}\}$ $ $ and some integer $ $ $i$ $\in$ $\llbracket  1,  \mbox{ } ...  \mbox{ } ,  \mbox{ } t \rrbracket $. $ $ Set $ $ $\Delta$ $\cap$ $\llbracket i + 1,  \mbox{ } ...  \mbox{ } ,  \mbox{ } t \rrbracket $ = $\{j_{c}$, $...$ , $j_{s}\}$ $ $ $(1 \mbox{ } \leq \mbox{ } c \mbox{ } \leq \mbox{ } s)$. $ $ Then the expression $ $ $s_{\alpha_{i}}$ $\circ$ $s_{\alpha_{j_{c}}}$ $\circ$ $...$ $\circ$ $s_{\alpha_{j_{s}}}$ $ $ is reduced. In particular, $s_{\alpha_{j_{1}}}$ $\circ$ $...$ $\circ$ $s_{\alpha_{j_{s}}}$ $ $ is a reduced expression of $ $ $\zeta (\Delta)$, 

\item Consider some $ $ $u$ $\in$ $W$ $ $ with $ $ $u$ $\leq$ $w$. $ $  Then, the only admissible diagram $ $ $\Delta$ $ $ such that $ $ $\zeta (\Delta)$ = $u$ is recursively defined as follows: \\
$\bullet$ $ $ $1$ $\in$ $\Delta$ $ $ $\Leftrightarrow$ $ $ $l(s_{\alpha_{1}}$ $\circ$ $u)$ = $l(u)$ $-$ $1$ $ $ $\Leftrightarrow$ $ $ $u^{-1}(\alpha_{1})$ $ $ is a negative root, \\

$\bullet$ $ $ Consider some integer $ $ $i$ $\in$ $\llbracket  1,  \mbox{ } ...  \mbox{ } ,  \mbox{ } t - 1 \rrbracket $, $ $ assume that $ $ $\Delta$ $\cap$ $\llbracket 1,  \mbox{ } ...  \mbox{ } ,  \mbox{ } i \rrbracket $ = $\{j_{1}$, $...$ , $j_{d}\}$ $ $ and set $ $ $u_{i}$ = $s_{\alpha_{j_{d}}}$ $\circ$ $...$ $\circ$ $s_{\alpha_{j_{2}}}$ $\circ$ $s_{\alpha_{j_{1}}}$ $\circ$ $u$ $ $ $(u_{i}$ = $u$ $ $ if $ $ $\Delta$ $\cap$ $\llbracket 1,  \mbox{ } ...  \mbox{ } ,  \mbox{ } i \rrbracket $  = $\emptyset)$. \\ Then 
$ $ ($i + 1$ $\in$ $\Delta$) $ $ $\Leftrightarrow$ $ $ ($l(s_{\alpha_{i+1}}$ $\circ$ $u_{i})$ = $l(u_{i})$ $-$ $1$) $ $ $\Leftrightarrow$ $ $ ($u_{i}^{-1}(\alpha_{i + 1})$ $ $ is a negative root).

\end{enumerate}
$ $ \\
It also implies (corollary 5.3. 2): \\ 

\begin{enumerate}

\item The map $ $ $\zeta^{\prime}$ : $ $ $ $ $\Delta$ = $\{j_{1}$ $<$ $...$ $<$ $j_{s}\}$ $ $ $\mapsto$ $ $ $u^{\prime}$ = $s_{\alpha_{j_{s}}}$ $\circ$ $...$ $\circ$ $s_{\alpha_{j_{1}}}$ $ $ $ $ is a bijection from the set of admissible diagrams onto the set $ $ $ $ $\{u \in W \mbox{ } | \mbox{ } u \mbox{ } \leq \mbox{ } v \mbox{ } = \mbox{ } w^{-1} \}$. 

\item Consider an admissible diagram $ $ $\Delta$ = $\{j_{1}$ $<$ $...$ $<$ $j_{s}\}$ $ $ and some integer $ $ $i$ $\in$ $\llbracket  1,  \mbox{ } ...  \mbox{ } ,  \mbox{ } t \rrbracket $. $ $ Set $ $ $\Delta$ $\cap$ $\llbracket i + 1,  \mbox{ } ...  \mbox{ } ,  \mbox{ } t \rrbracket $ = $\{j_{c}$, $...$ , $j_{s}\}$ $ $ $(1 \mbox{ } \leq \mbox{ } c \mbox{ } \leq \mbox{ } s)$. $ $ Then the expression $ $ $s_{\alpha_{j_{s}}}$ $\circ$ $...$ $\circ$ $s_{\alpha_{j_{c}}}$ $\circ$ $s_{\alpha_{i}}$ $ $ is reduced. In particular, $s_{\alpha_{j_{s}}}$ $\circ$ $...$ $\circ$ $s_{\alpha_{j_{1}}}$ $ $ is a reduced expression of $ $ $\zeta^{\prime} (\Delta)$. 

\end{enumerate}

In the particular case $ $ $\upw$ = $O_{q}(M_{p,m}(k))$ $ $  mentioned over, this last result gives back a theorem of A. Postnikov \cite{P} and a theorem of T. Lam and L. Williams \cite{LW}. \\ $ $ \\

In the general case, let us denote by $ $ $\Delta$ $\leftrightarrow$ $w^{\Delta}$ $ $ the one to one correspondence constructed by R. Marsh and K. Rietsch \cite{MR} between the positive diagrams and the elements of the Weyl group which are smaller or equal to $ $ $w$. \\
Let us also denote by $ $ $\Delta$ $\leftrightarrow$ $\mathcal{P}\it_{\Delta}$ $ $ the one to one correspondence constructed by G. Cauchon \cite{C} between the admissible diagrams and the prime ideals of $ $ $\upw$ $ $ which are invariant under the torus action. \\
By the theorem 5.3. 1 of this paper, it turns out that there exists a natural one to one correspondence

$$w^{\Delta} \mbox{ } \leftrightarrow \mbox{ } \mathcal{P}\it_{\Delta}$$

between the elements of the Weyl group which are smaller or equal to $ $ $w$ $ $ and the prime ideals of $ $ $\upw$ $ $ which are invariant under the torus action. \\
Let us recall that, in the particular case $ $ $\upw$ = $O_{q}(M_{p,m}(k))$ $ $  mentioned over, S. Launois \cite{L} has constructed (with quite different methods) such a one to one correspondence which, moreover, preserves the ordering (where the Weyl group is provided with the Bruhat order and the spectrum $ $ $Spec(O_{q}(M_{p,m}(k)))$ $ $ is provided with the inclusion order). This leads us to ask the following questions (unsolved at the moment): \\ $ $ \\

\begin{description}
\item[$\bullet$ $ $ \it Question 1.] Does this Launois correspondence coincide with ours in the particular case $ $ $\upw$ = $O_{q}(M_{p,m}(k))$? \\ $ $ \\
\item[$\bullet$ $ $ \it Question 2.] Does our correspondence $ $ $w^{\Delta} \mbox{ } \leftrightarrow \mbox{ } \mathcal{P}\it_{\Delta}$ $ $ preserve the ordering in the general case (where the Weyl group is provided with the Bruhat order and the spectrum $ $ $Spec(\upw)$ $ $ is provided with the inclusion order)? \\ $ $ \\
\end{description}

(A positive answer would supply the Bruhat order with a nice quantum interpretation.)

$ $ \\
\end{flushleft}

\newpage

\begin{flushleft}
\section{Background on Weyl groups.}
$ $\\
Following J. C. Jantzen (\cite{J}), we use the following conventions all along this paper. \\ $ $\\
\begin{itemizedot}
\item $\g$ is a complex simple Lie algebra.
\item $\Phi$ is the (irreducible) root system of $\g$ with respect to
a fixed Cartan subalgebra.
\item $\Pi$ is a fixed basis of $\Phi$, $ $ $\Phi^{+}$ denotes the subset of positive roots, and we set 
$n$ = $|\Pi|$, $ $ N = $|\Phi^{+}|$ $ $ ($1$ $\leq$ $n$ $\leq$ $N$).
\item $W$ is the Weyl group of $\Phi$ and ( , ) is the unique scalar product
on the real vector space $V$ generated by $\Phi$, such that $\|\beta\|^{2}$ = $2$ $ $ ($\|\beta\|$ := $\sqrt{(\beta,\beta})$) $ $ for all short roots $\beta$ in $\Phi$.
\item For any $\beta$ in $\Phi$, we set $ $ $d_{\beta}$ = $\displaystyle \frac{\|\beta\|^{2}}{2}$, $ $ $\beta^{\vee}$ = $\displaystyle \frac{\beta}{d_{\beta}}$,$ $ and $s_{\beta}$ denotes the reflection with respect to $\beta$ $ $ $ $ $ $ \\
($\displaystyle s_{\beta}(x)$ = $x$ $-$ ($\beta^{\vee}$, $x$)$\beta$ $ $ for any $x$ in $V$).
\item $k$ is a field with char($k$) $\neq$ $2$ and, in addition, char($k$) $\neq$ $3$ if $\Phi$ has type $G_{2}$.
\item $q$ $ $ $\in$ $ $ $k^{*} := k \setminus \{ 0 \}$, $ $ $q$ is not a root of unity.
\item The $k$-algebra $\uqg$ and it's canonical generators $E_{\alpha}$, $F_{\alpha}$, $K_{\alpha}^{\pm 1}$   ($\alpha$ $\in$ $\Pi$) are defined as in \cite{J}. We denote by $\upl$ $ $ (or $ $ $U_{q}^{+}(\mathfrak{g})$) $ $ the subalgebra generated by the $E_{\alpha}$ ($\alpha$ $\in$ $\Pi$).
\item Denote by \bbZ$\Pi$ the root lattice. The algebra $\uqg$ is \bbZ$\Pi$-graded and, if $ $ $\alpha$ in $\Pi$, $ $ deg($E_{\alpha}$) = $\alpha$, $ $ deg($F_{\alpha}$) = $-$ $\alpha$, $ $ deg($K_{\alpha}^{\pm 1}$) = $0$.
\item For any $\rho$ = $\displaystyle \sum_{\alpha \mbox{ } \in \mbox{ } \Pi} m_{\alpha}\alpha$ ($m_{\alpha}$ $\in$ \bbZ) in the root lattice \bbZ$\Pi$, we set $K_{\rho}$ = $\displaystyle \prod_{\alpha \mbox{ } \in \mbox{ } \Pi} K_{\alpha}^{m_{\alpha}}$. The multiplicative group $T$ = $\{K_{\rho}$ $|$ $\rho$ $\in$ \bbZ$\Pi \}$ is called the Torus (of $\uqg$). This group acts on the algebra $\uqg$ by
$$K_{\rho}.u \mbox{ }  = \mbox{ } K_{\rho}^{-1}uK_{\rho} \mbox{ } \mbox{ } \mbox{ } \mbox{ }  (\forall u \in U_{q}(\mathfrak{g}))$$
In particular, for any homogeneous element $u$ of $\uqg$ with degree deg($u$) = $\gamma$ $\in$ \bbZ$\Pi$, we have $K_{\rho}$.$u$ = $q^{-(\rho,\gamma)}u$.
\end{itemizedot}

\subsection{Quantum algebras $\upw$.}
Consider any $w$ $\in$ $W$, set $t$ = $l(w)$ and consider a reduced expression
\begin{equation}\label{eqn:expw}
w = s_{\alpha_{1}} \circ \mbox{ } ... \mbox{ } \circ s_{\alpha_{t}} \mbox{ } \mbox{ } \mbox{ } \mbox{ } \mbox{ } (\alpha_{i} \in \Pi \mbox{ }  \mbox{ }  for \mbox{ }  \mbox{ }  1 \leq i \leq t)
\end{equation}
It is well known that $$\beta_{1} = \alpha_{1},  \mbox{ } \beta_{2} = s_{\alpha_{1}}(\alpha_{2}), \mbox{ }  ...\mbox{ }, \mbox{ }  \beta_{t} = s_{\alpha_{1}} \mbox{ }  ... \mbox{ }  s_{\alpha_{t-1}}(\alpha_{t})$$ are distinct positive roots and that the set $\{\beta_{1},  \mbox{ } ... , \mbox{ }  \beta_{t}\}$ does not depend on the reduced expression (\ref{eqn:expw}) of $w$. For any $\alpha \mbox{ } \in \mbox{ } \Pi$, define the braid automorphism $T_{\alpha}$ of the algebra $\uqg$ as in (\cite{J}, p. 153), and set $$X_{\beta_{1}} = E_{\alpha_{1}}, \mbox{ }  X_{\beta_{2}} = T_{\alpha_{1}}(E_{\alpha_{2}}), \mbox{ }  ...\mbox{ }, \mbox{ }  X_{\beta_{t}} = T_{\alpha_{1}} \mbox{ }  ... \mbox{ }  T_{\alpha_{t-1}}(E_{\alpha_{t}}).$$ The following results are classical (\cite{J}, chapter 8): \\ $ $ \\
\begin{itemizedot}
\item $X_{\beta_{1}}$, ... , $X_{\beta_{t}}$ are all in $\upl$. Each $X_{\beta_{i}}$ ($ 1 \leq i \leq t$) is homogeneous of degree deg($X_{\beta_{i}}$) = $\beta_{i}$.
\item We denote by  $ $ $\upw$ $ $ (or $ $ $U_{q}^{w}(\mathfrak{g})$) $ $ the subalgebra of $ $ $\upl$ generated by $X_{\beta_{1}}$, ... , $X_{\beta_{t}}$. This algebra does not depend on the above reduced expression (\ref{eqn:expw}) of $w$ (although the variables $X_{\beta_{1}}$, ... , $X_{\beta_{t}}$ depend of (\ref{eqn:expw})).
\item The ordered monomials $ $ $X^{\underline{a}}$ $ $ := $X_{\beta_{1}}^{a_{1}}$ ...  $X_{\beta_{t}}^{a_{t}}$, $ $ $\underline{a}$ = ($a_{1}$, ... , $a_{1}$) $\in$ \bbN $^{t}$, are a basis of $\upw$ (as a $k$-vector space).
\item Since the above generators $X_{\beta_{i}}$ are homogeneous, the \bbZ$\Pi$-graduation of $\uqg$ induces a \bbZ$\Pi$-graduation of $\upw$ and the action of the Torus $T$ on $\uqg$ induces, by restriction, an action of $T$ on $\upw$.
\item If $1$ $\leq$ $i$ $<$ $j$ $\leq$ $t$, we have the following straightening formula due to Levendorskii and Soibelman:
\begin{equation}\label{LSb}
X_{\beta_{j}}X_{\beta_{i}} - q^{-(\beta_{i},\beta_{j})}X_{\beta_{i}}X_{\beta_{j}} \mbox{ } = \mbox{ } P_{j,i},
\end{equation} \\
\begin{equation}\label{LSPb}
P_{j,i} \mbox{ } = \mbox{ } \sum_{ \underline{a} \mbox{ } = \mbox{ } (a_{i+1}, ... ,a_{j-1})} c_{\underline{a}}X_{\beta_{i+1}}^{a_{i+1}} ... X_{\beta_{j-1}}^{a_{j-1}}.
\end{equation} \\
with $\underline{a}$ $\in$ \bbN$^{j-i-1}$, $\mbox{ }$  $c_{\underline{a}}$ $\in$ $k$, and $c_{\underline{a}} \neq 0$ for only finitely many $\underline{a}$. $P_{j,i}$ is homogeneous with degree $\beta_{i}$ + $\beta_{j}$ so that, if $j$ = $i$ + $1$, we have $P_{j,i}$ = $0$. \\ $ $ \\ The reader will observe a little difference between those formulas and the original Levendorskii - Soibeman's one (\cite{LS}, prop. 5.5.2.), in which the left member of (\ref{LSb}) is $ $ $ $ $ X_{\beta_{j}}X_{\beta_{i}} - q^{(\beta_{i},\beta_{j})}X_{\beta_{i}}X_{\beta_{j}}. $ $ $ $ $ The reason is that Levendorskii and Soibelman use a version of the quantum group $\uqg$ which slightly differs from our's. Under our conventions, a direct proof of formulas (\ref{LSb}) and (\ref{LSPb}) is given in \cite{M}.
\item When $w$ = $w_{0}$, we have \mbox{ } $t$ = $N$, \mbox{ } $\{\beta_{1},  \mbox{ } ... , \mbox{ }  \beta_{t}\}$ = $\Phi^{+}$ \mbox{ } and \mbox{ } $\upw$ = $\upl$.
\end{itemizedot}
$ $\\
\bf An example \rm \\$ $\\

Assume, for sake of simplicity, that $ $ $k$ = \bbC $ $ is the complex numbers field, that $ $ $\g$ $ $ has type $ $ $A_{n}$ $ $ with $ $ $n$ $\geq$ $3$ and that the simple roots $ $ $\epsilon_{1}$, $ $ ... $ $, $ $ $\epsilon_{n}$ $ $ are ordered such as the Dynkin diagram is
$$
\begin{array}{ccccccccccccc}
\epsilon_{1}&-\!\!\!-\!\!\!-&\epsilon_{2}&-\!\!\!-\!\!\!-&\epsilon_{3}& -\!\!\!-\!\!\!-& &\cdots& &-\!\!\!-\!\!\!- &\epsilon_{n-1}&-\!\!\!-\!\!\!-&\epsilon_{n}
\end{array}
$$

Consider the following particular reduced expression of the longest element in $ $ $W$:
\begin{equation}\label{expw0}
w_{0} \mbox{ } = s_{\epsilon_{1}} \circ (s_{\epsilon_{2}} \circ s_{\epsilon_{1}}) \circ \mbox{ } ... \mbox{ } \circ (s_{\epsilon_{j}} \circ s_{\epsilon_{j-1}} \circ \mbox{ } ... \mbox{ } \circ s_{\epsilon_{1}}) \circ \mbox{ } ... \mbox{ } \circ (s_{\epsilon_{n}} \circ s_{\epsilon_{n-1}} \circ \mbox{ } ... \mbox{ } \circ s_{\epsilon_{1}})
\end{equation}

Denote by $ $ 
$Y_{1,2}$, $Y_{1,3}$, $Y_{2,3}$, $ $ ... $ $, $ $ $Y_{1,j+1}$, $Y_{2,j+1}$, $ $ ... $ $, $ $ $Y_{j,j+1}$, $ $ ... $ $, $ $ $Y_{1,n+1}$, $Y_{2,n+1}$, $ $ ... $ $, $ $ $Y_{n,n+1}$ $ $ the canonical generators of $ $ $\upl$ $ $ with respect to this reduced decomposition and observe that:

\newtheorem{lem2.1.}{Lemma 2.1.}

\begin{lem2.1.} $ $ \\
\begin{enumerate}
\item $Y_{j,j+1}$ = $E_{\alpha_{j}}$ $ $ for $ $ $1$ $\leq$ $j$ $\leq$ $n$. 
\item $Y_{i,j+1}$ = $Y_{i,j} Y_{j,j+1}$ $-$ $q^{-1} Y_{j,j+1} Y_{i,j}$ $ $ for $ $ $1$ $\leq$ $i$ $<$ $j$ $\leq$ $n$.
\end{enumerate}
\end{lem2.1.}

\bf Proof \rm \\

\begin{enumerate}
\item This results from the equality
$$s_{\epsilon_{1}} \circ (s_{\epsilon_{2}} \circ s_{\epsilon_{1}}) \circ \mbox{ } ... \mbox{ } \circ (s_{\epsilon_{j}} \circ s_{\epsilon_{j-1}} \circ \mbox{ } ... \mbox{ } \circ s_{\epsilon_{2}}) (\epsilon_{1}) \mbox{ } = \mbox{ } \epsilon_{j}.$$
\item Recall (see section 2.1) that
$$Y_{i,j+1}  \mbox{ } = \mbox{ } T_{\epsilon_{1}} \circ (T_{\epsilon_{2}} \circ T_{\epsilon_{1}}) \circ \mbox{ } ... \mbox{ } \circ (T_{\epsilon_{j-1}} \circ T_{\epsilon_{j-2}} \circ \mbox{ } ... \mbox{ } \circ T_{\epsilon_{1}}) \circ (T_{\epsilon_{j}} \circ T_{\epsilon_{j-1}} \circ \mbox{ } ... \mbox{ } \circ T_{\epsilon_{j+2-i}}) (E_{\epsilon_{j+1-i}}).$$
Since $ $ $\{\epsilon_{j-i}, \mbox{ } ... \mbox{ } , \mbox{ } \epsilon_{1}\}$ $\bot$ $\{\epsilon_{j}, \mbox{ } ... \mbox{ } , \mbox{ } \epsilon_{j+2-i}\}$ $ $ and $ $ $\{\epsilon_{j-i-1}, \mbox{ } ... \mbox{ } , \mbox{ } \epsilon_{1}\}$ $\bot$ $\epsilon_{j+1-i},$ $ $ $Y_{i,j+1}$ $ $ is equal to
$$T_{\epsilon_{1}} \circ (T_{\epsilon_{2}} \circ T_{\epsilon_{1}}) \circ \mbox{ } ... \mbox{ } \circ (T_{\epsilon_{j-1}} \circ T_{\epsilon_{j-2}} \circ \mbox{ } ... \mbox{ } \circ T_{\epsilon_{j+1-i}}) \circ (T_{\epsilon_{j}} \circ T_{\epsilon_{j-1}} \circ \mbox{ } ... \mbox{ } \circ T_{\epsilon_{j+2-i}}) \circ T_{\epsilon_{j-i}} (E_{\epsilon_{j+1-i}})$$
with, by \cite{J},
$$T_{\epsilon_{j-i}} (E_{\epsilon_{j+1-i}}) \mbox{ } = \mbox{ } E_{\epsilon_{j-i}} E_{\epsilon_{j+1-i}} \mbox{ } - \mbox{ } q^{-1} E_{\epsilon_{j+1-i}} E_{\epsilon_{j-i}}.$$
Moreover, we have
$$s_{\epsilon_{1}} \circ (s_{\epsilon_{2}} \circ s_{\epsilon_{1}}) \circ \mbox{ } ... \mbox{ } \circ (s_{\epsilon_{j-1}} \circ s_{\epsilon_{j-2}} \circ \mbox{ } ... \mbox{ } \circ s_{\epsilon_{j+1-i}}) \circ (s_{\epsilon_{j}} \circ s_{\epsilon_{j-1}} \circ \mbox{ } ... \mbox{ } \circ s_{\epsilon_{j+2-i}}) (\epsilon_{j+1-i})$$
$$= \mbox{ } s_{\epsilon_{1}} \circ (s_{\epsilon_{2}} \circ s_{\epsilon_{1}}) \circ \mbox{ } ... \mbox{ } \circ (s_{\epsilon_{j-1}} \circ s_{\epsilon_{j-2}} \circ \mbox{ } ... \mbox{ } \circ s_{\epsilon_{j+1-i}}) (\epsilon_{j+1-i} + \epsilon_{j+2-i} + \mbox{ } ... \mbox{ } + \epsilon_{j})$$
$$= \mbox{ } s_{\epsilon_{1}} \circ (s_{\epsilon_{2}} \circ s_{\epsilon_{1}}) \circ \mbox{ } ... \mbox{ } \circ (s_{\epsilon_{j-2}} \circ s_{\epsilon_{j-3}} \circ \mbox{ } ... \mbox{ } \circ s_{\epsilon_{1}}) (\epsilon_{j}) \mbox{ } = \mbox{ } \epsilon_{j}.$$
This implies that
$$T_{\epsilon_{1}} \circ (T_{\epsilon_{2}} \circ T_{\epsilon_{1}}) \circ \mbox{ } ... \mbox{ } \circ (T_{\epsilon_{j-1}} \circ T_{\epsilon_{j-2}} \circ \mbox{ } ... \mbox{ } \circ T_{\epsilon_{j+1-i}}) \circ (T_{\epsilon_{j}} \circ T_{\epsilon_{j-1}} \circ \mbox{ } ... \mbox{ } \circ T_{\epsilon_{j+2-i}}) (E_{\epsilon_{j+1-i}})$$ $$ = \mbox{ } E_{\epsilon_{j}} \mbox{ } = \mbox{ } Y_{j,j+1}.$$
Since $ $ $\{\epsilon_{j}, \mbox{ } \epsilon_{j-1} \mbox{ } ... \mbox{ } , \epsilon_{j+2-i}\}$ $\bot$ $\epsilon_{j-i},$ $ $ we have

$$T_{\epsilon_{1}} \circ (T_{\epsilon_{2}} \circ T_{\epsilon_{1}}) \circ \mbox{ } ... \mbox{ } \circ (T_{\epsilon_{j-1}} \circ T_{\epsilon_{j-2}} \circ \mbox{ } ... \mbox{ } \circ T_{\epsilon_{j+1-i}}) \circ (T_{\epsilon_{j}} \circ T_{\epsilon_{j-1}} \circ \mbox{ } ... \mbox{ } \circ T_{\epsilon_{j+2-i}}) (E_{\epsilon_{j-i}})$$
$$= \mbox{ } T_{\epsilon_{1}} \circ (T_{\epsilon_{2}} \circ T_{\epsilon_{1}}) \circ \mbox{ } ... \mbox{ } \circ (T_{\epsilon_{j-1}} \circ T_{\epsilon_{j-2}} \circ \mbox{ } ... \mbox{ } \circ T_{\epsilon_{j+1-i}}) (E_{\epsilon_{j-i}}) \mbox{ } = \mbox{ } Y_{i,j}.$$ \bx
\end{enumerate}

Denote by $ $ $v$ $ $ a square root of $ $ $q$ and, as in \cite{AD}, denote by $ $ $e_{i,j}$ $ $ the generators of $ $ $U^{+}$ $ $ constructed by H. Yamane \cite{Y}. We know that: \\ $ $ \\
\begin{itemizedot}
\item $e_{j,j+1}$ = $E_{\epsilon_{j}}$ $ $ for $ $ $1$ $\leq$ $j$ $\leq$ $n$.
\item $e_{i,j+1}$ = $v e_{i,j} e_{j,j+1}$ $-$ $v^{-1} e_{j,j+1} e_{i,j}$ $ $ for $ $ $1$ $\leq$ $i$ $<$ $j$ $\leq$ $n$.
\end{itemizedot}
This implies: \\
\begin{itemizedot}
\item $v^{i-j}e_{i,j+1}$ = $(v^{i-j+1} e_{i,j}) e_{j,j+1}$ $-$ $e_{j,j+1} q^{-1} (v^{i-j+1} e_{i,j})$ $ $ for $ $ $1$ $\leq$ $i$ $<$ $j$ $\leq$ $n$ 
and, from lemma 2.1. 1, we deduce (by induction on $ $ $j-i$):
\end{itemizedot}

\begin{lem2.1.} $ $ \\
$Y_{i,j+1}$ = $v^{i-j}e_{i,j+1}$ $ $ for $ $ $1$ $\leq$ $i$ $\leq$ $j$ $\leq$ $n$.
\end{lem2.1.}
$ $ \\
Now, consider an integer $ $ $p$ $ $ with $ $ $1$ $<$ $p$ $<$ $n$ $ $ and set
\begin{equation}\label{expwp}
w \mbox{ } = \mbox{ } (s_{\epsilon_{p}} \circ s_{\epsilon_{p-1}} \circ \mbox{ } ... \mbox{ } \circ s_{\epsilon_{1}}) \circ (s_{\epsilon_{p+1}} \circ s_{\epsilon_{p}} \circ \mbox{ } ... \mbox{ } \circ s_{\epsilon_{2}})  \circ \mbox{ } ... \mbox{ } \circ (s_{\epsilon_{n}} \circ s_{\epsilon_{n-1}} \circ \mbox{ } ... \mbox{ } \circ s_{\epsilon_{n-p+1}})
\end{equation}

$ $ \\

If $ $ $p$ $<$ $j$ $<$ $n$, $ $ we observe that \\ $\{\epsilon_{j-p}, \mbox{ } ... \mbox{ } , \mbox{ } \epsilon_{1}\}$ $\bot$ $\{\epsilon_{j+1}, \mbox{ } ... \mbox{ } , \mbox{ } \epsilon_{j+1-p+1}, \mbox{ } ... \mbox{ } , \mbox{ } \epsilon_{n}, \mbox{ } ... \mbox{ } , \mbox{ } \epsilon_{n-p+1}\}.$ \\

This implies that
$$w_{0} \mbox{ } = \mbox{ } w_{1} \circ w \circ w_{2}$$
with

$$w_{1} \mbox{ } = s_{\epsilon_{1}} \circ (s_{\epsilon_{2}} \circ s_{\epsilon_{1}}) \circ \mbox{ } ... \mbox{ } \circ (s_{\epsilon_{p-1}} \circ s_{\epsilon_{p-2}} \circ \mbox{ } ... \mbox{ } \circ s_{\epsilon_{1}}),$$

$$w_{2} \mbox{ } = s_{\epsilon_{1}} \circ (s_{\epsilon_{2}} \circ s_{\epsilon_{1}}) \circ \mbox{ } ... \mbox{ } \circ (s_{\epsilon_{n-p}} \circ s_{\epsilon_{n-p-1}} \circ \mbox{ } ... \mbox{ } \circ s_{\epsilon_{1}}).$$

If $ $ $d$, $ $ $d_{1}$, $ $ $d_{2}$ $ $ denote the number of simple reflections which appear respectively in the expressions of $ $ $w$, $ $ $w_{1}$, $ $ $w_{2}$ $ $ given above, we observe that $ $ $d$ + $d_{1}$ + $d_{2}$ $ $ is equal to the length of
$ $ $w_{0}$. $ $ This implies:

\begin{lem2.1.} $ $ \\
\begin{enumerate}
\item \rm(\ref{expwp})\it $ $ is a reduced expression of $ $ $w$.
\item If $ $ $p$ $\leq$ $j$ $\leq$ $n$, $ $ then
$$T_{w_{1}} \circ (T_{\epsilon_{p}} \circ T_{\epsilon_{p-1}} \circ \mbox{ } ... \mbox{ } \circ T_{\epsilon_{1}}) \circ \mbox{ } ... \mbox{ } \circ (T_{\epsilon_{j-1}} \circ T_{\epsilon_{j-2}} \circ \mbox{ } ... \mbox{ } \circ T_{\epsilon_{j-p+1}}) (E_{\epsilon_{j}}) \mbox{ } = \mbox{ } Y_{j,j+1}.$$
\item  If $ $ $p$ $\leq$ $j$ $\leq$ $n$ $ $ and  $ $ $1$ $\leq$ $i$ $\leq$ $p$, $ $ then
$$T_{w_{1}} \circ (T_{\epsilon_{p}} \circ T_{\epsilon_{p-1}} \circ \mbox{ } ... \mbox{ } \circ T_{\epsilon_{1}}) \circ \mbox{ } ... \mbox{ } \circ (T_{\epsilon_{j}} \circ T_{\epsilon_{j-1}} \circ \mbox{ } ... \mbox{ } \circ T_{\epsilon_{j-i+2}}) (E_{\epsilon_{j-i+1}}) \mbox{ } = \mbox{ } Y_{i,j+1}.$$
\end{enumerate}
\end{lem2.1.}

This implies that the canonical generators $ $ $X_{i,j+1}$ $ $ of $ $ $\upw$ $ $ verify
$$X_{i,j+1}  \mbox{ } = \mbox{ } T_{w_{1}}^{-1} (Y_{i,j+1})$$
for $ $ $p$ $\leq$ $j$ $\leq$ $n$ $ $ and  $ $ $1$ $\leq$ $i$ $\leq$ $p$. \\ $ $ \\

Since the Yamane generators $ $ $e_{i,j+1}$ $ $ ($p$ $\leq$ $j$ $\leq$ $n$ $ $ , $ $ $1$ $\leq$ $i$ $\leq$ $p$) $ $ verify the commutation relations of quantum matrices, it results from lemma 2.1. 2 that the generators $ $ $Y_{i,j+1}$ $ $ verify the same property. Now, as $ $ $T_{w_{1}}^{-1}$ $ $is an automorphism of $ $ $\uqg$ $ $ and as the ordered monomials in the variables $ $ $X_{i,j+1}$ $ $ (where the order on these variables is defined by the inverse lexicographic order on the indexes $ $ $(i,j+1)$) are a basis of  $\upw $, we conclude that

\newtheorem{prop2.1.}{Proposition 2.1.}

\begin{prop2.1.} $ $ \\
$\upw$ $ $ is the algebra $ $ $O_{q}(M_{p,m}(k))$ $ $ with $ $ $m$ = $n-p+1$ and with canonical generators the variables $ $ $X_{i,j+1}$ $ $ $(p$ $\leq$ $j$ $\leq$ $n$ $ $ , $ $ $1$ $\leq$ $i$ $\leq$ $p)$.
\end{prop2.1.}

\subsection{Diagrams.}
$w$, $t$, the simple roots $\alpha_{i}$, the positive roots $\beta_{i}$ ($1 \leq i \leq t$), are defined as in section \bf 2.1., \rm and we set
\begin{equation}\label{eqn:expv}
v \mbox{ } = \mbox{ } w^{-1} \mbox{ } = \mbox{ } s_{\alpha_{t}} \circ \mbox{ } ... \mbox{ } \circ s_{\alpha_{1}}
\end{equation}

\newtheorem{def22}{Definition 2.2.}

\begin{def22}
A diagram with respect to $($\ref{eqn:expw}$)$ is any subset $ $ $\Delta$ $ $ of $ $ $\llbracket  1,$ ... ,$t \rrbracket $. $ $ $($If there is no possible ambiguity, we omit to precise "with respect to $($\ref{eqn:expw}$)".)$
\end{def22}

\it In the following, we sometimes also omit the symbol $ $ $\circ$ $ $ in the composition of maps. \rm \\ $ $ \\

Consider a diagram $ $ $\Delta$. $ $ For any $i$ $\in$ $\llbracket  1,$ ... , $t \rrbracket $, we set $$s_{\alpha_{i}}^{\Delta} = \left\{
\begin{array}{cc}
s_{\alpha_{i}} & if \mbox{ } i \mbox{ } \in \mbox{ } \Delta \\
Id             & if \mbox{ } i \mbox{ } \notin \mbox{ } \Delta \\
\end{array} \right\}
$$ and we denote \\ $ $ \\
\begin{itemizedot}
\item $w^{\Delta} \mbox{ } =  \mbox{ } s_{\alpha_{1}}^{\Delta} \mbox{ } ... \mbox{ } s_{\alpha_{t}}^{\Delta}$
\item $v^{\Delta}$ = $\mbox{ } (w^{\Delta})^{-1}$ =  $\mbox{ } s_{\alpha_{t}}^{\Delta}$ ... $s_{\alpha_{1}}^{\Delta}$ 
\item for any $i$ $\in$ $\llbracket  1$, ... , $t$, $t+1 \rrbracket $, $ $ $w_{i}^{\Delta} \mbox{ } = \mbox{ } s_{\alpha_{i}}^{\Delta} \mbox{ } ... \mbox{ } s_{\alpha_{t}}^{\Delta}$ $ $ $ $ $ $ ($w_{1}^{\Delta} \mbox{ }$ = $w^{\Delta} \mbox{ }$, $w_{t+1}^{\Delta}$ = $Id$) 
\item for any $i$ $\in$ $\llbracket  0,$ $1$, ... , $t \rrbracket $, $ $ $v_{i}^{\Delta}$ = $\mbox{ } (w_{t-i+1}^{\Delta})^{-1}$  = $\mbox{ } s_{\alpha_{t}}^{\Delta} \mbox{ } ... \mbox{ }s_{\alpha_{t-i+2}}^{\Delta} s_{\alpha_{t-i+1}}^{\Delta}$ $ $ $ $ $ $ ($v_{0}^{\Delta}$ = $\mbox{ } (w_{t+1}^{\Delta})^{-1}$ = $Id$, $v_{t}^{\Delta}$ = $\mbox{ } v^{\Delta}$)
\end{itemizedot} $ $ \\ $ $ \\

\newtheorem*{rem}{Remark}
\begin{rem}[The case of example 2.1] $ $ \\ $ $ \\ 
Assume we are in the situation of example 2.1, so that
\[\begin{array}{rl}
w &= s_{\alpha_1}\circ s_{\alpha_2} \circ ... \circ s_{\alpha_{t}}\\
  &= (s_{\epsilon_{p}} \circ s_{\epsilon_{p-1}} \circ \mbox{ } ... \mbox{ } \circ s_{\epsilon_{1}}) \circ (s_{\epsilon_{p+1}} \circ s_{\epsilon_{p}} \circ \mbox{ } ... \mbox{ } \circ s_{\epsilon_{2}})  \circ \mbox{ } ... \mbox{ } \circ (s_{\epsilon_{n}} \circ s_{\epsilon_{n-1}} \circ \mbox{ } ... \mbox{ } \circ s_{\epsilon_{n-p+1}})
\end{array}\]
Set $ $ $m = n-p+1$, $ $ so that $ $ $t$ = $mp$, $ $ and consider a rectangular tableau consisting in $p \times m$ boxes labeled from $1$ to $mp$ as mentioned in the following figure.

\begin{center}
\begin{pgfpicture}{-1cm}{-1cm}{7cm}{7cm}%
\pgfsetroundjoin%
\pgfsetlinewidth{0.4pt}
\pgfmoveto{\pgfxy(0,0)}\pgflineto{\pgfxy(6,0)}\pgflineto{\pgfxy(6,6)}
\pgflineto{\pgfxy(0,6)}\pgfclosepath
\pgfstroke
\pgfxyline(1,0)(1,6)
\pgfxyline(2,0)(2,6)
\pgfxyline(5,0)(5,6)
\pgfxyline(0,1)(6,1)
\pgfxyline(0,5)(6,5)
\pgfsetlinewidth{0.2pt}
\pgfputat{\pgfxy(0.5,5.5)}{\pgftext{\color{black}\small $1$}}\pgfstroke
\pgfputat{\pgfxy(1.5,5.5)}{\pgftext{\color{black}\small $p+1$}}\pgfstroke
\pgfputat{\pgfxy(0.5,0.5)}{\pgftext{\color{black}\small $p$}}\pgfstroke
\pgfputat{\pgfxy(1.5,0.5)}{\pgftext{\color{black}\small $2p$}}\pgfstroke
\pgfputat{\pgfxy(5.5,0.5)}{\pgftext{\color{black}\small $mp$}}\pgfstroke
\pgfputat{\pgfxy(5.5,6.5)}{\pgftext{\color{black}\small $(m-1)p+1$}}\pgfstroke
\pgfxyline(5.5,6.3)(5.5,5.5)
\pgfsetmiterjoin \pgfsetfillcolor{black}
\pgfmoveto{\pgfxy(5.5667,5.6155)}\pgflineto{\pgfxy(5.5,5.5)}\pgflineto{\pgfxy(5.4333,5.6155)}\pgflineto{\pgfxy(5.5,5.5577)}\pgfclosepath\pgffillstroke
\pgfsetroundjoin%
\end{pgfpicture}
\end{center}
Following A. Postnikov \rm \cite{P}, \it let us draw this tableau with $p$
wires going along the rows and $m$ wires going along the columns.
Label ends of the wires from $1$ to $n+1$ as in the following
figure.
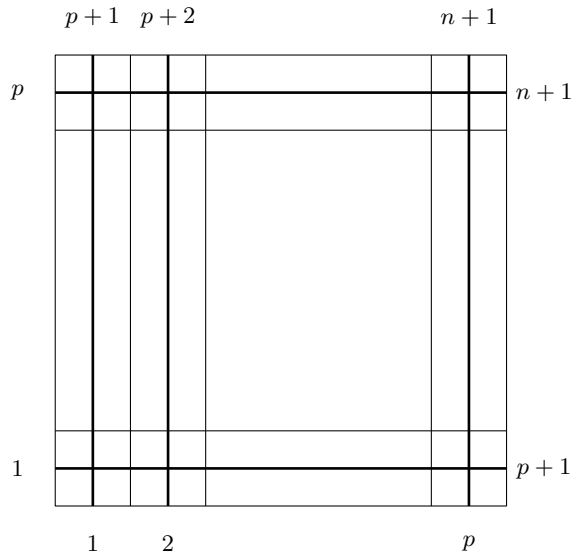
\begin{center}
\begin{pgfpicture}{-1.5cm}{-1.5cm}{7.5cm}{7.5cm}%
\pgfsetroundjoin%
\pgfsetlinewidth{0.4pt}
\pgfmoveto{\pgfxy(0,0)}\pgflineto{\pgfxy(6,0)}\pgflineto{\pgfxy(6,6)}
\pgflineto{\pgfxy(0,6)}\pgfclosepath
\pgfstroke
\pgfxyline(1,0)(1,6)
\pgfxyline(2,0)(2,6)
\pgfxyline(5,0)(5,6)
\pgfxyline(0,1)(6,1)
\pgfxyline(0,5)(6,5)
\pgfsetlinewidth{0.2pt}
\pgfputat{\pgfxy(-0.5,0.5)}{\pgftext{\color{black}\small $1$}}\pgfstroke
\pgfputat{\pgfxy(-0.5,5.5)}{\pgftext{\color{black}\small $p$}}\pgfstroke
\pgfputat{\pgfxy(0.5,6.5)}{\pgftext{\color{black}\small $p+1$}}\pgfstroke
\pgfputat{\pgfxy(1.5,6.5)}{\pgftext{\color{black}\small $p+2$}}\pgfstroke
\pgfputat{\pgfxy(5.5,6.5)}{\pgftext{\color{black}\small $n+1$}}\pgfstroke
\pgfputat{\pgfxy(6.5,5.5)}{\pgftext{\color{black}\small $n+1$}}\pgfstroke
\pgfputat{\pgfxy(6.5,0.5)}{\pgftext{\color{black}\small $p+1$}}\pgfstroke
\pgfputat{\pgfxy(5.5,-0.5)}{\pgftext{\color{black}\small $p$}}\pgfstroke
\pgfputat{\pgfxy(1.5,-0.5)}{\pgftext{\color{black}\small $2$}}\pgfstroke
\pgfputat{\pgfxy(0.5,-0.5)}{\pgftext{\color{black}\small $1$}}\pgfstroke
\pgfsetlinewidth{1pt}
\pgfxyline(0,0.5)(6,0.5)
\pgfxyline(0,5.5)(6,5.5)
\pgfxyline(0.5,6)(0.5,0)
\pgfxyline(1.5,6)(1.5,0)
\pgfxyline(5.5,6)(5.5,0)
\end{pgfpicture}
\end{center}
We observe that, if $W$ is identified as usual with the symmetric group $\mathcal{S}_{n+1}$ so that each $s_{\epsilon_i}$ is the transposition $(i ,  i+1)$, the wiring diagram defined over corresponds to the permutation $v$ = $w^{-1}$.\\
Now, consider any diagram $\Delta$ $($i.e. any subset of $\llbracket   1, mp \rrbracket)$, color in dark the boxes which label is in $\Delta$, and replace the corresponding crossing in non-colored boxes as a non crossing :\\
\begin{center}
\begin{tabular}{cc}
\begin{pgfpicture}{-1cm}{0cm}{2cm}{2cm}%
\pgfsetroundjoin%
\pgfsetlinewidth{0.4pt}
\pgfmoveto{\pgfxy(0,0)}\pgflineto{\pgfxy(1,0)}\pgflineto{\pgfxy(1,1)}
\pgflineto{\pgfxy(0,1)}\pgfclosepath
\pgfstroke
\pgfsetlinewidth{1pt}
\pgfmoveto{\pgfxy(0,0.5)}\pgfarc{90}{0}{0.5cm}
\pgfstroke
\pgfmoveto{\pgfxy(0.5,1)}\pgfarc{180}{270}{0.5cm}
\pgfstroke
\end{pgfpicture} &
\begin{pgfpicture}{-1cm}{0cm}{2cm}{2cm}%
\pgfsetroundjoin%
\pgfsetlinewidth{0.4pt}
\pgfmoveto{\pgfxy(0,0)}\pgflineto{\pgfxy(1,0)}\pgflineto{\pgfxy(1,1)}
\pgflineto{\pgfxy(0,1)}\pgfclosepath
\pgfstroke
\pgfsetlinewidth{1pt}
\pgfxyline(0.5,0)(0.5,1)
\pgfxyline(0,0.5)(1,0.5)
\end{pgfpicture}  \\
white box & colored box \\
\end{tabular}
\end{center}
Then, if $\Delta = \{i_1 < ... < i_l \}$, we observe that this new wiring diagram corresponds to the permutation $s_{\epsilon_{i_l}}\circ ... \circ s_{\epsilon_{i_1}} = v^{\Delta}$.\\
\end{rem}
So, using the \it A. Postnikov's conventions \rm (\cite{P}, p. 71-72) we can say:


\begin{enumerate}
\item $\Delta$ $ $ is a pipe dream  or a wiring diagram.
\item $u$ = $v^{\Delta}$ is the element of $W$ corresponding to the pipe dream $ $ $\Delta$ $ $ and $ $ $\Delta$ is a  pipe dream (or a wiring diagram) of $u$. \\ $ $ \\

In particular, we have: \\ 


\item $v$ is the element of $W$ corresponding to the full pipe dream $\Delta$ = $\llbracket  1$, ... , $t \rrbracket $.
\item $Id$ is the element of $W$ corresponding to the empty pipe dream $\Delta$ = $\emptyset$.
\end{enumerate}

Now, in the general case, let us denote by $\leq$ the Bruhat order on the Weyl group $W$. Since (\ref{eqn:expw}) is a reduced expression of $w$, we have that (\ref{eqn:expv}) is a reduced expression of $v$. \rm Let us recall (\cite{Jo}, corollary A.1.8) that an element $u$ of $W$ satisfies $u$ $\leq$ $w$ (resp. $u$ $\leq$ $v$) if and only if it can we written as a product of simple refections obtained by omitting some of the  $s_{\alpha_{i}}$ in the reduced expression (\ref{eqn:expw}) (resp. (\ref{eqn:expv})). So, we have immediately:

\newtheorem{lem22}{Lemma 2.2.}

\begin{lem22} $ $ \smallskip
\begin{enumerate}
\item The map: $ $ $ $ $\Delta$ $ $ $\mapsto$ $ $ $u$ = $w^{\Delta}$ $ $ $($resp. $u$ = $v^{\Delta})$ is surjective from the set of diagrams onto the set $\{u \in W \mbox{ } | \mbox{ } u \mbox{ } \leq \mbox{ } w \}$ $ $ $($resp. $\{u \in W \mbox{ } | \mbox{ } u \mbox{ } \leq \mbox{ } v \})$.
\item The map $ $ $ $ $u$ $ $ $\mapsto$ $ $ $u^{-1}$ is a bijection from $\{u \in W \mbox{ } | \mbox{ } u \mbox{ } \leq \mbox{ } w \}$ $ $ onto $\{u \in W \mbox{ } | \mbox{ } u \mbox{ } \leq \mbox{ } v \}$.
\end{enumerate}
\end{lem22}

\begin{lem22} $ $ \\ $ $ \smallskip
\bf v$^{\Delta}$ \rm $ $ := $ $  ($v_{0}^{\Delta}$, $v_{1}^{\Delta}$, ... , $v_{t}^{\Delta}$) $ $ $ $ is a subexpression of (\ref{eqn:expv}) in the sense of Marsh and Rietsch \cite{MR}.
\end{lem22}  \bf

\begin{flushleft}
Proof
\end{flushleft} \rm
This is because $ $ $v_{i-1}^{\Delta}$ = $s_{\alpha_{t}}^{\Delta} \mbox{ } ... \mbox{ } s_{\alpha_{t-i+2}}^{\Delta}$ $ $ $ $ $\Longrightarrow$ $ $ $ $ $v_{i}^{\Delta}$ = $v_{i-1}^{\Delta}s_{\alpha_{t-i+1}}^{\Delta}$ $ $ $ $ with $$s_{\alpha_{t-i+1}}^{\Delta} = \left\{
\begin{array}{cc}
s_{\alpha_{t-i+1}} & if \mbox{ } \mbox{ } t-i+1 \mbox{ } \in \mbox{ } \Delta \\
Id             & if \mbox{ } \mbox{ } t-i+1 \mbox{ } \notin \mbox{ } \Delta \\
\end{array} \right\}.
$$ \bx

\begin{lem22} $ $ \\ $ $ \smallskip
If $ $ \bf u \rm $ $ = $ $  ($u_{0}$, $u_{1}$, ... , $u_{t}$) \it is a subexpression of $v$, there exists a unique diagram $ $ $\Delta$ $ $ such that $ $ \bf v$^{\Delta}$ \rm $ $ = \bf u \rm. \it So, the map $ $ $ $ $\Delta$ $ $ $\mapsto$ $ $ \bf v$^{\Delta}$ \it is a bijection from the set of diagrams onto the set of subexpressions of $v$.
\end{lem22}  \bf

\begin{flushleft}
Proof
\end{flushleft} \rm
Since \bf u \rm $ $ is a subexpression of $v$, we have $ $ $u_{0}$ = $Id$ $ $ and
$$ (\forall i \in \mbox{ } \llbracket  1, ... ,t \rrbracket ) \mbox{ } \mbox{ } \mbox{ } (u_{i-1})^{-1}u_{i} \mbox{ } \in  \mbox{ } \{s_{t-i+1}, \mbox{ } Id\}.$$
Given any diagram $\Delta$, we have $ $ \bf v$^{\Delta}$ \rm = ($v_{0}^{\Delta}$, $v_{1}^{\Delta}$, ... , $v_{t}^{\Delta}$). $ $ Since $v_{0}^{\Delta}$ = $u_{0}$ = $Id$, we have  $ $ \bf v$^{\Delta}$ \rm $ $ = \bf u \rm $ $ if and only if
$$(\forall i \in \mbox{ } \llbracket  1, ... ,t \rrbracket ) \mbox{ }\mbox{ }\mbox{ }\mbox{ }(u_{i-1})^{-1}u_{i} \mbox{ } = \mbox{ } (v_{i-1}^{\Delta})^{-1}v_{i}^{\Delta} \mbox{ } = \mbox{ } s_{\alpha_{t-i+1}}^{\Delta}. $$
Now, there exists a unique diagram $\Delta$ which satisfies those conditions. It is defined by
$$ (\forall i \in \llbracket  1, ... , t \rrbracket ) \mbox{ } \mbox{ } \mbox{ } (t-i+1 \in \Delta) \mbox{ } \Leftrightarrow \mbox{ } (u_{i-1})^{-1}u_{i} \mbox{ } = \mbox{ } {s_{t-i+1}}. $$ \bx

\begin{def22}
We say that a diagram $\Delta$ is positive with respect to the reduced decomposition $(\ref{eqn:expw})$ $($or that $\Delta$ is positive if there is no ambiguity$)$ if $ $ \bf v$^{\Delta}$ \rm  =   $(v_{0}^{\Delta}$, $v_{1}^{\Delta}$, ... , $v_{t}^{\Delta})$ \it is a positive subexpression of $(\ref{eqn:expv})$ in the sense of Marsh and Rietsch \rm \cite{MR}, \it namely:
$$ v_{i-1}^{\Delta} \mbox{ } < \mbox{ } v_{i-1}^{\Delta}s_{\alpha_{t-i+1}} \mbox{ } \mbox{ } \mbox{ } \mbox{ } (\forall \mbox{ } i \mbox{ } \in \mbox{ } \llbracket  1, ... ,t \rrbracket ) $$
\end{def22} $ $

\newtheorem{obs22}{Observation 2.2.}

\begin{obs22}
The full diagram and the empty diagram are both positive.
\end{obs22}

\bf Proof \rm \\
\begin{itemizedot}
\item If $ $ $\Delta$ = $\llbracket  1,$ ... ,$t \rrbracket $ $ $ is the full diagram, we have $ $ $v_{i-1}^{\Delta}$ =  $\mbox{ } s_{\alpha_{t}} \mbox{ } ... \mbox{ }s_{\alpha_{t-i+2}}$ $ $ for any $ $ $i$ $\in$ $\llbracket  1,$ ... ,$t \rrbracket $ $ $ ($v_{0}^{\Delta}$ =  $Id$). $ $ Since (\ref{eqn:expv}) is a reduced expression of $ $ $v$, $ $  $\mbox{ } s_{\alpha_{t}} \mbox{ } ... \mbox{ }s_{\alpha_{t-i+2}} s_{\alpha_{t-i+1}}$ $ $ is a reduced expression of $ $ $v_{i-1}^{\Delta} s_{\alpha_{t-i+1}}$. $ $ So, clearly, we have $ $ $v_{i-1}^{\Delta}$ $<$ $v_{i-1}^{\Delta} s_{\alpha_{t-i+1}}$.
\item If $ $ $\Delta$ = $\emptyset$, $ $ we have $ $ $v_{i-1}^{\Delta}$ = $Id$, so that $ $ $v_{i-1}^{\Delta}$ $<$ $v_{i-1}^{\Delta} s_{\alpha_{t-i+1}}$.
\end{itemizedot}
\bx

Consider any $i \mbox{ } \in \mbox{ } \llbracket  1,  \mbox{ } ...  \mbox{ } ,t \rrbracket $, $ $ set $j$ = $t-i+1$ $ $ and $ $ $\{j_{1}$ $<$ ... $<$  $j_{s}\}$ $ $ =  $\Delta$ $\cap$ $ \mbox{ } \llbracket  j+1,  \mbox{ } ...  \mbox{ } ,t \rrbracket $ $ $ = $ $ $\Delta$ $\cap$ $\mbox{ } \llbracket  t-i+2,  \mbox{ } ...  \mbox{ } ,t \rrbracket $, $ $ $ $ so that $v_{i-1}^{\Delta}$ = $s_{\alpha_{t}}^{\Delta} \mbox{ } ... \mbox{ } s_{\alpha_{t-i+2}}^{\Delta}$ $ $ = $ $ $s_{\alpha_{j_{s}}} \mbox{ } ... \mbox{ } s_{\alpha_{j_{1}}}$ $ $ $ $ and $ $  $ $ $v_{i-1}^{\Delta}s_{\alpha_{t-i+1}}$ $ $ $ $  = $ $ $ $ $s_{\alpha_{j_{s}}} \mbox{ } ... \mbox{ } s_{\alpha_{j_{1}}}s_{\alpha_{j}}$. \\

\begin{flushleft}
If $w_{1}$ and $w_{2}$ are in $W$, it results from (\cite{Jo}, corollary A.1.8) that \\
\begin{tabular}{ll}
 $ $ $ $ $ $ ($w_{1}$ $<$ $w_{2}$)& $\Leftrightarrow$ $ $ $ $ ($w_{1}$ $\neq$ $w_{2}$ and $w_{1}$ is obtained by omitting some simple reflections in a reduced expression of $w_{2}$) $ $ $ $ $ $ \\ 
 &  $\Leftrightarrow$ $ $ $ $ ($w_{1}^{-1}$ $<$ $w_{2}^{-1}$).\\ 
\end{tabular} \\
So, with $w_{1}$ $ $ = $ $ $v_{i-1}^{\Delta}$ $ $ and $ $ $w_{2}$ $ $ = $ $ $v_{i-1}^{\Delta}s_{\alpha_{t-i+1}}$, $ $ we obtain: \end{flushleft}
$$ (v_{i-1}^{\Delta} \mbox{ } < \mbox{ } v_{i-1}^{\Delta}s_{\alpha_{t-i+1}}) \mbox{ } \mbox{ } \mbox{ } \mbox{ } \Leftrightarrow \mbox{ } \mbox{ } \mbox{ } \mbox{ } (s_{\alpha_{j_{1}}} \mbox{ } ... \mbox{ } s_{\alpha_{j_{s}}} \mbox{ } < \mbox{ } s_{\alpha_{j}}s_{\alpha_{j_{1}}} \mbox{ } ... \mbox{ } s_{\alpha_{j_{s}}})$$

\begin{flushleft}
Let us set $ $ $u$ $ $ = $ $ $s_{\alpha_{j_{1}}} \mbox{ } ... \mbox{ } s_{\alpha_{j_{s}}}$ $ $ and recall that $ $ $l(s_{\alpha_{j}}u)$ $ $ = $ $ $l(u)$ $ $ $\pm$ $ $ $1$. $ $ If $ $ $l(s_{\alpha_{j}}u)$ $ $ = $ $ $l(u)$ + $1$, $ $ we have that $ $ $u$ $ $ $<$ $ $ $s_{\alpha_{j}}u$ $ $ (\cite{Jo}, A.1.6). $ $  Conversely, if $ $ $u$ $ $ $<$ $ $ $s_{\alpha_{j}}u$, $ $ we have $ $ $l(u)$ $ $ $<$ $ $ $l(s_{\alpha_{j}}u)$ (\cite{Jo}, A.1.6) $ $  and, consequently, $ $  $l(s_{\alpha_{j}}u)$ $ $ = $ $ $l(u)$ + $1$. $ $ So, we have:
$$ (v_{i-1}^{\Delta} \mbox{ } < \mbox{ } v_{i-1}^{\Delta}s_{\alpha_{t-i+1}}) \mbox{ } \mbox{ } \mbox{ } \mbox{ } \Leftrightarrow \mbox{ } \mbox{ } \mbox{ } \mbox{ } l(s_{\alpha_{j}}s_{\alpha_{j_{1}}} \mbox{ } ... \mbox{ } s_{\alpha_{j_{s}}}) \mbox{ } = \mbox{ } 1 \mbox{ } + \mbox{ } l(s_{\alpha_{j_{1}}} \mbox{ } ... \mbox{ } s_{\alpha_{j_{s}}}) $$

Since the mapping $ $ $ $ $i$ $ $ $\mapsto$ $ $ $j$ $ $ = $ $ $t-i+1$ $ $ $ $ is a bijection of $ $  $\llbracket  1,  ...  , t \rrbracket $, $ $  we conclude:
\end{flushleft} \smallskip

\begin{lem22}

The diagram $\Delta$ is positive if and only if, for any $j \mbox{ } \in \mbox{ } \llbracket  1, ... ,t \rrbracket $,
$$ \Delta \mbox{ } \cap \mbox{ } \llbracket  j+1, ... ,t \rrbracket  = \{j_{1} < ... <  j_{s}\} \mbox{ } \Rightarrow \mbox{ } l(s_{\alpha_{j}}s_{\alpha_{j_{1}}} \mbox{ } ... \mbox{ } s_{\alpha_{j_{s}}}) = 1 \mbox{ }  + \mbox{ }  l(s_{\alpha_{j_{1}}} \mbox{ } ... \mbox{ } s_{\alpha_{j_{s}}}). $$

\end{lem22}

$ $ \\
Now, we can prove the following characterization of positive diagrams.

\newtheorem{prop22}{Proposition 2.2.}

\begin{prop22}

The diagram $\Delta$ is positive if and only if, for any $j \mbox{ } \in \mbox{ } \llbracket  1, ... ,t \rrbracket $,
$$ \Delta \mbox{ } \cap \mbox{ } \llbracket  j+1, ... ,t \rrbracket  = \{j_{1} < ... <  j_{s}\} \mbox{ } \Rightarrow \mbox{ } l(s_{\alpha_{j}}s_{\alpha_{j_{1}}} \mbox{ } ... \mbox{ } s_{\alpha_{j_{s}}}) = 1 \mbox{ }  + \mbox{ } s. $$

\end{prop22}

\begin{flushleft}
\bf Proof
\end{flushleft} \rm  Assume that $ $ $\Delta$ $ $ is positive, choose $ $ $j, \mbox{ } j_{1}, \mbox{ } ... \mbox{ } , \mbox{ } j_{s}\mbox{ }$  as above, and consider any $l$ with  $ $ $1 \mbox{ } < \mbox{ } l \mbox{ } < \mbox{ } s$. We have $\Delta \mbox{ } \cap \mbox{ } \llbracket  j_{l}+1, \mbox{ } ... \mbox{ } , \mbox{ } t \rrbracket $ $ $ = $ $ $\{j_{l+1} \mbox{ } < \mbox{ } ... \mbox{ } <  \mbox{ } j_{s}\}$ $ $ and, by lemma 2.2. 4, $ $ $l$($s_{\alpha_{j_{l}}}s_{\alpha_{j_{l+1}}} \mbox{ } ... \mbox{ } s_{\alpha_{j_{s}}}$) = $1$ + $l$($s_{\alpha_{j_{l+1}}} \mbox{ } ... \mbox{ } s_{\alpha_{j_{s}}}$). $ $ From this, we get that $ $ $l$($s_{\alpha_{j_{1}}} \mbox{ } ... \mbox{ } s_{\alpha_{j_{s}}})$ $ $ = $ $ $s$ $ $ and, again by lemma 2.2. 4, $ $ $l$($s_{\alpha_{j}}s_{\alpha_{j_{1}}} \mbox{ } ... \mbox{ } s_{\alpha_{j_{s}}})$ $ $ = $ $ $1$ + $s$. \\
Conversely, assume that this equality holds for any $j$ and choose $ $ $j, \mbox{ } j_{1}, \mbox{ } ... \mbox{ } , \mbox{ } j_{s}\mbox{ }$  as above. We have $ $ $\Delta \mbox{ }  \cap \mbox{ } \llbracket  j_{1}+1, \mbox{ } ... \mbox{ } , \mbox{ } t \rrbracket $ $ $ = $ $ $\{j_{2} \mbox{ } < \mbox{ } ... \mbox{ } <  \mbox{ } j_{s}\}$ $ $ and, by assumption, $ $ $l$($s_{\alpha_{j_{1}}}s_{\alpha_{j_{2}}} \mbox{ } ... \mbox{ } s_{\alpha_{j_{s}}}$) $ $ = $ $ $1$ $ $ + $ $ ($s$ $-$  $1$) $ $ = $ $ $s$. $ $ So, $ $ $l$($s_{\alpha_{j}}s_{\alpha_{j_{1}}} \mbox{ } ... \mbox{ } s_{\alpha_{j_{s}}}$) $ $ = $ $ $1$ $ $ + $ $ $l$($s_{\alpha_{j_{1}}} \mbox{ } ... \mbox{ } s_{\alpha_{j_{s}}}$) $ $ and $ $ $\Delta$ $ $ is positive by lemma 2.2. 4. \bx

\subsection{Some properties of positive diagrams.}

In this section, we use the same conventions as in section 2.2.

\newtheorem{prop23}{Proposition 2.3.}

\begin{prop23} $ $ \smallskip
\begin{enumerate}
\item The map: $ $ $ $ $\Delta$ $ $ $\mapsto$ $ $ $u$ = $v^{\Delta}$ $ $ $ $ is a bijection from the set of positive diagrams onto the set $ $ $ $ $\{u \in W \mbox{ } | \mbox{ } u \mbox{ } \leq \mbox{ } v \}$.
\item The map: $ $ $ $ $\Delta$ $ $ $\mapsto$ $ $ $u$ = $w^{\Delta}$ $ $ $ $ is a bijection from the set of positive diagrams onto the set $ $ $ $ $\{u \in W \mbox{ } | \mbox{ } u \mbox{ } \leq \mbox{ } w \}$.
\end{enumerate}
\end{prop23}


\bf Proof \rm
\begin{enumerate}
\item By lemma 2.2. 1, $f$: $ $ $\Delta$ $ $ $\mapsto$ $ $ $u$ = $v^{\Delta}$ $ $ $ $ is a map from the set of positive diagrams in the set $\{u \in W \mbox{ } | \mbox{ } u \mbox{ } \leq \mbox{ } v \}$.  Consider any $u$ $\in$ $W$ with $ $ $u$ $\leq \mbox{ } v $. By (\cite{MR}, Lemma 3.5.), there exists a unique positive subexpression \bf u$^{+}$ \rm = $ $ ($u_{0}$, $u_{1}$, ... , $u_{t}$) of (\ref{eqn:expv}) with $u_{t}$ = $u$. By lemma 2.2. 3, there exists a unique diagram $ $ $\Delta$ $ $ which satisfies  $ $ \bf v$^{\Delta}$ \rm  =  ($v_{0}^{\Delta}$, $v_{1}^{\Delta}$, ... , $v_{t}^{\Delta}$) =  \bf u$^{+}$ \rm. So, $ $ $\Delta$ $ $ is positive (def. 2.2. 2) and  $v^{\Delta}$ = $v_{t}^{\Delta}$ =  $u_{t}$ = $u$. This proves that $f$ is surjective. \\
If $ $ $\Delta^{\prime}$ $ $ is a positive diagram such that $ $ $v^{\Delta^{\prime}}$ $ $ = $u$, then $ $ \bf v$^{\Delta^{\prime}}$ \rm  =  ($v_{0}^{\Delta^{\prime}}$, $v_{1}^{\Delta^{\prime}}$, ... , $v_{t}^{\Delta^{\prime}}$) $ $ is a positive subexpression of  (\ref{eqn:expv}) with $v_{t}^{\Delta^{\prime}}$ = $v^{\Delta^{\prime}}$ = $u$. So, we have $ $ \bf v$^{\Delta^{\prime}}$ \rm  =   \bf u$^{+}$ \rm $ $ $\Rightarrow$  $ $ \bf v$^{\Delta^{\prime}}$ \rm  =   $ $ \bf v$^{\Delta}$ \rm and, by lemma 2.2. 3,  $ $ $\Delta^{\prime}$ = $\Delta$. This proves that $f$ is bijective.
\item Denote by $g$ the bijection $ $ $ $ $u$ $ $ $\mapsto$ $ $ $u^{-1}$ from $\{u \in W \mbox{ } | \mbox{ } u \mbox{ } \leq \mbox{ } w \}$ $ $ onto $\{u \in W \mbox{ } | \mbox{ } u \mbox{ } \leq \mbox{ } v \}$ $ $ (lemma 2.2. 1). \\ Then $g^{-1}$ $\circ$ $f$ : $ $ $\Delta$ $ $ $\mapsto$ $ $ $u$ = $(v^{\Delta})^{-1}$ = $w^{\Delta}$ $ $ $ $ is a bijection from the set of positive diagrams onto the set $\{u \in W \mbox{ } | \mbox{ } u \mbox{ } \leq \mbox{ } w \}$.
\end{enumerate}  \bx

Consider any $p$  $\in$ $\llbracket  1,$ ... , $t \rrbracket $, $p$  $\neq$ $1$, $t$, and set

\begin{equation}\label{eqn:expw1}
w_{1} = s_{\alpha_{1}} \mbox{ } ... \mbox{ } s_{\alpha_{p}}
\end{equation}
\begin{equation}\label{eqn:expw2}
w_{2} = s_{\alpha_{p+1}} \mbox{ } ... \mbox{ } s_{\alpha_{t}}
\end{equation} \\

We have $ $ $w$ = $w_{1}w_{2}$ $ $ and, since (\ref{eqn:expw}) is a reduced expression of $w$, (\ref{eqn:expw1}) and (\ref{eqn:expw2}) are reduced expressions of $w_{1}$ and $w_{2}$ respectively. \\ Denote by $\Delta$ any subset of $\llbracket  1,$ ... , $p \rrbracket $, so that $\Delta$ is a diagram with respect to (\ref{eqn:expw1}) and with respect to (\ref{eqn:expw}) both. \\ $ $ \\
\begin{prop23} $ $ \\
$(\Delta$ is positive with respect to $(\ref{eqn:expw1}))$ $ $ $ $ $ $ $\Leftrightarrow$ $ $ $ $ $ $ $(\Delta$ is positive with respect to $(\ref{eqn:expw}))$
\end{prop23}

\bf Proof \rm \\
Assume that $\Delta$ is positive with respect to (\ref{eqn:expw}) and consider any $j$  $\in$ $\llbracket  1,$ ... , $p \rrbracket $. We have $\Delta$ $ $ $\cap$ $\llbracket  j+1,$ ... , $p \rrbracket $ $ $ = $ $ $\Delta$ $ $ $\cap$ $\llbracket  j+1,$ ... , $t \rrbracket $ and, using the characterization of positive diagrams given in proposition 2.2. 1, we obtain that $\Delta$ is positive with respect to (\ref{eqn:expw1}). \\
Assume that $\Delta$ is positive with respect to (\ref{eqn:expw1}) and consider any $j$  $\in$ $\llbracket  1,$ ... , $t \rrbracket $.
If $j$ $\leq$ $p$, we have $\Delta$ $\cap$ $\llbracket  j+1,$ ... , $t \rrbracket $ $ $ = $ $ $\Delta$ $\cap$ $\llbracket  j+1,$ ... , $p \rrbracket $, so that the characterization given in proposition 2.2. 1 is satisfied. If $j$ $>$ $p$, we have $\Delta$ $\cap$ $\llbracket  j+1,$ ... , $t \rrbracket $ $ $ = $ $ $\emptyset$, so that the characterization given in proposition 2.2. 1 is again satisfied. \bx

Now, consider a non empty diagram $\Delta$ (with respect to (\ref{eqn:expw})) and two integers $j$, $m$ in $\llbracket  1,$ ... , $t \rrbracket $ with $j$ $<$ $m$ and $m$ $\in$ $\Delta$. Let us recall that ($\beta_{1}$, ... , $\beta_{t}$) is the sequence of positive roots associated to the reduced expression (\ref{eqn:expw})) of $w$ (section 2.1) and let us denote: \\

\begin{itemize}
\item $\overline{\Delta}$ = $\llbracket  1,$ ... , $t \rrbracket $ $\setminus$ $\Delta$,
\item $\Delta$ $\cap$ $\llbracket  j+1,$ ... , $m-1 \rrbracket $ = $\{j_{1}$ $<$ $...$ $<$ $j_{r}\}$ (unless this set is empty),
\item $\overline{\Delta}$ $\cap$ $\llbracket  j+1,$ $...$ , $m-1 \rrbracket $ = $\{l_{1}$ $<$ $...$ $<$ $l_{p}\}$ (unless this set is empty),
\item ($\gamma_{1}$, $...$ , $\gamma_{p}$, $\gamma_{p+1}$) is the sequence of (non necessarily positive) roots defined recursively by $ $ $\gamma_{p+1}$ = $\beta_{m}$ $ $ and, for \\ $1$ $\leq$ $i$ $\leq$ $p$, $\gamma_{i} = s_{\beta_{l_{i}}}(\gamma_{i+1})$,
\item For $1$ $\leq$ $i$ $\leq$ $p$, we set $a_{i} = (\beta_{l_{i}}^{\vee} , \gamma_{i+1})$.
\end{itemize}

Until the end of this section, we assume that $ $ $\overline{\Delta}$ $\cap$ $\llbracket  j+1,$ $...$ , $m-1 \rrbracket $ $ $ is nonempty, so that $ $ $p$ $ $ and the roots $ $ $\gamma_{i}$ $ $ ($1$ $\leq$ $i$ $\leq$ $p + 1$) $ $ are well defined.

\newtheorem{lem23}{Lemma 2.3.}

\begin{lem23}
$$\gamma_{1} = \beta_{m} - a_{p}\beta_{l_{p}} - \mbox{ } ... - \mbox{ } a_{1}\beta_{l_{1}}$$
\end{lem23}

\bf Proof \rm \\
For $1$ $\leq$ $i$ $\leq$ $p$, we have $ $ $\gamma_{i}$ = $\gamma_{i+1} \mbox{ } - \mbox{ } a_{i}\beta_{l_{i}}$. $ $ Summing these equalities, we get $ $ $\gamma_{1} = \gamma_{p+1} - a_{p}\beta_{l_{p}} - \mbox{ } ... \mbox{ } - a_{1}\beta_{l_{1}}$ $ $ and, by definition, $ $ $\gamma_{p+1}$ = $\beta_{m}$. \bx

Set
\begin{equation}\label{eqn:expwpr}
w^{\prime} = s_{\alpha_{1}} \mbox{ } ... \mbox{ } s_{\alpha_{m-1}}
\end{equation}
and observe that, since (\ref{eqn:expw}) is a reduced expression of $w$, (\ref{eqn:expwpr}) is a reduced expression of $w^{\prime}$.

\begin{lem23} $ $ \\
For any $i$ $\in$ $\llbracket  1,$ $...$ , $p \rrbracket $, $ $ denote by $ $ $w^{\prime}_{i}$ $ $ the element of $W$ obtained by omitting $ $ $s_{\alpha_{l_{i}}}$, $...$ , $s_{\alpha_{l_{p}}}$ $ $ in $(\ref{eqn:expwpr})$. \\ Then, we have $\gamma_{i}$ = $w^{\prime}_{i}(\alpha_{m})$.
\end{lem23}

\bf Proof \rm \\
Assume that $i$ = $p$ and set $u_{p}$ = $s_{\alpha_{1}} \mbox{ } ... \mbox{ } s_{\alpha_{l_{p}}-1}$, so that: \\ $w^{\prime}$ = $u_{p}s_{\alpha_{l_{p}}}s_{\alpha_{l_{p}+1}}  \mbox{ } ... \mbox{ } s_{\alpha_{m-1}}$ $ $ $\Rightarrow$ $ $ $u_{p}^{-1}w^{\prime}$ = $s_{\alpha_{l_{p}}}s_{\alpha_{l_{p}+1}}  \mbox{ } ... \mbox{ } s_{\alpha_{m-1}}$ $ $ $\Rightarrow$ $ $ $u_{p}s_{\alpha_{l_{p}}}u_{p}^{-1}w^{\prime}$ = $u_{p}s_{\alpha_{l_{p}+1}}  \mbox{ } ... \mbox{ } s_{\alpha_{m-1}}$ = $w^{\prime}_{p}$. \\ As $u_{p}(\alpha_{l_{p}})$ = $\beta_{l_{p}}$, we conclude that $w^{\prime}_{p}$ = $s_{\beta_{l_{p}}}w^{\prime}$. $ $ So, $ $ $w^{\prime}_{p}$($\alpha_{m}$) =  $s_{\beta_{l_{p}}}w^{\prime}$($\alpha_{m}$) = $s_{\beta_{l_{p}}}$($\beta_{m}$) = $\gamma_{p}$. \\
Assume $i$ $<$ $p$, $ $ $\gamma_{i+1}$ = $w^{\prime}_{i+1}$($\alpha_{m}$), $ $ and set $ $ $u_{i}$ = $s_{\alpha_{1}} \mbox{ } ... \mbox{ } s_{\alpha_{l_{i}}-1}$. We have $w^{\prime}_{i+1}$ = $u_{i}s_{\alpha_{l_{i}}}\sigma$ with $\sigma$ $\in$ $W$ and $w^{\prime}_{i}$ = $u_{i}\sigma$. Form this, we deduce (as in the case $i$ = $p$) that $u_{i}s_{\alpha_{l_{i}}}u_{i}^{-1}w^{\prime}_{i+1}$ = $w^{\prime}_{i}$. $ $ So, $ $ $w^{\prime}_{i}$($\alpha_{m}$) =  $s_{\beta_{l_{i}}}w^{\prime}_{i+1}$($\alpha_{m}$) = $s_{\beta_{l_{i}}}$($\gamma_{i+1}$) = $\gamma_{i}$. \bx \smallskip

Now, we can prove:

\begin{prop23} Assume that
$$\beta_{j} + \beta_{m} \mbox{ } = \mbox{ } a_{p}\beta_{l_{p}} + \mbox{ } ... + \mbox{ } a_{1}\beta_{l_{1}}.$$
Then $\Delta$ is not positive $($with respect to $(\ref{eqn:expw}))$.
\end{prop23}

\bf Proof \rm \\
By lemma 2.3. 1, we have $\beta_{j}$ = $-$ $\gamma_{1}$ and, by lemma 2.3. 2, we can write $ $ $\beta_{j}$ = $-$ $w^{\prime}_{1}$($\alpha_{m}$). Recall that $w^{\prime}_{1}$ is obtained by omitting $ $ $s_{\alpha_{l_{1}}}$, $ $ $...$ , $ $ $s_{\alpha_{l_{p}}}$ $ $ in the reduced expression (\ref{eqn:expwpr}) of $w^{\prime}$ $ $ and that \\ $\{l_{1}$ $<$ $...$ $<$ $l_{p}\}$ $ $ = $ $ $\overline{\Delta}$ $\cap$ $\llbracket  j+1,$ $...$ , $m-1 \rrbracket $ $ $ = $ $  $\llbracket  j+1,$ $...$ , $m-1 \rrbracket $ $ $ $\setminus$ $ $ $\{j_{1}$ $<$ $...$ $<$ $j_{r}\}$. $ $So, we have
$$w^{\prime}_{1} \mbox{ } = \mbox{ } s_{\alpha_{1}} \mbox{ } ... \mbox{ } s_{\alpha_{j-1}}s_{\alpha_{j}}s_{\alpha_{j_{1}}} \mbox{ } ... \mbox{ } s_{\alpha_{j_{r}}} \mbox{ } = \mbox{ } \eta s_{\alpha_{j}}s_{\alpha_{j_{1}}} \mbox{ } ... \mbox{ } s_{\alpha_{j_{r}}}$$
with $\eta \mbox{ } $ = $\mbox{ } s_{\alpha_{1}} \mbox{ } ... \mbox{ } s_{\alpha_{j-1}}$. \\
Since $ $ $\beta_{j}$ = $\eta$($\alpha_{j}$), $ $ we have $ $ $\beta_{j}$ = $-$ $w^{\prime}_{1}$($\alpha_{m}$) $ $ $\Rightarrow$ $ $ $\alpha_{j}$ = $-$ $s_{\alpha_{j}}s_{\alpha_{j_{1}}} \mbox{ } ... \mbox{ } s_{\alpha_{j_{r}}}$($\alpha_{m}$) $ $ $\Rightarrow$ $ $ $s_{\alpha_{j}}s_{\alpha_{j_{1}}} \mbox{ } ... \mbox{ } s_{\alpha_{j_{r}}}$($\alpha_{m}$) = $- \alpha_{j}$ $ $ is a negative root. By (\cite{H}, lemma 10.2.C), this implies that  $ $ $l$($s_{\alpha_{j}}s_{\alpha_{j_{1}}} \mbox{ } ... \mbox{ } s_{\alpha_{j_{r}}}s_{\alpha_{m}}$) $<$ $r+2$. \\
As $ $ $m$ $\in$ $\Delta$, $ $ we have  $ $ $\Delta$ $\cap$ $\llbracket  j+1,$ ... , $t \rrbracket $ = $\{j_{1}$ $<$ $...$ $<$ $j_{s}\}$ $ $ with $s$ $>$ $r$ and $j_{r+1}$ = $m$. So, $ $ $l$($s_{\alpha_{j}}s_{\alpha_{j_{1}}} \mbox{ } ... \mbox{ } s_{\alpha_{j_{s}}}$) $ $ $\leq$  $ $ $l$($s_{\alpha_{j}}s_{\alpha_{j_{1}}} \mbox{ } ... \mbox{ } s_{\alpha_{j_{r}}}s_{\alpha_{m}}$) $+$ $s-(r+1)$ $<$ $r+2+s-(r+1)$ = $s+1$ and, by proposition 2.2. 1, $\Delta$ is not positive. \bx

\section{Background on the deleting derivations algorithm.}

\subsection{Conventions.}

In this section, we use the conventions of section 2.1 and we set $R$ $ $ = $ $ $\upw$. In order to simplify a little bit the notations, we set \\
\begin{itemize}
\item $X_{i} \mbox{ } = \mbox{ } X_{\beta_{i}}$ for any $i \mbox{ } \in \mbox{ } \llbracket  1, \mbox{ } ... , \mbox{ } t \rrbracket $, so that $ $ $R$ = $k<X_{1},  \mbox{ } ... \mbox{ }, \mbox{ } X_{t}>$. $ $ Moreover, $ $ $X_{1},  \mbox{ } ... \mbox{ }, \mbox{ } X_{t}$ $ $ are called the canonical generators \rm (with respect to the reduced decomposition (\ref{eqn:expw})) of $ $ $R$.
\item Recall (section 2.1) that, for each $\rho$ in the root lattice \bbZ$\Pi$, the map $h_{\rho}$: $u$ $ $ $\mapsto$ $ $ $K_{\rho}.u$ is in $Aut(R)$, the group of automorphisms of the algebra $R$.
\item Let us set $H$ = $\{h_{\rho} \mbox{ } | \mbox{ } \rho \in$ \bbZ$\Pi\}$ and observe that $H$ is an abelian subgroup of $Aut(R)$.
\item Recall that $R$ is  \bbZ$\Pi$-graded and that, for any homogeneous element $u$ of degree $\gamma$ in $R$, for any  $\rho$ $\in$ \bbZ$\Pi$, we have $h_{\rho}$($u$) = $q^{-(\rho,\gamma)}u$.
\item For each $ $ ($i,j$) $\in$ ($\llbracket  1,$ ... , $t \rrbracket $)$^{2}$, $ $ set $ $ $\lambda_{i,j}$ $=$ $q^{-(\beta_{i},\beta_{j})}$, $ $ $q_{i}$ $=$ $\lambda_{i,i}$ $=$ $q^{-\|\beta_{i}\|^{2}}$, $ $ and observe that $q_{i}$ is not a root of unity.
\item If $1$ $\leq$ $i$ $<$ $j$ $\leq$ $t$, the Levendorskii-Soibeman formula can be written
\begin{equation}\label{LSR}
X_{j}X_{i} - \lambda_{j,i}X_{i}X_{j} \mbox{ } = \mbox{ } P_{j,i},
\end{equation} \\
\begin{equation}\label{LSP}
P_{j,i} \mbox{ } = \mbox{ } \sum_{ \underline{a} \mbox{ } = \mbox{ } (a_{i+1}, ... ,a_{j-1})} c_{\underline{a}}X_{i+1}^{a_{i+1}} ... X_{j-1}^{a_{j-1}}.
\end{equation} \\
with $ $ $\underline{a}$ $\in$ \bbN$^{j-i-1}$, $\mbox{ }$  $c_{\underline{a}}$ $\in$ $k$, $ $ and $c_{\underline{a}} \neq 0$ for only finitely many $\underline{a}$ . Moreover, $P_{j,i}$ is homogeneous with degree $\beta_{i}$ + $\beta_{j}$ $ $ so that, if $ $ $j$ = $i$ + $1$, $ $ we have $ $ $P_{j,i}$ = $0$. $ $ In the general case, this also implies that, if $ $ $c_{\underline{a}}$ $\neq$ $0$, $ $ then $ $ $a_{i+1} \beta_{i+1} \mbox{ } + \mbox{ }  ... \mbox{ } + \mbox{ } a_{j-1} \beta_{j-1}$ = $\beta_{i}$ + $\beta_{j}$.
\item Since $ $ $\beta_{i}$ $ $ and $ $ $\beta_{j}$ $ $ are positive roots, $ $ $\beta_{i}$ + $\beta_{j}$ $ $ is non zero. So, if  $ $ $c_{\underline{a}}$ $\neq$ $0$, $ $ then $ $ $\underline{a}$ $ $ is nonzero.
\item From this, we get that the generators $ $ $X_{1},  \mbox{ } ... \mbox{ }, \mbox{ } X_{t}$ $ $ of $R$ satisfy the equalities and the assumption 6.1.1. of (\cite{C}, section 6.1.).
\item Moreover, since the ordered monomials  $ $ $X^{\underline{a}}$ $ $ := $X_{1}^{a_{1}}$ ... $X_{t}^{a_{t}}$, $ $ $\underline{a}$ = ($a_{1}$, ... , $a_{1}$) $\in$ \bbN $^{t}$, are a basis of $R$, it results from  (\cite{C}, propositions 6.1.1. and 6.1.2.) that $R$ satisfies the conventions of (\cite{C}, section 3.1.). In particular, $R$ is an iterated Ore extension of the ground field $k$ and so, $R$ is a noetherian domain.
\item We denote by $ $ $F$ = $Fract(R)$ $ $ the division ring of fractions of $R$.
\item By (\cite{C}, proposition 6.1.1.), we also get that $R$ is the $k$-algebra generated by the "variables" $ $ $X_{1},  \mbox{ } ... \mbox{ }, \mbox{ } X_{t}$ $ $ submitted to the relations (\ref{LSR}).
\item For each $ $ $l$ $\in$ $\llbracket  1,$ ... , $t \rrbracket $, we set  $h_{l}$ = $h_{\beta_{l}}$  $\in$ $H$, and observe that, if $i$ $\in$ $\llbracket  1,$ ... , $t \rrbracket $, we have $h_{l}$($X_{i}$) = $\lambda_{l,i}$ $X_{i}$.
\item Since each $X_{i}$ is homogeneous, it is an $H$-eigenvector and, since $\lambda_{1,1}$ = $q^{-(\beta_{1},\beta_{1})}$ is not a root of unity, the assumption 4.1.2. of \cite{C} is satisfied. Since each $\lambda_{i,j}$ is a power of $q$, the assumption 4.1.1. of \cite{C} is also satisfied. As explained in \cite{C} (proof of lemma 4.2.2.) this implies that each prime ideal of $ $ $R$ $ $ is completely prime.
\item Recall that each automorphism $h$ $\in$ $H$ can be extended in a (unique) automorphism (denoted $h$ also) of $F$, so that $H$ can be seen as a subgroup of $Aut(F)$.
\end{itemize}

\subsection{The algebras $R^{(m)}$.}

Recall (\cite{C}, section 3.) that, for any $m$ $\in$ $\llbracket  2,$ ... , $t+1 \rrbracket $, there exists a family ($X_{1}^{(m)}$, $ $ ... $ $ , $ $ $X_{t}^{(m)}$) of new "variables" in $ $ $F$ $ $, called the canonical generators (with respect to the reduced decomposition $(1))$ of the algebra $R^{(m)}$ = $k<X_{1}^{(m)}, \mbox{ } ... \mbox{ } , \mbox{ } X_{t}^{(m)}>$, and which satisfies the following properties: \\
\begin{itemize}
\item If $ $ $1$ $\leq$ $i$ $<$ $j$ $\leq$ $t$, we have the following simplified Levendorskii-Soibelman formula: \\ $ $ \\
\begin{equation}\label{LSRm}
X_{j}^{(m)}X_{i}^{(m)} - \lambda_{j,i}X_{i}^{(m)}X_{j}^{(m)} \mbox{ } = \mbox{ } P_{j,i}^{(m)}
\end{equation}
\begin{center}
\it with
\end{center}
\begin{eqnarray}
&&\diamond \mbox{ } \mbox{ } \mbox{ } \mbox{ }  \mbox{ } \mbox{ } m \leq j \mbox{ } \Rightarrow \mbox{ } P_{j,i}^{(m)} \mbox{ } = \mbox{ } 0. \label{LSPm1}\\
&&\diamond \mbox{ } \mbox{ }  \mbox{ } \mbox{ } \mbox{ }  \mbox{ } j < m \mbox{ } \Rightarrow \mbox{ } P_{j,i}^{(m)} \mbox{ } = \mbox{ } \sum_{ \underline{a} \mbox{ } = \mbox{ } (a_{i+1}, ... ,a_{j-1})} c_{\underline{a}}(X_{i+1}^{(m)})^{a_{i+1}} ... (X_{j-1}^{(m)})^{a_{j-1}}  \label{LSPm2}
\end{eqnarray}
\begin{center}
\it where the coefficients $ $ $c_{\underline{a}}$ $ $ are the same as in $(\ref{LSP})$. \\
\rm($So,$ \it $ $ $P_{j,i}^{(m)} \mbox{ } = \mbox{ } 0$ $ $ in the case $ $ $i + 1$ = $j$ $<$ $m$ \rm $ $ and, in the general case, if $ $ $c_{\underline{a}}$ $\neq$ $0$, \\ we have $ $ $a_{i+1} \beta_{i+1}$ + $...$ + $a_{j-1} \beta_{j-1}$ = $\beta_{i}$ + $\beta_{i}$ $ $, so that $ $ $\underline{a}$ $ $ is nonzero.)
\end{center} \rm  $ $
\item This implies that the generators $ $ $X_{1}^{(m)},  \mbox{ } ... \mbox{ }, \mbox{ } X_{t}^{(m)}$ $ $ of $ $ $R^{(m)}$ $ $ still satisfy the equalities and the assumption 6.1.1. of (\cite{C}, section 6.1.).
\item The ordered monomials  $ $ $(X^{(m)})^{\underline{a}}$ $ $ := $(X_{1}^{(m)})^{a_{1}}$ ...  $(X_{t}^{(m)})^{a_{t}}$, $ $ $\underline{a}$ = ($a_{1}$, ... , $a_{1}$) $\in$ \bbN $^{t}$, are a basis of $R^{(m)}$. So, $R^{(m)}$ still satisfies the conventions of (\cite{C}, section 3.1.). In particular, $R^{(m)}$ is an iterated Ore extension of the ground field $k$ and so, $R^{(m)}$ is a noetherian domain.
\item We have $Fract$($R^{(m)}$) $ $ $=$ $ $ $Fract$($R$) $ $ $=$ $F$.
\item By (\cite{C}, proposition 6.1.1.), $R^{(m)}$ is the $k$-algebra generated by the "variables" $ $ $X_{1}^{(m)},  \mbox{ } ... \mbox{ }, \mbox{ } X_{t}^{(m)}$ $ $ submitted to the relations (\ref{LSRm}).
\item For each $\rho$ in the root lattice \bbZ$\Pi$ and for each $i$ $\in$ $\llbracket  1,$ ... , $t \rrbracket $, we still have$ $ $h_{\rho}$($X_{i}^{(m)}$) = $q^{-(\rho,\beta_{i})}X_{i}^{(m)}$. $ $ This implies that $h_{\rho}$($R^{(m)}$) = $R^{(m)}$. $ $ So $H$ can still be seen as a subgroup of $Aut(R^{(m)})$ and each $X_{i}^{(m)}$ is a $H$ eigenvector. If $l$ and $i$ are in $\llbracket  1,$ ... , $t \rrbracket $, we have $h_{l}$($X_{i}^{(m)}$) = $\lambda_{l,i}$ $X_{i}^{(m)}$ and, as above, this implies that $R^{(m)}$ satisfies the assumptions 4.4.2. and 4.4.1. of \cite{C}. As explained in (\cite{C}, proof of lemma 4.2.2), this implies that each prime ideal of $ $ $R^{(m)}$ $ $ is completely prime.
\item If $ $ $u$ $\in$ $F$, $ $ and if $ $ $\gamma$ $ $ $\in$ $ $ \bbZ$\Pi$, $ $ we say that $ $ $u$  $ $ \it is homogeneous of degree $ $ $\gamma$, \rm $ $ if $ $  $h_{\rho}$($u$) = $q^{-(\rho,\gamma)}u$ $ $ for all $ $ $\rho$ $ $ in $ $ \bbZ$\Pi$. $ $ So, for example, each $ $ $X_{i}^{(m)}$ $ $ is homogeneous of degree $ $ $\beta_{i}$.
\item If $ $ $u_{1}, \mbox{ } ... \mbox{ } , \mbox{ } u_{r}$ $ $ are homogeneous of same degree $ $ $\gamma$, $ $ then any linear combination of $ $ $u_{1}, \mbox{ } ... \mbox{ } , \mbox{ } u_{r}$ $ $ (with coefficients in $ $ $k$) $ $ is homogeneous of degree $ $ $\gamma$.
\item Clearly, if $ $ $u_{1}$  $ $ is homogeneous of degree $ $ $\gamma_{1}$ $ $ and $ $ $u_{2}$  $ $  is homogeneous of degree $ $ $\gamma_{2}$, $ $ then $ $ $u_{1} u_{2}$  $ $ is homogeneous of degree $ $ $\gamma_{1}$ +  $\gamma_{2}$.
\item Likewise, if $ $ $u$ $ $ is nonzero and homogeneous of degree $ $ $\gamma$, $ $ then $ $ $u^{-1}$ $ $ is homogeneous of degree $ $ $-$ $\gamma$.
\item So, if $ $ $u_{1}, \mbox{ } ... \mbox{ } , \mbox{ } u_{r}$ $ $ are nonzero and homogeneous of degrees $ $ $\gamma_{1}, \mbox{ } ... \mbox{ } , \mbox{ } \gamma_{r}$ $ $ respectively, if $ $ $\underline{a}$ = $(a_{1}, \mbox{ } ... \mbox{ } , \mbox{ } a_{r})$ $\in$ $ $ \bbZ$^{r}$,$ $ then $ $ $u^{\underline{a}}$ := $u_{1}^{a_{1}} \mbox{ } ... \mbox{ } u_{r}^{a_{r}}$ $ $ is homogeneous of degrees $ $ $a_{1} \gamma_{1} + \mbox{ } ... \mbox{ } + a_{r} \gamma_{r}$.
\item For each $j$ $\in$ $\llbracket  1,$ ... , $t \rrbracket $, denote by $ $ $\delta_{j}^{(m)}$ $ $ the left $h_{j}$ - derivation of $R^{(m)}$ defined by $X_{j}^{(m)}$. This means that, for each $a$ $\in$ $R^{(m)}$, $ $ $\delta_{j}^{(m)}$($a$) := $X_{j}^{(m)}a$ $-$ $h_{j}(a)X_{j}^{(m)}$.
\item If $j$ $\in$ $\llbracket  1,$ ... , $t \rrbracket $, we have
$$h_{j}(k<X_{1}^{(m)},  \mbox{ } ... \mbox{ }, \mbox{ } X_{j-1}^{(m)}>)  \mbox{ } =  \mbox{ } k<X_{1}^{(m)},  \mbox{ } ... \mbox{ }, \mbox{ } X_{j-1}^{(m)}>,$$
$$\delta_{j}^{(m)}(k<X_{1}^{(m)},  \mbox{ } ... \mbox{ }, \mbox{ } X_{j-1}^{(m)}>)  \mbox{ } \subset  \mbox{ } k<X_{1}^{(m)},  \mbox{ } ... \mbox{ }, \mbox{ } X_{j-1}^{(m)}>.$$
\end{itemize}
So, $h_{j}$ induces an automorphism, still denoted $h_{j}$, of the algebra $ $ $k<X_{1}^{(m)},  \mbox{ } ... \mbox{ }, \mbox{ } X_{j-1}^{(m)}>$. $ $ It is the only automorphism of $ $ $k<X_{1}^{(m)},  \mbox{ } ... \mbox{ }, \mbox{ } X_{j-1}^{(m)}>$ $ $ which satisfies $ $ $h_{j} (X_{i}^{(m)})$ =  $\lambda_{j,i} X_{i}^{(m)}$ $ $ for any $ $ $i$.\\
Likewise, $ $ $\delta_{j}^{(m)}$ $ $ induces a left $h_{j}$ - derivation of $ $ $k<X_{1}^{(m)},  \mbox{ } ... \mbox{ }, \mbox{ } X_{j-1}^{(m)}>$, $ $ still denoted $ $ $\delta_{j}^{(m)}$. $ $ It satisfies the following properties:
\begin{enumerate}
\item If $1$ $\leq$ $i$ $<$ $j$, we have $ $ $\delta_{j}^{(m)}$($X_{i}^{(m)}$) = $P_{j,i}^{(m)}$. \\ This implies that $ $ $\delta_{j}^{(m)}$ $ $ is zero on $ $ $k<X_{1}^{(m)},  \mbox{ } ... \mbox{ }, \mbox{ } X_{j-1}^{(m)}>$ $ $  and so, that $ $ $X_{j}^{(m)} a$ = $h_{j}(a) X_{j}^{(m)}$ $ $ for any $ $ $a$ $ $ in $ $ $k<X_{1}^{(m)},  \mbox{ } ... \mbox{ }, \mbox{ } X_{j-1}^{(m)}>$, $ $ as soon as $j$ $\geq$ $m$. \\
\item $\delta_{j}^{(m)}$ $ $ is locally nilpotent on $ $ $k<X_{1}^{(m)},  \mbox{ } ... \mbox{ }, \mbox{ } X_{j-1}^{(m)}>$.
\item We have
$$R^{(m)} =  k[X_{1}^{(m)}][X_{2}^{(m)};  h_{2}, \delta_{2}^{(m)}]  \mbox{ } ... \mbox{ } [X_{m-1}^{(m)}; h_{m-1}, \delta_{m-1}^{(m)}] \mbox{ } [X_{m}^{(m)}; h_{m}] \mbox{ } ... \mbox{ } [X_{t}^{(m)}; h_{t}].$$
\end{enumerate}
$ $ $\bullet$ $ $ $ $ For each $i$ $\in$ $\llbracket  1,$ ... , $t \rrbracket $, $ $ we have $ $ $ $ $X_{i}^{(t+1)}$ = $X_{i}$, $ $ so that $ $ $ $ $R^{(t+1)}$ = $R$. \\ $ $ \\

Until the end of this section, we assume that $ $ $m$ $\in$ $\llbracket  2,$ ... , $t \rrbracket $. \\  \smallskip

\begin{itemize}
\item Recall that $ $ $q_{m}$ = $\lambda_{m,m}$ $ $ is not a root of unity and define quantum integers $ $ $[l]_{q_{m}}$ $ $ (resp. quantum factorial $ $ $[l]!_{q_{m}}$) $ $ for $ $ $l$ $\in$ \bbN, $ $ as in (\cite{C}, section 2.).
\item For each $i$ $\in$ $\llbracket  1,$ ... , $t \rrbracket $, we have
\begin{enumerate}
\item $ $ $ $ $ $ $ $ $m \leq i \mbox{ } \Rightarrow \mbox{ } X_{i}^{(m)}  = \mbox{ } X_{i}^{(m+1)}$. $ $ In particular, we have $ $ $X_{m}^{(m)} = \mbox{ } X_{m}^{(m+1)}$. \\
\item $$i < m  \mbox{ } \Rightarrow \mbox{ } X_{i}^{(m)} = \mbox{ } X_{i}^{(m+1)} \mbox{ } \mbox{ } +  \mbox{ } \mbox{ } \displaystyle \sum_{l = 1}^{+\infty} C_{l}^{(m+1)}(X_{m}^{(m+1)})^{-l}$$
with
$$C_{l}^{(m+1)} = \mbox{ }  \displaystyle \frac{(1-q_{m})^{-l}}{[l]!_{q_{m}}} \mbox{ } \lambda_{m,i}^{-l}  \mbox{ } (\delta_{m}^{(m+1)})^{l} \mbox{ } (X_{i}^{(m+1)}).$$
\end{enumerate}
\item We observe that, since $ $ $k<X_{1}^{(m+1)},  \mbox{ } ... \mbox{ }, \mbox{ } X_{m-1}^{(m+1)}>$ $ $ is $ $ $\delta_{m}^{(m+1)}$ - stable, we have $ $ $C_{l}^{(m+1)}$ $\in$ $ $ $k<X_{1}^{(m+1)},  \mbox{ } ... \mbox{ }, \mbox{ } X_{m-1}^{(m+1)}>$ $ $ and, since $ $ $\delta_{m}^{(m+1)}$ $ $ is locally nilpotent on  $ $ $k<X_{1}^{(m+1)},  \mbox{ } ... \mbox{ }, \mbox{ } X_{m-1}^{(m+1)}>$, $ $ only finitely many  $ $ $C_{l}^{(m+1)}$ $ $ are non zero.
\item If $ $ $i < m$ $ $ and $ $ $X_{m}^{(m+1)}X_{i}^{(m+1)} - \lambda_{m,j}X_{i}^{(m+1)}X_{m}^{(m+1)}  = \mbox{ } P_{m,i}^{(m+1)} = \mbox{ } 0$, $ $ then $ $ $\delta_{m}^{(m+1)}$($X_{i}^{(m+1)}$) $ = \mbox{ } 0$ $ $ and so, $ $ $X_{i}^{(m)} = \mbox{ } X_{i}^{(m+1)}$.
\item There exists a unique homomorphism of $k$-algebras
$$ \Theta^{(m)}: \mbox{ } k<X_{1}^{(m+1)}, \mbox{ } ... \mbox{ }, \mbox{ } X_{m-1}^{(m+1)}> \mbox{ } \mbox{ } \rightarrow \mbox{ } \mbox{ } k<X_{1}^{(m)}, \mbox{ } ... \mbox{ }, \mbox{ } X_{m-1}^{(m)}> $$
which verifies $ $ $ $ $ $ $\Theta^{(m)}$($X_{i}^{(m+1)})$ $ $ = $ $ $X_{i}^{(m)}$ $ $ $ $ $ $ for $1$ $\leq$ $i$ $\leq$ $m-1$.
\item When $m$ = $2$, we set
\begin{enumerate}
\item $ $ $ $ $ $ $ $ $Z_{i}$ = $X_{i}^{(2)}$ for $i$ $\in$ $\llbracket  1,$ ... , $t \rrbracket $.
\item $ $ $ $ $ $ $ $ $\overline{R}$ = $R^{(2)}$.
\end{enumerate}
\item $\overline{R}$ $ $ = $ $ $k<Z_{1}, \mbox{ } ... \mbox{ }, \mbox{ } Z_{t}>$ is the $k$ - algebra generated by the variables  $Z_{1}, \mbox{ } ... \mbox{ }, \mbox{ } Z_{t}$ submitted to the relations:
\begin{equation}\label{LSRb}
(1 \leq i < j \leq t) \mbox{ } \Rightarrow \mbox{ } Z_{j}Z_{i} \mbox{ } = \mbox{ } \lambda_{j,i}Z_{i}Z_{j}
\end{equation}
\item For any $i$ $\in$ $\llbracket  1,$ ... , $t \rrbracket $, we have $X_{i}^{(2)}$ = $ $ $...$ $=$ $X_{i}^{(i+1)}$ = $Z_{i}$.
\item $S_{m}$ = $\{(X_{m}^{(m+1)})^{l} \mbox{ } | \mbox{ } l \in$ \bbN $\}$ $ $ = $ $ $\{(X_{m}^{(m)})^{l} \mbox{ } | \mbox{ } l \in \mbox{ } $ \bbN $\}$ $ $  is a multiplicative system of regular elements in $R^{(m)}$ and in $R^{(m+1)}$. It satisfies the Ore condition (on both sides) in each one of those rings and we have $ $ $R^{(m)}S_{m}^{-1}$ $ $ = $ $ $R^{(m+1)}S_{m}^{-1}$.
\end{itemize}
 
\subsection{Prime and $H$-prime spectrum of $ $ $R^{(m)}$, admissible diagrams}

Consider any integer $ $ $m$ $\in$ $\llbracket  2,$ ... , $t+1 \rrbracket $. $ $ As usual, we denote by $ $ $Spec(R^{(m)})$ $ $ (resp. $ $ $H-Spec(R^{(m)})$) $ $ the set of all prime ideals (resp. $H$- invariant prime ideals) of $ $ $R^{(m)}$. $ $ Recall (see section 3.2) that each prime ideal of $ $ $R^{(m)}$ $ $ is completely prime. If $ $ $m$ $\in$ $\llbracket  2,$ ... , $t \rrbracket $, $ $ we denote by $$\phi_{m}: \mbox{ }  Spec(R^{(m+1)}) \mbox{ } \mbox{ } \rightarrow \mbox{ } \mbox{ } Spec(R^{(m)})$$ the canonical injection defined in (\cite{C}, section 4.3). \\
Moreover, if $ $ $m^{\prime}$ $\leq$ $m + 1$ and if $ $ $\mathcal{P}$ $\in$ $Spec(R^{(m+1)})$, $ $ then $ $ $\mathcal{P}^{\prime}$ = $\phi_{m^{\prime}}$ $\circ$ $ $ ... $ $ $\circ$ $\phi_{m}$($\mathcal{P}$) is called the canonical image of $ $ $\mathcal{P}$ $ $ in $ $ $Spec(R^{(m^{\prime})})$. $ $ (If $ $ $m^{\prime}$ = $m+1$, $ $ this canonical image is $ $ $\mathcal{P}$ $ $ itself.) \\ $ $ \\

Let us now recall the main properties of the maps $ $ $\phi_{m}$. \\  $ $ \\

\begin{itemize}
\item Denote by $ $ $\mathcal{P}^{(m+1)}$ $ $ any prime ideal in $ $ $Spec(R^{(m+1)})$ $ $ and set $ $ $\mathcal{P}$$^{(m)}$ $ $ = $ $ $\phi_{m}(\mathcal{P}$$^{(m+1)})$.
\begin{enumerate}
\item If $X_{m}^{(m+1)}$ $ $ $\notin$ $ $ $\mathcal{P}$$^{(m+1)}$, $ $ then $\mathcal{P}$$^{(m+1)}$ $\cap$ $S_{m}$ $ $ = $ $ $\emptyset$ $ $ and  $$\mathcal{P}\it^{(m)} \mbox{ } = \mbox{ } R^{(m)} \cap \mathcal{P}\it^{(m+1)}S_{m}^{-1}.$$
In this case, $ $ $X_{m}^{(m)}$ $=$ $X_{m}^{(m+1)}$ $\notin$ $\mathcal{P}\it^{(m)}$.
\item If $X_{m}^{(m+1)}$ $ $ $\in$ $ $ $\mathcal{P}$$^{(m+1)}$, $ $ there exists a unique algebra homomorphism $ $ $ $ $g:  \mbox{ } R^{(m)} \mbox{ } \mbox{ } \rightarrow \mbox{ } \mbox{ } R^{(m+1)}/\mathcal{P}\it^{(m+1)}$ $ $ $ $ which verifies $ $ $g$($X_{i}^{(m)}$) $ $ = $ $ $X_{i}^{(m+1)}$ $ $ $+$ $ $ $\mathcal{P}\it^{(m+1)}$ $ $ for all $ $ $i$ $\in$ $\llbracket  1,$ ... , $t \rrbracket $, $ $ and $ $ $ $ $\mathcal{P}\it^{(m)}  \mbox{ } = \mbox{ } Ker(g).$ \\
In this case, $ $ $X_{m}^{(m)}$ $=$ $X_{m}^{(m+1)}$ $\in$ $\mathcal{P}\it^{(m)}$.
\item For any $ $ $h$ $\in$ $H$, $ $ we have  $ $ $h(\mathcal{P}$$^{(m)})$ $ $ = $ $ $\phi_{m}(h(\mathcal{P}\it^{(m+1)}))$ $ $ (\cite{C}, lemma 5.5.5.) so that, by the injectivity of $ $ $\phi_{m}$, $ $ $\mathcal{P}\it^{(m+1)}$ $ $ is $H$ - invariant if and only if $ $ $\mathcal{P}\it^{(m)}$ $ $ is $H$ - invariant.
\end{enumerate}
\item A diagram $ $ $\Delta$ $ $ is said admissible (with respect to the reduced decomposition (\ref{eqn:expw})) if there exists a prime ideal $ $ $\mathcal{P}$ $ $ of $ $ $R$ ( = $R^{(t+1)}$) whose canonical image $ $ $\mathcal{P}^{(2)}$ $ $ in $ $ $Spec(\overline{R})$ ( = $Spec(R^{(2)})$) verifies:

$$\mathcal{P}^{(2)} \mbox{ } \cap \mbox{ } \{Z_{i} \mbox{ } | 1 \mbox{ } \leq \mbox{ }  i \mbox{ } \leq \mbox{ } t\} \mbox{ } = \mbox{ } \{Z_{i} \mbox{ } | \mbox{ } i \mbox{ } \in \mbox{ } \Delta\}.$$

\item Consider a diagram $ $ $\Delta$ $ $ and denote by $ $ $\mathcal{P}^{(2)}_{\Delta}$ $ $ = $ $ ($\{Z_{i} \mbox{ } | \mbox{ } i \mbox{ } \in \mbox{ } \Delta\}$) $ $ the ideal generated by the variables $ $ $Z_{i}$ $ $ with $ $ $i \mbox{ } \in \mbox{ } \Delta$. $ $ Then, by (\cite{C}, proposition 5.5.1.), we have
\begin{enumerate}
\item $\mathcal{P}^{(2)}_{\Delta}$ $ $ $\in$ $ $ $H-Spec(R^{(2)})$.
\item Conversely, for any $ $ $\cal Q$ $ $ $\in$ $ $ $H-Spec(R^{(2)})$, there exists a (unique) diagram $ $ $\Delta$ $ $ such that $ $ $\cal Q$ $ $ = $ $ $\mathcal{P}^{(2)}_{\Delta}$, namely $\Delta = \{i \in \llbracket 1, t \rrbracket \ | \ Z_i \in Q\}$.
\end{enumerate}
\item A diagram $\Delta$ is admissible if and only if there exists $ $ $\mathcal{P}\it_{\Delta}$ $ $ $\in$ $ $ $Spec(R)$ $ $ (= $ $ $Spec(R^{(t+1)})$) $ $ such that $$\mathcal{P}^{(2)}_{\Delta} \mbox{ } = \mbox{ } \phi_{2} \circ \mbox{ } ... \mbox{ } \circ \phi_{t}(\mathcal{P}\it_{\Delta}).$$ (See \cite{C}, theorem 5.5.1. and observe that, since each $ $ $\phi_{i}$ $ $ is injective, $ $ $\mathcal{P}\it_{\Delta}$ $ $ is unique.)
\item The map $ $ $\Delta$ $ $ $\mapsto$ $ $ $\mathcal{P}\it_{\Delta}$ $ $ is a bijection from the set of admissible diagrams onto the set $ $ $H-Spec(R)$ $ $ (= $ $ $H-Spec(R^{(t+1)})$). In fact, if $ $ $\Delta$ $ $ is admissible, then $ $ $\mathcal{P}\it_{\Delta}$ $ $ is $H$ - invariant because $ $ $\mathcal{P}^{(2)}_{\Delta} \mbox{ } = \mbox{ } \phi_{2} \circ \mbox{ } ... \mbox{ } \circ \phi_{t}(\mathcal{P}\it_{\Delta})$ $ $ is. So, $ $ $\Delta$ $ $ $\mapsto$ $ $ $\mathcal{P}\it_{\Delta}$ $ $ is a map from the set of admissible diagrams into $ $ $H-Spec(R)$. $ $ It is injective because the  map $ $ $\Delta$ $ $ $\mapsto$ $ $ $\mathcal{P}^{(2)}_{\Delta}$ $ $ is injective. If $ $ $\mathcal{P}$ $ $ $\in$ $ $ $H-Spec(R)$ $ $, then  $ $ $\phi_{2} \circ \mbox{ } ... \mbox{ } \circ \phi_{t}(\mathcal{P})$  $ $ $\in$ $ $ $H-Spec(R^{(2)})$. $ $ So  $ $ $\phi_{2} \circ \mbox{ } ... \mbox{ } \circ \phi_{t}(\mathcal{P})$ $ $ = $ $  $\mathcal{P}^{(2)}_{\Delta}$ $ $ with $ $ $\Delta$ $ $ an admissible diagram such that  $ $  $\mathcal{P}$ $ $ = $ $  $\mathcal{P}\it_{\Delta}$.
\end{itemize}
 
\section{New results on $ $ $H$-invariant prime ideals.}

In this section, we consider an integer  $ $ $m$ $\in$ $\llbracket  2,$ ... , $t+1 \rrbracket $, $ $ and we denote by $ $ $\mathcal{P}\it^{(m)}$ $ $ an $ $ $H$-invariant prime ideal. We set $ $ $A^{(m)}$ $ $ = $ $ $R^{(m)}/P\it^{(m)}$ $ $ and we observe that this algebra is a noetherian domain (since $ $ $R^{(m)}$ $ $ is noetherian and $ $ $P\it^{(m)}$ $ $ is completely prime). Set $ $ $D_{m}$ $ $ = $ $ $Frac(A^{(m)})$ $ $ it's division ring of fractions. Denote by $ $ $f_{m}: \mbox{ }  R^{(m)} \mbox{ } \rightarrow \mbox{ } A^{(m)}$ the canonical homomorphism and, for each  $ $ $i$ $\in$ $\llbracket  1,$ ... , $t \rrbracket $, $ $ set $ $ $x_{i}^{(m)}$ $ $ = $ $ $f_{m}(X_{i}^{(m)})$ $ $ the canonical image of $ $ $X_{i}^{(m)}$ $ $ in $ $ $A^{(m)}$. $ $ If $ $ $h$ $ $ is an element of the group $ $ $H$, $ $ we have $ $ $h$($\mathcal{P}\it^{(m)}$) $ $ = $ $ $\mathcal{P}\it^{(m)}$ $ $ and, consequently, $ $ $h$ $ $ induces an automorphism of the algebra $ $ $A^{(m)}$, $ $ denoted $ $ $\overline{h}$, $ $ which satisfies $ $ $\overline{h} \circ f_{m}$ $ $ = $ $ $f_{m} \circ h$. This automorphism can be extended in an automorphism, still denoted $ $ $\overline{h}$ $ $ of the division ring $ $ $D_{m}$.

\subsection{A necessary and sufficient condition for $ $ $\mathcal{P}\it^{(m)}$ $ $ to be in $ $ $Im(\phi_{m})$.}

Assume that $ $ $m$ $ $ $\leq$ $ $ $t$ $ $ and recall (section 3.2) that $ $ $R^{(m)}$ $ $ is the $k$ - algebra generated by the variables $ $ $X_{1}^{(m)}, \mbox{ } ... \mbox{ } , \mbox{ } X_{t}^{(m)}$ $ $ submitted to the relations \ref{LSRm}:
$$X_{j}^{(m)}X_{i}^{(m)} - \lambda_{j,i}X_{i}^{(m)}X_{j}^{(m)} \mbox{ } = \mbox{ } P_{j,i}^{(m)}(X_{i+1}^{(m)}, \mbox{ } ... \mbox{ } , \mbox{ } X_{j-1}^{(m)})$$ for $ $ $i$ $ $ $<$ $j$, and with $ $ $P_{j,i}^{(m)}(X_{i+1}^{(m)}, \mbox{ } ... \mbox{ } , \mbox{ } X_{j-1}^{(m)})$ $ $ $\in$ $ $ $k<X_{i+1}^{(m)}, \mbox{ } ... \mbox{ } , \mbox{ } X_{j-1}^{(m)}>$.

\newtheorem{prop4.1.}{Proposition 4.1.}
\begin{prop4.1.} $ $ \\
$\mathcal{P}\it^{(m)}$ $ $ $\in$ $ $ $Im(\phi_{m})$ $ $ if and only if one of two following conditions is satisfied \\
a) $ $ $X_{m}^{(m)}$ $ $ $\notin$ $ $ $\mathcal{P}\it^{(m)}$. \\
b) $ $ $X_{m}^{(m)}$ $ $ $\in$ $ $ $\mathcal{P}\it^{(m)}$ $ $ and $ $ $\Theta^{(m)}(\delta_{m}^{(m+1)}(X_{i}^{(m+1)})) $  $\in$ $ $ $\mathcal{P}\it^{(m)}$ $ $ for $ $ $1$  $\leq$ $i$ $\leq$ $m-1$.
\end{prop4.1.}

\bf Proof. \rm \\ $ $ \\
Assume that $ $ $\mathcal{P}\it^{(m)}$ $ $ $\in$ $ $ $Im(\phi_{m})$, $ $ so that $ $ $\mathcal{P}\it^{(m)}$ $ $ = $ $ $\phi_{m}$($\mathcal{P}\it^{(m+1)}$) $ $ with  $ $ $\mathcal{P}\it^{(m+1)}$ $ $ $\in$ $ $ $Spec$($R^{(m+1)}$), $ $ and assume that condition a) is not satisfied. This implies that $ $ $\mathcal{P}\it^{(m)}$ $ $ = $ $ $ker(g)$ $ $ where $ $ $g: \mbox{ } R^{(m)} \mbox{ } \rightarrow \mbox{ } R^{(m+1)} / \mathcal{P}\it^{(m+1)}$ $ $ is the homomorphism which transforms each $ $ $X_{i}^{(m)}$ $ $ in $ $ $x_{i}^{(m+1)}$ $ $ = $ $ $X_{i}^{(m+1)}$ $ $ + $ $ $\mathcal{P}\it^{(m+1)}$.
Consider $ $ $1$  $\leq$ $i$ $\leq$ $m-1$. \\ Recall that $ $ $\delta_{m}^{(m+1)}$($X_{i}^{(m+1)}$) $ $ = $ $ $P_{m,i}^{(m+1)}$($X_{i+1}^{(m+1)}, \mbox{ } ... \mbox{ } , \mbox{ } X_{m-1}^{(m+1)}$) and $ $ and that $ $ $\Theta^{(m)}: \mbox{ } k<X_{1}^{(m+1)}, \mbox{ } ... \mbox{ } , \mbox{ } X_{m-1}^{(m+1)}> \mbox{ } \rightarrow \mbox{ } k<X_{1}^{(m)}, \mbox{ } ... \mbox{ } , \mbox{ } X_{m-1}^{(m)}>$ $ $ is the homomorphism which transforms each $ $ $X_{l}^{(m+1)}$ $ $ in $ $ $X_{l}^{(m)}$. $ $ Since $ $ $X_{m}^{(m)}$ $ $ $\in$ $ $ $\mathcal{P}\it^{(m)}$, $ $ we have $ $ $X_{m}^{(m+1)}$ $ $ $\in$ $ $ $\mathcal{P}\it^{(m+1)}$ $ $ (\cite{C}, proposition 4.3.1.) and so, $ $ $\delta_{m}^{(m+1)}$($X_{i}^{(m+1)}$) $ $ $\in$ $ $ $\mathcal{P}\it^{(m+1)}$. $ $ Now, we have \\ $g$($\Theta^{(m)}$($\delta_{m}^{(m+1)}$($X_{i}^{(m+1)}$))) $ $ = $ $ $g$($\Theta^{(m)}$($P_{m,i}^{(m+1)}$($X_{i+1}^{(m+1)}, \mbox{ } ... \mbox{ } , \mbox{ } X_{m-1}^{(m+1)}$))) $ $ = $ $ $g$($P_{m,i}^{(m+1)}$($X_{i+1}^{(m)}, \mbox{ } ... \mbox{ } , \mbox{ } X_{m-1}^{(m)}$)) $ $ = $ $ $P_{m,i}^{(m+1)}$($x_{i+1}^{(m+1)}, \mbox{ } ... \mbox{ } , \mbox{ } x_{m-1}^{(m+1)})$ $ $ = $ $ $P_{m,i}^{(m+1)}$($X_{i+1}^{(m+1)}, \mbox{ } ... \mbox{ } , \mbox{ } X_{m-1}^{(m+1)})$ $ $ + $ $ $\mathcal{P}\it^{(m+1)}$ $ $ = $ $ $0$. \\ This implies that $ $ $\Theta^{(m)}$($\delta_{m}^{(m+1)}$($X_{i}^{(m+1)}$))) $ $ $\in$ $ $ $ker(g)$ $ $ = $ $ $\mathcal{P}\it^{(m)}$. \\ $ $ \\
If condition a) is satisfied, then $ $ $\mathcal{P}\it^{(m)}$ $ $ $\in$ $ $ $Im(\phi_{m})$ $ $ by (\cite{C}, lemma 4.3.1.). \\ $ $ \\
Assume that condition b) is satisfied. So, if $ $ $1$  $\leq$ $i$ $\leq$ $m-1$, $ $ we have, as previously, $ $  $P_{m,i}^{(m+1)}$($X_{i+1}^{(m)}, \mbox{ } ... \mbox{ } , \mbox{ } X_{m-1}^{(m)}$) $ $ = $ $ $\Theta^{(m)}$($\delta_{m}^{(m+1)}$($X_{i}^{(m+1)}$)) $ $ $\in$ $ $ $\mathcal{P}\it^{(m)}$. $ $ So, in $ $ $A^{(m)}$ $ $ = $ $ $R^{(m)}/\mathcal{P}\it^{(m)}$, $ $  we have  $ $  $P_{m,i}^{(m+1)}$($x_{i+1}^{(m)}, \mbox{ } ... \mbox{ } , \mbox{ } x_{m-1}^{(m)}$) $ $ = $ $ $0$. \\ Since $ $ $P_{m,i}^{(m)}$ $ $ = $ $ $0$ $ $ (see section 3.2), we can write $ $ $ $ $x_{m}^{(m)}x_{i}^{(m)}$ $ $ $-$ $ $ $\lambda_{m,i}x_{i}^{(m)}x_{m}^{(m)}$ $ $ = $ $ $P_{m,i}^{(m)}$($x_{i+1}^{(m)}, \mbox{ } ... \mbox{ } , \mbox{ } x_{m-1}^{(m)}$) $ $ = $ $ $0$ $ $ = $ $ $P_{m,i}^{(m+1)}$($x_{i+1}^{(m)}, \mbox{ } ... \mbox{ } , \mbox{ } x_{m-1}^{(m)}$). \\  If  $ $ $1$  $\leq$ $i$ $\leq$ $j-1$ $ $ with $j$ $ $ $\neq$ $ $ $m$, $ $ we have (see section 3.2) \\ $ $ $x_{j}^{(m)}x_{i}^{(m)}$ $ $ $-$ $ $ $\lambda_{j,i}x_{i}^{(m)}x_{j}^{(m)}$ $ $ = $ $ $P_{j,i}^{(m)}$($x_{i+1}^{(m)}, \mbox{ } ... \mbox{ } , \mbox{ } x_{j-1}^{(m)}$) $ $ = $ $ $P_{j,i}^{(m+1)}$($x_{i+1}^{(m)}, \mbox{ } ... \mbox{ } , \mbox{ } x_{j-1}^{(m)}$). \\ So, by the universal property of algebras defined by generators and relations, there exists a (unique) homomorphism $ $ $\epsilon: \mbox{ } R^{(m+1)} \mbox{ } \rightarrow \mbox{ } A^{(m)}$ $ $ which transforms each $ $ $X_{l}^{(m+1)}$ $ $ in $ $ $x_{l}^{(m)}$. $ $ This homomorphism is surjective,  and its kernel $ $ $ker(\epsilon)$ $ $ = $ $ $\mathcal{P}\it^{(m+1)}$ $ $ is a prime ideal of $ $ $R^{(m+1)}$. $ $ We observe that, since $ $ $X_{m}^{(m)}$ $ $ $\in$ $ $ $\mathcal{P}\it^{(m)}$, $ $ we have $ $ $X_{m}^{(m+1)}$ $ $ $\in$ $ $ $\mathcal{P}\it^{(m+1)}$, $ $  and that $ $ $\epsilon$ $ $ induces an automorphism
$$ \overline{\epsilon}: \mbox{ } R^{(m+1)} / \mathcal{P}\it^{(m+1)} \mbox{ } \rightarrow \mbox{ } A^{(m)} \mbox{ } = \mbox{ } R^{(m)} / \mathcal{P}\it^{(m)}$$ which transforms each  $ $ $x_{l}^{(m+1)}$ $ $ in $ $ $x_{l}^{(m)}$. $ $ Recall that $ $ $f_{m}: \mbox{ } R^{(m)} \mbox{ } \rightarrow \mbox{ } R^{(m)} / \mathcal{P}\it^{(m)}$  denotes the canonical homomorphism. \\ So, $g \mbox{ } = \mbox{ } (\overline{\epsilon})^{-1} \circ f_{m}: \mbox{ } R^{(m)} \mbox{ } \rightarrow \mbox{ } R^{(m+1)} / \mathcal{P}\it^{(m+1)}$ $ $ is the homomorphism which transforms each $ $ $X_{l}^{(m)}$ $ $ in $ $ $x_{l}^{(m+1)}$. $ $ As $ $ $ker(g)$ $ $ = $ $ $ker(f_{m})$ $ $ = $ $ $\mathcal{P}\it^{(m)}$, $ $ we conclude that  $ $ $\mathcal{P}\it^{(m)}$ $ $ = $ $ $\phi_{m}$($\mathcal{P}\it^{(m+1)}$). \bx

\newtheorem{cor4.1.}{Corollary 4.1.}

\begin{cor4.1.} $ $ \\
Assume $ $ $2$ $\leq$ $m$ $\leq$ $t+1$ $ $ and consider some $ $ $\cal Q$ $\in$ $Spec(R^{(m)})$. \\ Assume that $ $ $X_{m}^{(m)}$ $\notin$ $\cal Q$, $ $ ... $ $, $ $ $X_{t}^{(m)}$ $\notin$ $\cal Q$. \\
Then there exists $ $ $\mathcal{P}$ $\in$ $Spec(R)$ $ $ such that $ $ $\cal Q$ $ $ is the canonical image of $ $ $\mathcal{P}$ in $Spec(R^{(m)})$. $ $ Moreover, if $ $ $\mathcal{P}^{(2)}$ $ $ is the canonical image of $ $ $\mathcal{P}$ $ $ in $Spec(R^{(2)})$, $ $ then
$$\mathcal{P}^{(2)} \mbox{ } \cap \mbox{ } \{Z_{m}, \mbox{ } ... \mbox{ } , \mbox{ } Z_{t}\} \mbox{ } = \mbox{ } \emptyset.$$
\end{cor4.1.}

\bf Proof. \rm \\
We prove this by decreasing induction on $ $ $m$. \\ $ $ \\
\begin{itemize}
\item If $ $ $m$ = $t+1$, $ $ we have $ $ $\cal Q$ $\in$ $Spec(R)$ $ $ and $ $ $\mathcal{P}$ = $\cal Q$ $ $ satisfies the required properties.
\item Assume $ $ $2$ $\leq$ $m$ $\leq$ $t$. $ $ By proposition 4.1. 1, we have $ $ $\cal Q$ $=$ $\phi_{m} (\cal Q\it^{\prime})$ $ $ with $ $ $\cal Q\it^{\prime}$ $\in$ $Spec(R^{(m+1)})$. $ $ Since $ $ $X_{m}^{(m)}$ $\notin$ $\cal Q$, $ $ we have (see section 3.3) $ $ $X_{m}^{(m+1)}$ $\notin$ $\cal Q\it^{\prime}$ $ $ and, if we set $ $ $S_{m}$ = $\{(X_{m}^{(m+1)})^{d} \mbox{ } | \mbox{ } d$ $\in$ \bbN $\}$, $ $

$$\cal Q \it \mbox{ } = \mbox{ } R^{(m)} \mbox{ } \cap \mbox{ } \cal Q\it^{\prime} S_{m}^{-1}.$$

Assume $ $ $X_{j}^{(m+1)}$ $\in$ $\cal Q\it^{\prime}$ $ $ with $ $ $j$ $\geq$ $m$. $ $ Since $ $ $X_{j}^{(m+1)}$ = $X_{j}^{(m)}$ $ $ (see section 3.2), we have $ $ $X_{j}^{(m)}$ $\in$ $R^{(m)} \mbox{ } \cap \mbox{ } \cal Q\it^{\prime}$ $\subseteq$ $R^{(m)} \mbox{ } \cap \mbox{ } \cal Q\it^{\prime} S_{m}^{-1}$ = $\cal Q$, $ $ which is false. So, we have:
\begin{itemize}
\item $X_{m}^{(m+1)}$ $\notin$ $\cal Q\it^{\prime}$.
\item $X_{m+1}^{(m+1)}$ $\notin$ $\cal Q\it^{\prime}$, $ $ ... $ $, $ $ $X_{t}^{(m+1)}$ $\notin$ $\cal Q\it^{\prime}$ $ $ and, by the induction assumption, there exists $ $ $\mathcal{P}$ $\in$ $Spec(R)$ $ $ such that $ $ $\cal Q\it^{\prime}$ $ $ is the canonical image of $ $ $\mathcal{P}$ $ $ in $Spec(R^{(m+1)})$. $ $ This implies that $ $ $\cal Q$ $=$ $\phi_{m} (\cal Q\it^{\prime})$ $ $ is the canonical image of $ $ $\mathcal{P}$ $ $ in $Spec(R^{(m)})$. $ $  Moreover, again by the induction assumption, we have
$$\mathcal{P}^{(2)} \mbox{ } \cap \mbox{ } \{Z_{m+1}, \mbox{ } ... \mbox{ } , \mbox{ } Z_{t}\} \mbox{ } = \mbox{ } \emptyset.$$
As $ $ $X_{m}^{(m+1)}$ $\notin$ $\cal Q\it^{\prime}$, $ $ we know (\cite{C}, proposition 5.2.1.) that $ $ $Z_{m}$ $\notin$ $\mathcal{P}^{(2)}$.
\end{itemize}
\end{itemize}
 \bx
 
\subsection{Some  properties of $ $ $A^{(m) }$.}
 
For each integer  $ $ $j$ $\in$ $\llbracket  2,$ ... , $m-1 \rrbracket $, $ $ we set $ $ $\mathcal{P}\it^{(j)}$ $ $ = $ $ $\phi_{j} \circ \mbox{ } ... \mbox{ } \circ \phi_{m-1}$($\mathcal{P}\it^{(m)}$). $ $ Let us, first, construct in $ $ $A^{(m)}$, $ $ a new version of the deleting derivations algorithm.

\newtheorem{prop4.2.}{Proposition 4.2.}
\begin{prop4.2.} $ $ \\
For each integer  $ $ $j$ $\in$ $\llbracket  2,$ ... , $m-1 \rrbracket $, $ $ there exists a sequence $ $ $(x_{1}^{(j)}, \mbox{ } ... \mbox{ } , \mbox{ } x_{t}^{(j)}) $ in $ $ $D_{m}$ $ $ = $ $ $Frac(A^{(m)})$ $ $ and a subalgebra $A^{(j)}$ = $k<x_{1}^{(j)}, \mbox{ } ... \mbox{ } , \mbox{ } x_{t}^{(j)}>$ of $D_{m}$ which satisfy the following properties.
\begin{enumerate}
\item If $ $ $l$ $ $ and $ $ $i$ $ $ are in $ $ $\llbracket  1,$ ... , $t \rrbracket $, $ $ we have $ $ $\overline{h_{l}}(x_{i}^{(j)})$ $ $ = $ $ $\lambda_{l,i}x_{i}^{(j)}$, $ $ $\overline{h_{l}}(A^{(j)})$ $ $ = $ $ $A^{(j)}$ $ $ and $ $ $\overline{h_{l}}$ $ $ induces, by restriction, an automorphism $($still denoted $ $ $\overline{h_{l}}) $ of $ $ $A^{(j)}$. \\
\item Choose $ $ $l$ $\in$ $\llbracket  1,$ ... , $t \rrbracket $ $ $ and denote by $ $ $d_{l}^{(j)}$ $ $ the left $ $ $\overline{h_{l}}$ - derivation of $ $ $A^{(j)}$ $ $ associated to $ $ $x_{l}^{(j)}$ $ $ $(d_{l}^{(j)}(a)$ $ $ = $ $  $x_{l}^{(j)}a$ $-$ $\overline{h_{l}}(a)x_{l}^{(j)}$ $ $ for all $ $ $a$ $ $ $\in$ $ $ $A^{(j)})$. $ $ Then $ $ $d_{l}^{(j)}(k<x_{1}^{(j)}, \mbox{ } ... \mbox{ } , \mbox{ } x_{l-1}^{(j)}>)$ $ $ $\subset$ $ $ $k<x_{1}^{(j)}, \mbox{ } ... \mbox{ } , \mbox{ } x_{l-1}^{(j)}>$ $ $ and $ $ $d_{l}^{(j)}$ $ $ is locally nilpotent on $ $ $k<x_{1}^{(j)}, \mbox{ } ... \mbox{ } , \mbox{ } x_{l-1}^{(j)}>$. \\
\item There exists a unique homomorphism $ $ $f_{j}: \mbox{ } R^{(j)} \mbox{ } \rightarrow \mbox{ } A^{(j)}$ $ $ which transforms each  $ $ $X_{i}^{(j)}$ $ $ in $ $ $x_{i}^{(j)}$. $ $ $f_{j}$ $ $ is surjective and $ $ $ker(f_{j})$ $ $ = $ $ $\mathcal{P}\it^{(j)}$. \\
\item If $ $ $l$ $\in$ $\llbracket  1,$ ... , $t \rrbracket $, $ $ we have $ $ $d_{l}^{(j)} \circ f_{j}$ $ $ = $ $ $f_{j} \circ \delta_{l}^{(j)}$ $ $ and $ $ $\overline{h_{l}} \circ f_{j}$ $ $ = $ $ $f_{j} \circ h_{l}$. \\
\item  If $ $ $l, e$ $\in$ $\llbracket  1,$ ... , $t \rrbracket $, $ $ we have $ $ $\overline{h_{l}} \circ \overline{h_{e}}$ $ $ = $ $ $\overline{h_{e}} \circ \overline{h_{l}}$ $ $ and $ $ $\overline{h_{l}} \circ d_{e}^{(j)}$ $ $ = $ $ $\lambda_{l,e} d_{e}^{(j)} \circ \overline{h_{l}}.$ \\
\item If $ $ $x_{j}^{(j+1)}$ $ $ = $ $ $0$, $ $ then $ $ $x_{i}^{(j)}$ $ $ = $ $ $x_{i}^{(j+1)}$ $ $ for each $ $ $i$ $\in$ $\llbracket  1,$ ... , $t \rrbracket $. \\
\item If $ $ $x_{j}^{(j+1)}$ $ $ $\neq$ $ $ $0$, $ $ and $ $ $i$ $\in$ $\llbracket  1,$ ... , $t \rrbracket $, $ $ then

\begin{enumerate}
\item $j \leq i \mbox{ } \Rightarrow \mbox{ } x_{i}^{(j)}  = \mbox{ } x_{i}^{(j+1)}$. $ $ In particular, we have $ $ $x_{j}^{(j)} = \mbox{ } x_{j}^{(j+1)}$. \\
\item $i < j  \mbox{ } \Rightarrow \mbox{ } x_{i}^{(j)} = \mbox{ } x_{i}^{(j+1)} \mbox{ } \mbox{ } +  \mbox{ } \mbox{ } \displaystyle \sum_{d = 1}^{+\infty} c_{d}^{(j+1)}(x_{j}^{(j+1)})^{-d} \mbox{ } \mbox{ }  \mbox{ } $
$with \mbox{ } \mbox{ } \mbox{ } \mbox{ } c_{d}^{(j+1)} = \mbox{ }  \displaystyle \frac{(1-q_{j})^{-d}}{[d]!_{q_{j}}} \mbox{ } \lambda_{j,i}^{-d}  \mbox{ } (d_{j}^{(j+1)})^{d} \mbox{ } (x_{i}^{(j+1)}).$
\end{enumerate}
\rm
(We observe that, since $ $ $k<x_{1}^{(j+1)},  \mbox{ } ... \mbox{ }, \mbox{ } x_{j-1}^{(j+1)}>$ $ $ is $ $ $d_{j}^{(j+1)}$ - stable, we have $ $ $c_{d}^{(j+1)}$ $\in$ $ $ $k<x_{1}^{(j+1)},  \mbox{ } ... \mbox{ }, \mbox{ } x_{j-1}^{(j+1)}>$ $ $ and, since $ $ $d_{j}^{(j+1)}$ $ $ is locally nilpotent on  $ $ $k<x_{1}^{(j+1)},  \mbox{ } ... \mbox{ }, \mbox{ } x_{j-1}^{(j+1)}>$, $ $ only finitely many  $ $ $c_{d}^{(j+1)}$ $ $ are non zero.)

\end{enumerate}
\end{prop4.2.}

\bf Proof. \rm \\ $ $ \\

First observe that, if 3. is satisfied and if $ $ $\overline{h_{l}} \circ f_{j}$ $ $ = $ $ $f_{j} \circ h_{l}$ $ $ (as homomorphism from $ $ $R^{(j)}$ $ $ towards $ $ $D^{(m)}$) $ $for any $ $ $l$, $ $ then properties 1. to 5. are all satisfied. In fact, \\
\begin{description}
\item[1. ] For all $ $ $l$ $ $ and $ $ $i$ $ $ in $ $ $\llbracket  1,$ ... , $t \rrbracket $, $ $ we have $ $ $\overline{h_{l}}(x_{i}^{(j)})$ $ $ = $ $ $\overline{h_{l}} \circ f_{j}(X_{i}^{(j)})$ $ $ = $ $ $f_{j} \circ h_{l}(X_{i}^{(j)})$ $ $ = $ $ $\lambda_{l,i}f_{j}(X_{i}^{(j)})$ $ $ = $ $ $\lambda_{l,i}x_{i}^{(j)}$. This immediately implies that $ $ $\overline{h_{l}}(A^{(j)})$ $ $ = $ $ $A^{(j)}$ $ $ and $ $ $\overline{h_{l}}$ $ $ induces, by restriction, an automorphism of $ $ $A^{(j)}$. 
\item[4. ] If $ $ $P$ $\in$ $R^{(j)}$, $ $ we have $ $ $d_{l}^{(j)}$($f_{j}(P)$) $ $ = $ $  $x_{l}^{(j)}f_{j}(P)$ $-$ $\overline{h_{l}}(f_{j}(P))x_{l}^{(j)}$, $ $ = $ $  $f_{j}$($X_{l}^{(j)}P$ $-$ $h_{l}(P)X_{l}^{(j)}$) $ $ = $ $  $f_{j} \circ \delta_{l}^{(j)}$($P$). $ $  This proves that $ $ $d_{l}^{(j)} \circ f_{j}$ $ $ = $ $  $f_{j} \circ \delta_{l}^{(j)}$.
\item[2. ] By 4., we have $ $ $d_{l}^{(j)}$($k<x_{1}^{(j)}, \mbox{ } ... \mbox{ } , \mbox{ } x_{l-1}^{(j)}>$) $ $ = $ $ $d_{l}^{(j)}$($f_{j}$($k<X_{1}^{(j)}, \mbox{ } ... \mbox{ } , \mbox{ } X_{l-1}^{(j)}>$)) $ $ =  $ $ $f_{j}$($\delta_{l}^{(j)}$($k<X_{1}^{(j)}, \mbox{ } ... \mbox{ } , \mbox{ } X_{l-1}^{(j)}>$)) $ $  $\subset$ $ $ $f_{j}$($k<X_{1}^{(j)}, \mbox{ } ... \mbox{ } , \mbox{ } X_{l-1}^{(j)}>$) $ $ = $ $ $k<x_{1}^{(j)}, \mbox{ } ... \mbox{ } , \mbox{ } x_{l-1}^{(j)}>$. $ $ Moreover, if $ $ $a$ $\in$ $k<x_{1}^{(j)}, \mbox{ } ... \mbox{ } , \mbox{ } x_{l-1}^{(j)}>$, $ $ we can write $ $ $a$ $ $ = $ $ $f_{j}(P)$ $ $ with $ $ $P$ $\in$ $ $ $k<X_{1}^{(j)}, \mbox{ } ... \mbox{ } , \mbox{ } X_{l-1}^{(j)}>$. $ $ Since $ $ $\delta_{l}^{(j)}$ $ $ is locally nilpotent on  $ $ $k<X_{1}^{(j)}, \mbox{ } ... \mbox{ } , \mbox{ } X_{l-1}^{(j)}>$, $ $ we have $ $ $(\delta_{l}^{(j)})^{r}$($P$) $ $ = $ $ $0$ $ $ for some $ $ $r$ $\in$ $ $ \bbN , $ $ and $ $ $(d_{l}^{(j)})^{r}$($a$) $ $ =  $ $ $(d_{l}^{(j)})^{r} \circ f_{j}$($P$) $ $ = $ $  $f_{j} \circ (\delta_{l}^{(j)})^{r}$($P$) $ $ =  $ $ $0$. $ $ So, $ $ $d_{l}^{(j)}$ $ $ is locally nilpotent on $ $ $k<x_{1}^{(j)}, \mbox{ } ... \mbox{ } , \mbox{ } x_{l-1}^{(j)}>$.
\item[5. ] If $ $ $i$ $ $ in $ $ $\llbracket  1,$ ... , $t \rrbracket $, $ $ we have $ $ $\overline{h_{l}} \circ \overline{h_{e}} (x_{i}^{(j)})$ $ $ = $ $ $\lambda_{l,i} \lambda_{e,i} x_{i}^{(j)}$ $ $ =  $ $ $\overline{h_{e}} \circ \overline{h_{l}} (x_{i}^{(j)})$ and, since $ $ $A^{(j)}$ $ $ is generated by the elements $ $ $x_{i}^{(j)}$, $ $ we have $ $ $\overline{h_{l}} \circ \overline{h_{e}}$ $ $ = $ $ $\overline{h_{e}} \circ \overline{h_{l}}$. \\
Consider two elements $ $ $a, \mbox{ } b$ $ $ in  $ $ $A^{(j)}$ $ $ such that the maps $ $ $\overline{h_{l}} \circ d_{e}^{(j)}$ $ $ and $ $ $\lambda_{l,e} d_{e}^{(j)} \circ \overline{h_{l}}$ $ $ coincide on $ $ $a$ $ $ and  $ $ $b$. $ $ Then, we have $ $ $\overline{h_{l}} \circ d_{e}^{(j)} (ab)$ $ $ = $ $ $\overline{h_{l}}(d_{e}^{(j)}(a)b \mbox{ } + \mbox{ } \overline{h_{e}}(a)d_{e}^{(j)}(b))$ $ $ = $ $ $\lambda_{l,e} d_{e}^{(j)}(\overline{h_{l}}(a))\overline{h_{l}}(b) \mbox{ } + \mbox{ } \lambda_{l,e} \overline{h_{e}}(\overline{h_{l}}(a)) d_{e}^{(j)}(\overline{h_{l}}(b))$ $ $ = $ $ $\lambda_{l,e} d_{e}^{(j)} \circ \overline{h_{l}} (ab).$ $ $ So, we just have to prove that $ $ $\overline{h_{l}} \circ d_{e}^{(j)}$ $ $ and $ $ $\lambda_{l,e} d_{e}^{(j)} \circ \overline{h_{l}}$ $ $ coincide on each $ $ $x_{i}^{(j)}$. Now,  $ $
$\overline{h_{l}} \circ d_{e}^{(j)} (x_{i}^{(j)})$ $ $ = $ $ $\overline{h_{l}}$($x_{e}^{(j)}x_{i}^{(j)} \mbox{ } - \mbox{ } \lambda_{e,i}x_{i}^{(j)}x_{e}^{(j)}$) $ $ = $ $ $\lambda_{l,e} \lambda_{l,i}$($x_{e}^{(j)}x_{i}^{(j)} \mbox{ } - \mbox{ } \lambda_{e,i}x_{i}^{(j)}x_{j}^{(j)}$) $ $ = $ $ $\lambda_{l,e} \lambda_{l,i} d_{e}^{(j)}(x_{i}^{(j)})$ $ $ = $ $ $ \lambda_{l,e} d_{e}^{(j)}(\overline{h_{l}}(x_{i}^{(j)}))$. So we can conclude that $ $ $\overline{h_{l}} \circ d_{e}^{(j)}$ $ $ = $ $ $\lambda_{l,e} d_{e}^{(j)} \circ \overline{h_{l}}$.
\end{description}
$ $ \\

For  $ $ $j$ $\in$ $\llbracket  2,$ ... , $m-1 \rrbracket $, $ $ denote by $ $ $\cal H\it_{j}$ $ $ the following assumption: \\ $ $ \\
\it There exists a sequence $ $ $(x_{1}^{(j+1)}, \mbox{ } ... \mbox{ } , \mbox{ } x_{t}^{(j+1)})$ $ $ $\in$ $ $ $(D_{m})^{t}$ $ $ and a subalgebra $A^{(j+1)}$ = $k<x_{1}^{(j+1)}, \mbox{ } ... \mbox{ } , \mbox{ } x_{t}^{(j+1)}>$ of $D_{m}$ which satisfy the analog of properties 1. to 5. obtained by changing $j$ in $j+1$. \rm \\ $ $ \\
Observe that the assumption  $ $ $\cal H\it_{m-1}$ $ $ is satisfied. In fact, the constructions of  $ $ $f_{m}$ $ $ and of the variables $ $ $x_{i}^{(m)}$ $ $ (see the beginning of section 4) imply straightway the property 3. and the equality $ $ $\overline{h} \circ f_{m}$ $ $ = $ $ $f_{m} \circ h$ for any $ $ $h$ $ $ in $ $ $H$. \\
Now, assume that $ $ $\cal H\it_{j}$ $ $ is satisfied, define the elements $ $ $x_{i}^{(j)}$ $ $ as in 6. and 7. and set $ $ $A^{(j)}$ $ $ = $ $ $ $ $k<x_{1}^{(j)}, \mbox{ } ... \mbox{ } , \mbox{ } x_{t}^{(j)}>$. $ $ It remains to prove that property 3. is satisfied and that $ $ $\overline{h_{l}} \circ f_{j}$ $ $ = $ $ $f_{j} \circ h_{l}$. \\ $ $ \\

\begin{itemize}
\item Assume that $ $ $x_{j}^{(j+1)}$ $ $ = $ $ $0$, $ $ so that $ $ $X_{j}^{(j+1)}$ $ $ $\in$ $ $ $Ker(f_{j+1})$ $ $ = $ $ $\mathcal{P}\it^{(j+1)}$. $ $ Recall that, in this case, $ $ $\mathcal{P}\it^{(j)}$ $ $ is the kernel of the homomorphism $ $ $g: \mbox{ } R^{(j)} \mbox{ } \rightarrow \mbox{ } R^{(j+1)} / \mathcal{P}\it^{(j+1)}$ $ $ which transforms each $ $ $X_{i}^{(j)}$ $ $ in $ $ $X_{i}^{(j+1)}$ + $\mathcal{P}\it^{(j+1)}$. \\ By the assumption $ $ $\cal H\it_{j}$, $ $ $f_{j+1}:  \mbox{ } R^{(j+1)} \mbox{ } \rightarrow \mbox{ } A^{(j+1)}$ $ $ induces an isomorphism $ $ $\widehat{f_{j+1}}:  \mbox{ } R^{(j+1)} / \mathcal{P}\it^{(j+1)} \mbox{ } \rightarrow \mbox{ } A^{(j+1)}$ $ $ which transforms each $ $ $X_{i}^{(j+1)}$ + $\mathcal{P}\it^{(j+1)}$ $ $ in $ $ $x_{i}^{(j+1)}$. $ $ So $ $ $f_{j}$ $ $ = $ $ $\widehat{f_{j+1}} \circ g: \mbox{ } R^{(j)} \mbox{ } \rightarrow \mbox{ } A^{(j+1)}$ $ $ is the homomorphism which transforms each $ $ $X_{i}^{(j)}$ $ $ in $ $ $x_{i}^{(j+1)}$, $ $ and we have $ $ $ker(f_{j})$ $ $ = $ $ $ker(g)$ $ $ = $ $ $\mathcal{P}\it^{(j)}$. $ $ Since, in this case, we have $ $ $x_{i}^{(j)}$ $ $ = $ $ $x_{i}^{(j+1)}$ $ $ for each $ $ $i$, we obtain that $f_{j}$ satisfies the property 3. Moreover, if  $l$  and  $i$ are in $\llbracket  1,$ ... , $t \rrbracket $, $ $ we have $ $ $\overline{h_{l}} \circ f_{j}  (X_{i}^{(j)})$ $ $ = $ $ $\overline{h_{l}}(x_{i}^{(j)})$ $ $ = $ $ $\overline{h_{l}}(x_{i}^{(j+1)})$ $ $ = $ $ $\lambda_{l,i}(x_{i}^{(j+1)})$ $ $ =  $ $ $\lambda_{l,i}(x_{i}^{(j)})$ $ $ = $ $ $f_{j} \circ  h_{l} (X_{i}^{(j)})$. $ $ So, $ $ $\overline{h_{l}} \circ  f_{j}$ $ $ = $ $ $f_{j}  \circ  h_{l}$.
\item Assume that $ $ $x_{j}^{(j+1)}$ $ $ $\neq$ $ $ $0$, $ $ so that $ $ $X_{j}^{(j+1)}$ $ $ $\notin$ $ $ $Ker(f_{j+1})$ $ $ = $ $ $\mathcal{P}\it^{(j+1)}$. $ $ Set $ $ $S_{j}$ $ $ = $ $ $\{(X_{j}^{(j+1)})^{l} \mbox{ } | \mbox{ } l \mbox{ } \in $ $ $ \bbN $\}$ $ $ and recall that
$$ \mathcal{P}\it^{(j)} \mbox{ } = \mbox{ } \mathcal{P}\it^{(j+1)}S_{j}^{-1} \cap R^{j}.$$
As $ $ $f_{j+1}(X_{j}^{(j+1)})$ $ $ = $ $ $x_{j}^{(j+1)}$ $ $ is invertible in $ $ $D_{m}$, $ $ $ $ $f_{j+1}$ $ $ can be extended in an homomorphism $ $ $\widetilde{f_{j+1}}: \mbox{ } R^{(j+1)}S_{j}^{-1} \mbox{ } \rightarrow \mbox{ } D_{m}$ $ $ $ $ $ $ ($\widetilde{f_{j+1}}(Q(X_{j}^{(j+1)})^{-l})$ $ $ = $ $ $f_{j+1}(Q)(x_{j}^{(j+1)})^{-l}$ $ $ for all $Q$ $ $ in $ $ $R^{(j+1)}$ and $l$ $ $ in $ $ \bbN) $ $ $ $ $ $ and we have $ $ $ker(\widetilde{f_{j+1}})$ $ $ = $ $ $\mathcal{P}\it^{(j+1)}S_{j}^{-1}$. $ $ Recall that $ $ $R^{(j)}$ $ $ $\subset$ $ $ $R^{(j+1)}S_{j}^{-1}$ $ $ and let us compute $ $ $\widetilde{f_{j+1}}(X_{i}^{(j)})$ $ $ for each $ $ $i$ $\in$ $\llbracket 1,$ ... , $t \rrbracket $. \\ $ $ \\
\begin{itemize}
\item If $ $ $j$ $\leq$ $i$, then we have $ $ $X_{i}^{(j)}$ $ $ = $ $ $X_{i}^{(j+1)}$ $ $ (see section 3.2) and $ $ $\widetilde{f_{j+1}}(X_{i}^{(j)})$ $ $ = $ $ $x_{i}^{(j+1)}$ $ $ = $ $ $x_{i}^{(j)}$.
\item If $ $ $i$ $<$ $j$, then we have
$$X_{i}^{(j)} = \mbox{ } X_{i}^{(j+1)} \mbox{ } \mbox{ } +  \mbox{ } \mbox{ } \displaystyle \sum_{d = 1}^{+\infty} C_{d}^{(j+1)}(X_{j}^{(j+1)})^{-d}$$
with
$$C_{d}^{(j+1)} = \mbox{ }  \displaystyle \frac{(1-q_{j})^{-d}}{[d]!_{q_{j}}} \mbox{ } \lambda_{j,i}^{-d}  \mbox{ } (\delta_{j}^{(j+1)})^{d} \mbox{ } (X_{i}^{(j+1)}).$$
As $ $ $C_{d}^{(j+1)}$ $\in$ $R^{(j+1)}$, $ $ we have $ $ $\widetilde{f_{j+1}}(C_{d}^{(j)})$ $ $ = $ $ $f_{j+1}(C_{d}^{(j)})$ $ $ = $ $ $\displaystyle \frac{(1-q_{j})^{-d}}{[d]!_{q_{j}}} \mbox{ } \lambda_{j,i}^{-d}  \mbox{ } f_{j+1} \circ (\delta_{j}^{(j+1)})^{d} \mbox{ } (X_{i}^{(j+1)})$ $ $ =  $ $ $\displaystyle \frac{(1-q_{j})^{-d}}{[d]!_{q_{j}}} \mbox{ } \lambda_{j,i}^{-d}  \mbox{ } (d_{j}^{(j+1)})^{d} \circ  f_{j+1} (X_{i}^{(j+1)})$ $ $ = $ $ $\displaystyle \frac{(1-q_{j})^{-d}}{[d]!_{q_{j}}} \mbox{ } \lambda_{j,i}^{-d}  \mbox{ } (d_{j}^{(j+1)})^{d} \mbox{ } (x_{i}^{(j+1)})$ $ $ = $ $ $c_{d}^{(j+1)}$. $ $ So, we get $ $ $$\widetilde{f_{j+1}}(X_{i}^{(j)}) \mbox{ } = \mbox{ } x_{i}^{(j+1)} \mbox{ } \mbox{ } +  \mbox{ } \mbox{ } \displaystyle \sum_{d = 1}^{+\infty} c_{d}^{(j+1)}(x_{j}^{(j+1)})^{-d} \mbox{ } = \mbox{ } x_{i}^{(j)}$$
\end{itemize}
So, we always have $ $ $\widetilde{f_{j+1}}(X_{i}^{(j)})$ $ $ = $ $ $x_{i}^{(j)}$ $ $ and the restriction $ $ $f_{j}$ $ $ of  $ $ $\widetilde{f_{j+1}}$ $ $ to $ $ $R^{(j)}$ $ $ is the unique homomorphism from $ $ $R^{(j)}$ $ $ to $ $ $A^{(j)}$ $ $ which transforms each $ $ $X_{i}^{(j)}$ $ $ in $ $ $x_{i}^{(j)}$. $ $ $f_{j}$ $ $ is surjective and $ $ $ker(f_{j})$ $ $ = $ $ $ker(\widetilde{f_{j+1}}) \mbox{ } \cap \mbox{ } R^{(j)}$ $ $ = $ $ $\mathcal{P}\it^{(j+1)}S_{j}^{-1} \mbox{ } \cap \mbox{ } R^{(j)}$ $ $ = $ $ $\mathcal{P}\it^{(j)}$. This proves that $ $ $f_{j}$ $ $ satisfies the property 3 and it remains to verify that, for any  $ $ $l$ $ $ in $\llbracket  1,$ ... , $t \rrbracket $, $ $ we have $ $ $\overline{h_{l}} \circ f_{j}$ $ $ = $ $ $f_{j} \circ h_{l}$. $ $ As $ $ $h_{l}$($X_{i}^{(j)}$) $ $ = $ $ $\lambda_{l,i} X_{i}^{(j)}$, $ $ it is enough to check that $ $ $\overline{h_{l}}$($x_{i}^{(j)}$) $ $ = $ $ $\lambda_{l,i} x_{i}^{(j)}$ $ $ for each $ $ $i$. \\ 
\end{itemize}
If $j \leq i$, we have $\overline{h_{l}}(x_i^{(j)}) = \overline{h_{l}}(x_i^{(j+1)}) = \lambda_{l,i} x_i^{(j+1)} = \lambda_{l,i} x_i^{(j)}$. \\ 
Now, assume $i < j$.

By the assumption  $ $ $\cal H\it_{j}$, $ $ we can use property 5. at the rank $j+1$, so that
$$\overline{h_{l}} \mbox{ } \circ  \mbox{ } d_{j}^{(j+1)} \mbox{ } = \mbox{ } \lambda_{l,j} d_{j}^{(j+1)} \mbox{ } \circ \mbox{ } \overline{h_{l}}.$$
For each $d$ $ $ $\in$  $ $ \bbN, $ $ we have  
$$\overline{h_{l}}(c_{d}^{(j+1)}) \mbox{ } =  \mbox{ } \overline{h_{l}}(\displaystyle \frac{(1-q_{j})^{-d}}{[d]!_{q_{j}}} \mbox{ } \lambda_{j,i}^{-d}  \mbox{ } (d_{j}^{(j+1)})^{d} \mbox{ } (x_{i}^{(j+1)}))$$
$$=  \mbox{ } \displaystyle \frac{(1-q_{j})^{-d}}{[d]!_{q_{j}}} \mbox{ } \lambda_{j,i}^{-d}  \mbox{ } \lambda_{l,j}^{d}(d_{j}^{(j+1)})^{d} \circ \overline{h_{l}}(x_{i}^{(j+1)})  \mbox{ } =   \mbox{ } \lambda_{l,i} \lambda_{l,j}^{d} c_{d}^{(j+1)}.$$
  So, $ $ $\overline{h_{l}}$($ $ $c_{d}^{(j+1)}(x_{j}^{(j+1)})^{-d})$ $ $ = $ $ $\lambda_{l,i} c_{d}^{(j+1)} (x_{j}^{(j+1)})^{-d}$ $ $, which implies that $ $ $\overline{h_{l}}$($x_{i}^{(j)}$) $ $ = $ $ $\lambda_{l,i} x_{i}^{(j)}$. \bx  $ $ \\

Recall (see section 3.3) that, since $ $ $\mathcal{P}\it^{(m)}$ $ $ is $H$ - invariant, each $ $ $\mathcal{P}\it^{(j)}$ $ $ ($2 \mbox{ } \leq \mbox{ } j \mbox{ } \leq \mbox{ } m$) $ $ is $H$ - invariant. In particular, $ $ $\mathcal{P}^{(2)}$ $ $ is $H$ - invariant. Set

$$\Delta  \mbox{ } = \mbox{ } \{i \in \mbox{ } \rm [\it\!|1,  \mbox{ } ...  \mbox{ } ,  \mbox{ } t|\!\rm] \it \mbox{ } | \mbox{ } Z_{i} \in \mbox{ } \mathcal{P}^{(2)} \} \mbox{ } = \mbox{ } \{j_{1} \mbox{ } < \mbox{ } ...\mbox{ } <\mbox{ } j_{s}\} \rm \mbox{ } (unless \mbox{ } this \mbox{ } set \mbox{ } is \mbox{ } empty)$$
and
$$\overline{\Delta} \mbox{ }  = \mbox{ } \rm [\it\!|1,  \mbox{ } ...  \mbox{ } ,  \mbox{ } t|\!\rm] \it \mbox{ } \setminus \mbox{ } \Delta \mbox{ } = \mbox{ } \{l_{1} \mbox{ } < \mbox{ } ... \mbox{ } < \mbox{ } l_{e}\} \rm \mbox{ } (unless \mbox{ } this \mbox{ } set \mbox{ } is \mbox{ } empty).$$

Recall that $\mathcal{P}^{(2)} = \mathcal{P}_{\Delta}^{(2)}$ (see section 3.3), $ $ $\overline{R}$ $ $ = $ $ $R^{(2)}$ $ $ and set $ $ $\overline{A}$ $ $ = $ $ $A^{(2)}$. $ $ For each $ $ $i \in \mbox{ } \rm [\it\!|1,  \mbox{ } ...  \mbox{ } ,  \mbox{ } t|\!\rm]$, $ $ set $ $ $z_{i}$ $ $ = $ $ $x_{i}^{(2)}$ $ $ and observe that $ $ $f_{2}: \mbox{ } \overline{R} \mbox{ } \rightarrow \mbox{ } \overline{A}$ $ $  transforms each  $ $ $Z_{i}$ = $X_{i}^{(2)}$ $ $ in $ $ $z_{i}$. $ $ So, we can also describe $ $ $\Delta $ $ $ and $ $ $\overline{\Delta}$ $ $ as follows:

$$\Delta  \mbox{ } = \mbox{ } \{i \in \mbox{ } \rm [\it\!|1,  \mbox{ } ...  \mbox{ } ,  \mbox{ } t|\!\rm] \it \mbox{ } | \mbox{ } z_{i} \mbox{ } = \mbox{ } 0\}$$

$$\overline{\Delta}  \mbox{ } = \mbox{ } \{i \in \mbox{ } \rm [\it\!|1,  \mbox{ } ...  \mbox{ } ,  \mbox{ } t|\!\rm] \it \mbox{ } | \mbox{ } z_{i} \mbox{ } \neq \mbox{ } 0\}$$

\newtheorem{obs4.2.}{Observation 4.2.}

\begin{obs4.2.} $ $ \\
Assume that $\overline{\Delta}$ $ $ is empty. Then each $ $ $x_{i}^{(m)}$ $ $ $(i$ $\in$  $[\it\!|1,  \mbox{ } ...  \mbox{ } ,  \mbox{ } t|\!\rm])$ $ $ is zero.
\end{obs4.2.}

\newpage

\bf Proof \rm \\
If $ $ $\overline{\Delta}$ $ $ is empty, we have $ $ $x_{j}^{(j+1)}$ = $z_{j}$ = $0$ $ $ for each $ $ $j \in \mbox{ } \rm [\it\!|1,  \mbox{ } ...  \mbox{ } ,  \mbox{ } m-1|\!\rm]$. $ $ By (proposition 4.2. 1, 6.), this implies that, if $ $  $i$ $\in$  $[\it\!|1,  \mbox{ } ...  \mbox{ } ,  \mbox{ } t|\!\rm]$ $ $ and $ $ $j \in \mbox{ } \rm [\it\!|1,  \mbox{ } ...  \mbox{ } ,  \mbox{ } m-1|\!\rm]$, $ $ we have $ $ $x_{i}^{(j+1)}$ =  $x_{i}^{(j)}$. $ $ So, each $ $ $x_{i}^{(m)}$  = $x_{i}^{(2)}$ = $z_{i}$ = $0$. \bx

Until the end of section 4.2, we assume that $ $ $\overline{\Delta}$ $ $ is nonempty. So, we have:

\begin{prop4.2.} $ $ \\
\begin{enumerate}
\item $ $ $\overline{A}$ $ $ = $ $ $k<z_{l_{1}}, \mbox{ } ... \mbox{ } , \mbox{ } z_{l_{e}}>.$
\item If $ $ $1 \mbox{ } \leq i < \mbox{ } d \leq \mbox{ } e$, $ $ we have $ $ $z_{l_{d}}z_{l_{i}}$ $ $ = $ $ $\lambda_{l_{d},l_{i}}z_{l_{i}}z_{l_{d}}$.
\item $z_{l_{1}}, \mbox{ } ... \mbox{ } , \mbox{ } z_{l_{e}}$ $ $ are all nonzero and the Laurent $(ordered)$ monomials $ $ $z^{\underline{a}}$ $ $ = $ $ $z_{l_{1}}^{a_{1}} \mbox{ } ... \mbox{ } z_{l_{e}}^{a_{e}}$ $ $ $(\underline{a}$ $ $ =$ $ $(a_{1}, \mbox{ } ... \mbox{ } , \mbox{ } a_{e})$ $ $ $\in$ $ $ \bbZ$^{e})$ are $ $ $k$ - linearly independent.
\end{enumerate}
\end{prop4.2.} $ $ \\

\bf Proof. \rm \\ $ $ \\

\begin{enumerate}
\item $\overline{A}$ $ $ = $ $ $k<z_{1}, \mbox{ } ... \mbox{ } , \mbox{ } z_{t}>$ $ $ (see proposition 4.2. 1) and, for $ $ $i$ $ $ $\notin$ $ $ $\{l_{1}, \mbox{ } ... \mbox{ } , \mbox{ } l_{e}\}$, $ $ we have $ $ $i$ $ $ $\in$ $ $ $\Delta$ $ $ $\Rightarrow$ $ $ $Z_{i}$ $ $ $\in$ $ $ $\mathcal{P}^{(2)}$ $ $ = $ $ $ker(f_{2})$ $ $ $\Rightarrow$ $ $ $z_{i}$ $ $ = $ $ $0$.
\item We know (section 3.2) that $ $ $Z_{l_{d}}Z_{l_{i}}$ $ $ = $ $ $\lambda_{l_{d},l_{i}}Z_{l_{i}}Z_{l_{d}}$. $ $ If we transform by $ $ $f_{2}$, $ $ we obtain the required equality.
\item Denote by $ $ $\widehat{R}$ $ $ = $ $ $k<Z_{l_{1}}, \mbox{ } ... \mbox{ } , \mbox{ } Z_{l_{e}}>$ $ $ the subalgebra of $ $ $\overline{R}$ $ $ generated by $ $ $Z_{l_{1}}, \mbox{ } ... \mbox{ } , \mbox{ } Z_{l_{e}}$. $ $ By the property 1., $ $ $f_{2}$ $ $ induces (by restriction) a surjective homomorphism  $ $ $\widehat{f_{2}}: \mbox{ } \widehat{R} \mbox{ } \rightarrow \mbox{ } \overline{A}$ $ $ and $ $ $ker(\widehat{f_{2}}) \mbox{ } = \mbox{ }  \widehat{R} \cap \mathcal{P}^{(2)}$. $ $ As each $ $ $P$ $ $ $\in$ $ $ $\mathcal{P}^{(2)}$ = $\mathcal{P}_\Delta^{(2)}$ = $(\{Z_{j_1},...,Z_{j_s}\})$ $ $ is a linear combination of monomials in which at least one of the variables $ $ $Z_{j_{i}}$ $ $ with $ $ $1 \mbox{ } \leq i \mbox{ } \leq s$ $ $ appears, this intersection is reduced to zero and then, $ $ $\widehat{f_{2}}$ $ $ is an isomorphism which transforms each $ $ $Z_{l_{i}}$ $ $ in $ $ $z_{l_{i}}$. $ $ As the ordered monomials in $ $ $Z_{l_{1}}, \mbox{ } ... \mbox{ } , \mbox{ } Z_{l_{e}}$ $ $ are linearly independent, we have the same property for the ordered monomials in $ $ $z_{l_{1}}, \mbox{ } ... \mbox{ } , \mbox{ } z_{l_{e}}$. $ $ This easily implies that $ $ $z_{l_{1}}, \mbox{ } ... \mbox{ } , \mbox{ } z_{l_{e}}$ $ $ are nonzero and that the Laurent ordered monomials in $ $ $z_{l_{1}}, \mbox{ } ... \mbox{ } , \mbox{ } z_{l_{e}}$ $ $ are also linearly independent.
\end{enumerate}
\bx

\subsection{Each $ $ $x_{i}^{(m)}$ $ $ is a Laurent polynomial in $ $ $z_{l_{1}}, \mbox{ } ... \mbox{ } , \mbox{ } z_{l_{e}}$.}

The conventions are the same as in section 4.2. and we still assume that $ $ $\overline{\Delta}$ $ $ is nonempty. \\ $ $ \\

Let us consider some $ $ $i$ $ $ in $ $ $\llbracket  1,$ ... , $t \rrbracket $. \\ $ $ \\

If $ $ $u$ $\in$ $D_{m}$ $ $ and if $ $ $\gamma$ $ $ $\in$ $ $ \bbZ$\Pi$, $ $ we say (as in section 3.2) that $ $ $u$  $ $ \it is homogeneous of degree $ $ $\gamma$ \rm $ $ if $ $  $\overline{h_{\rho}}$($u$) = $q^{-(\rho,\gamma)}u$ $ $ for all $ $ $\rho$ $ $ in $ $ \bbZ$\Pi$. Since $ $ $q$ $ $ is not a root of unity and $ $ $V$ $ $ is spanned by $ $ $\Pi$, $ $ the degree of a nonzero homogeneous element is uniquely defined. Moreover, we immediately have the following properties: \\ $ $ \\

\begin{itemize}
\item If $ $ $u_{1}, \mbox{ } ... \mbox{ } , \mbox{ } u_{r}$ $ $ are homogeneous of same degree $ $ $\gamma$, $ $ then any linear combination of $ $ $u_{1}, \mbox{ } ... \mbox{ } , \mbox{ } u_{r}$ $ $ (with coefficients in $ $ $k$) $ $ is homogeneous of degree $ $ $\gamma$.
\item If $ $ $u_{1}$  $ $ is homogeneous of degree $ $ $\gamma_{1}$ $ $ and $ $ $u_{2}$  $ $  is homogeneous of degree $ $ $\gamma_{2}$, $ $ then $ $ $u_{1} u_{2}$  $ $ is homogeneous of degree $ $ $\gamma_{1}$ +  $\gamma_{2}$. 
\item Likewise, if $ $ $u$ $ $ is nonzero and homogeneous of degree $ $ $\gamma$, $ $ then $ $ $u^{-1}$ $ $ is homogeneous of degree $ $ $-$ $\gamma$.
\item So, if $ $ $u_{1}, \mbox{ } ... \mbox{ } , \mbox{ } u_{r}$ $ $ are nonzero and homogeneous of degrees $ $ $\gamma_{1}, \mbox{ } ... \mbox{ } , \mbox{ } \gamma_{r}$ $ $ respectively, if $ $ $\underline{a}$ = $(a_{1}, \mbox{ } ... \mbox{ } , \mbox{ } a_{r})$ $\in$ $ $ \bbZ$^{r}$,$ $ then $ $ $u^{\underline{a}}$ := $u_{1}^{a_{1}} \mbox{ } ... \mbox{ } u_{r}^{a_{r}}$ $ $ is homogeneous of degree $ $ $a_{1} \gamma_{1} + \mbox{ } ... \mbox{ } + a_{r} \gamma_{r}$.
\end{itemize}

\newpage

\newtheorem{lem4.3.}{Lemma 4.3.}

\begin{lem4.3.} $ $ \\
Consider some $ $ $j$ $ $ in $ $ $\llbracket  2,$ ... , $m \rrbracket $. \\
\begin{enumerate}
\item If $ $ $j$ $ $ $\leq$ $ $ $i + 1$, $ $ then $ $ $x_{i}^{(j)}$ $ $ = $ $ $z_{i}$. $ $ In particular, $x_{i}^{(i+1)}$ $ $ = $ $ $z_{i}$.
\item If $ $ $1$ $\leq$ $i$ $<$ $l$ $\leq$ $t$, $ $ then $ $ $x_{l}^{(j)}x_{i}^{(j)}$ $-$ $\lambda_{l,i}x_{i}^{(j)}x_{l}^{(j)}$ = $ $ $d_{l}^{(j)}(x_{i}^{(j)})$ = $ $ $f_{j}(P_{l,i}^{(j)})$ $\in$ $ $ $k<x_{i+1}^{(j)}, \mbox{ } ... \mbox{ } , \mbox{ } x_{l-1}^{(j)}>$. $ $ Moreover, if $ $ $ $ $ $ $l$ $\geq$ $j$ or $ $ $ $ $l$ = $i + 1$, $ $ then $ $ $d_{l}^{(j)}(x_{i}^{(j)})$ = $0$.
\item If $ $ $U$ $\in$ $R^{(j)}$ is homogeneous of degree $ $ $\gamma$, $ $ then $ $ $f_{j}(U)$ $ $ is homogeneous of same degree $ $ $\gamma$. $ $ In particular, $ $ $x_{i}^{(j)}$ $ $ is homogeneous of degree $ $ $\beta_{i}$.
\item If $ $ $j$ $\leq$ $l$ $\leq$ $t$ $ $ and $ $ $u$ $\in$ $k<x_{1}^{(j)}, \mbox{ } ... \mbox{ } , \mbox{ } x_{l-1}^{(j)}>$, $ $ then  $ $ $z_{l} u$ = $\overline{h_{l}}(u) z_{l}$.
\end{enumerate}
\end{lem4.3.}

\bf Proof \rm \\

\begin{enumerate}
\item By proposition 4.2. 1, we have $ $ $x_{i}^{(i+1)}$ $ $ = $ $ $x_{i}^{(i)}$ $ $ = $ $ $...$ $ $ = $ $ $x_{i}^{(2)}$ $ $ = $ $ $z_{i}$.
\item By proposition 4.2. 1, we have $ $ $x_{l}^{(j)}x_{i}^{(j)}$ $-$ $\lambda_{l,i}x_{i}^{(j)}x_{l}^{(j)}$ =  $ $ $x_{l}^{(j)}x_{i}^{(j)}$ $-$ $\overline{h_{l}}(x_{i}^{(j)})x_{l}^{(j)}$ =  $ $ $d_{l}^{(j)}$($x_{i}^{(j)}$) = $ $ $d_{l}^{(j)} \circ \mbox{ } f_{j}$($X_{i}^{(j)}$) = $ $ $f_{j} \circ \mbox{ } \delta_{l}^{(j)}$($X_{i}^{(j)}$) $ $ and we know (see section 3.2) that $ $ $\delta_{l}^{(j)}$($X_{i}^{(j)}$) =  $P_{l,i}^{(j)}$ $ $ $\in$ $ $ $k<X_{i+1}^{(j)}, \mbox{ } ... \mbox{ } , \mbox{ } X_{l-1}^{(j)}>$. $ $ Moreover, if $ $ $l$ $\geq$ $j$ or $ $ $l$ = $i + 1$, $ $ we know (see section 3.2) that $ $ $P_{l,i}^{(j)}$ = $0$. Transforming by $ $ $f_{j}$, $ $ we get the required result.
\item Consider $ $ $\rho$ $ $ in $ $ \bbZ$\Pi$. $ $ By proposition 4.2. 1, we have $ $  $\overline{h_{\rho}}$($f_{j}(U)$) = $f_{j}$($h_{\rho}(U)$) = $q^{-(\rho,\gamma)}f_{j}(U)$.
\item By 1., we have $ $ $z_{l}$ = $x_{l}^{(j)}$. $ $ So, $ $ $z_{l}u$ $-$ $\overline{h_{l}}(u)z_{l}$ = $ $ $x_{l}^{(j)}u$ $-$ $\overline{h_{l}}(u)x_{l}^{(j)}$ = $d_{l}^{(j)}$($u$). \\ By 2., if $ $ $i$ $<$ $l$, $ $ we have $ $ $d_{l}^{(j)}$($x_{i}^{(j)}$) = $0$ $ $ (since $ $ $j$ $\leq$ $l$). $ $ This implies that $ $ $d_{l}^{(j)}$ $ $ is zero on $ $ $k<x_{1}^{(j)}, \mbox{ } ... \mbox{ } , \mbox{ } x_{l-1}^{(j)}>$ $ $ and, in particular, $ $ $d_{l}^{(j)}$($u$) = $0$.
\end{enumerate}
  \bx

$ $ \\

\begin{lem4.3.} $ $  \\
Consider some $ $ $j$ $ $ in $ $ $\llbracket  2,$ ... , $m \rrbracket $. \\
\begin{enumerate}
\item If $ $ $j \mbox{ } \leq \mbox{ } l_{1}$, $ $ then $ $ $x_{i}^{(j)}$ $ $ = $ $ $z_{i}$.
\item Assume that $ $ $l_{1} \mbox{ } < \mbox{ } j$ $ $ and denote by $ $ $d$ $ $ the greatest integer such that  $ $ $l_{d} \mbox{ } < \mbox{ } j$ $ $ $(1 \mbox{ } \leq \mbox{ } d \mbox{ } \leq \mbox{ } p)$. $ $  Then $ $ $x_{i}^{(j)}$ $ $ = $ $ $x_{i}^{(l_{d}+1)}$.
\end{enumerate}
\end{lem4.3.}

\bf Proof \rm \\
\begin{description}
\item[2. ] We prove this by induction on $ $ $j$. $ $ If $ $ $j$ $ $ = $ $ $l_{d}+1$, $ $ there is nothing to prove. Assume $ $ $j$ $ $ $>$ $ $ $l_{d}+1$ $ $ and set $ $ $j^{\prime}$ $ $ = $ $ $j-1$ $ $ $>$ $ $ $l_{d}$. $ $ By lemma 4.3. 1, we have $ $ $x_{j^{\prime}}^{(j^{\prime}+1)}$ $ $ = $ $ $z_{j^{\prime}}$ $ $ = $ $ $0$ $ $ (since $ $ $j^{\prime}$ $\in$ $\Delta$). $ $ So, by proposition 4.2. 1, $ $ $x_{i}^{(j)}$ $ $ = $ $ $x_{i}^{(j^{\prime}+1)}$ $ $ = $ $ $x_{i}^{(j^{\prime})}$ $ $ = $ $ $x_{i}^{(l_{d}+1)}$ $ $ (by the induction assumption).
\item[1. ] The proof is the same (observing that there is nothing to prove if $ $ $j$ $ $ = $ $ $2$).
\end{description}
\bx

Now, assume that $ $ $l_{1} \mbox{ } < \mbox{ } m$ $ $ and denote by $ $ $p$ $ $ the greatest integer such that $ $ $l_{p} \mbox{ } < \mbox{ } m$.

\newtheorem{prop4.3.}{Proposition 4.3.}

\begin{prop4.3.} $ $ \\
Consider some $ $ $j$ $ $ in $ $ $\llbracket  2,$ ... , $m \rrbracket $ $ $ and assume that $ $ $i \mbox{ } < \mbox{ } j$. \\
\begin{enumerate}
\item If $ $ $i \mbox{ } < \mbox{ } j \mbox{ } \leq \mbox{ } l_{1}$, $ $ then $ $ $x_{i}^{(j)}$ $ $ = $ $ $z_{i}$ = $0$. \\ $ $ \\
Assume that $ $ $l_{1} \mbox{ } < \mbox{ } j$ $ $ and denote by $ $ $d$ $ $ the greatest integer such that  $ $ $l_{d} \mbox{ } < \mbox{ } j$ $ $ $(1 \mbox{ } \leq \mbox{ } d \mbox{ } \leq \mbox{ } p)$.

\newpage

\item 
\begin{enumerate}
\item If $ $ $l_{d} \mbox{ } \leq \mbox{ } i \mbox{ } < \mbox{ } j$, $ $ then $ $ $x_{i}^{(j)}$ $ $ = $ $ $z_{i}$.
\item If $ $ $i \mbox{ } < \mbox{ } l_{d} \mbox{ } < \mbox{ } j$, $ $ then
$$x_{i}^{(j)} \mbox{ } = \mbox{ } x_{i}^{(l_{d})} \mbox{ } + \mbox{ } Q_{1}z_{l_{d}}^{-1} + \mbox{ } ... \mbox{ }+ Q_{K}z_{l_{d}}^{-K}$$
with $ $ $K$ $\geq$ $1$ and: 
\begin{itemize}
\item If $ $ $d$ = $1$, $ $  then each $ $ $Q_{l}$ $ $ $\in$ $ $ $k$.
\item If $ $ $d$ $>$ $1$ and $ $ $i$ $\geq$ $l_{d-1}$, $ $ then each $ $ $Q_{l}$ $ $ $\in$ $ $ $k$.
\item If $ $ $d$ $>$ $1$ and $ $ $i$ $<$ $l_{d-1}$, $ $ then each $ $ $Q_{l}$ $ $ $\in$ $ $ $k<x_{i+1}^{(l_{d})}, \mbox{ } ... \mbox{ } , \mbox{ } x_{l_{d-1}}^{(l_{d})}>$.
\end{itemize}
\end{enumerate}
\end{enumerate}
\end{prop4.3.}

\bf Proof \rm \\ $ $ \\
\begin{enumerate}
\item By (lemma 4.3. 2, 1.), we have $ $ $x_{i}^{(j)}$ $ $ = $ $ $z_{i}$ $ $ and, since $ $ $i \mbox{ } < \mbox{ } l_{1}$, $ $ we have $ $ $i \mbox{ } \in \mbox{ } \Delta$, $ $ so that $ $ $z_{i}$ = $0$.
\item
\begin{enumerate}
\item By (lemma 4.3. 2, 2.), we have $ $ $x_{i}^{(j)}$ $ $ = $ $ $x_{i}^{(l_{d}+1)}$ $ $ and, since $ $ $l_{d}+1$ $\leq$ $i$, $ $ this is also equal to $ $ $z_{i}$.
\item We prove this result by decreasing induction on $ $ $i$. \\
As in 2.a., we have $ $ $x_{i}^{(j)}$ $ $ = $ $ $x_{i}^{(l_{d}+1)}$ $ $ and, since $ $ $x_{l_{d}}^{(l_{d}+1)}$ = $z_{l_{d}}$ $\neq$ $0$, $ $ it results from (proposition 4.2. 1, 7.(b)) that
$$x_{i}^{(j)} \mbox{ } = \mbox{ } x_{i}^{(l_{d})} + P_{1}z_{l_{d}}^{-1} + \mbox{ } ... \mbox{ }+ P_{L}z_{l_{d}}^{-L}$$
with $ $ $L$ $\geq$ $1$ $ $ and each $ $ $P_{l}$ $ $ = $ $ $\mu_{l}$ $(d_{l_{d}}^{(l_{d}+1)})^{l}$ $(x_{i}^{(l_{d}+1)})$ $ $ ($\mu_{l}$ $\in$ $k^{*}$). $ $ So, by lemma 4.3. 1,  each $ $ $P_{l}$ $ $ $\in$ $ $ $k<x_{i+1}^{(l_{d}+1)}, \mbox{ } ... \mbox{ }, \mbox{ } x_{l_{d}-1}^{(l_{d}+1)}>$ $ $ = $ $ $k<x_{i+1}^{(j)}, \mbox{ } ... \mbox{ }, \mbox{ } x_{l_{d}-1}^{(j)}>$ $ $ (by lemma 4.3. 2). \\ $ $ \\

Assume that $ $ $i$ = $l_{d}-1$. $ $ Then, by lemma 4.3. 1, we have $ $ $d_{l_{d}}^{(l_{d}+1)} \mbox{ } (x_{i}^{(l_{d}+1)})$ $ $ = $0$. $ $ This implies that $ $ $P_{l}$ = $0$ $ $ for each $ $ $l$. $ $ So, $ $ $x_{i}^{(j)}$ $ $ = $ $ $x_{i}^{(l_{d})}$ $ $ and the proof is over.\\ $ $ \\

Now, assume that $ $ $i$ $<$ $l_{d}-1$ and, for $ $ $i+1$ $\leq$ $h$ $\leq$ $l_{d}-1$,
$$x_{h}^{(j)} \mbox{ } = \mbox{ } x_{h}^{(l_{d})} \mbox{ } + \mbox{ } Q_{1,h}z_{l_{d}}^{-1} + \mbox{ } ... \mbox{ }+ \mbox{ } Q_{M,h}z_{l_{d}}^{-M}$$
with $ $ $M$ $\geq$ $1$, and: \\
\begin{itemize}
\item each $Q_{l,h}$ $ $ $\in$ $ $ $k$ $ $ if $ $ $d$ = $1$,
\item $Q_{l,h}$ $ $ $\in$ $ $ $k$ $ $ if $ $ $d$ $>$ $1$ $ $ and $ $ $h$ $\geq$ $l_{d-1}$, 
\item $Q_{l,h}$ $ $ $\in$ $ $ $k<x_{h+1}^{(l_{d})}, \mbox{ } ... \mbox{ } x_{l_{d-1}}^{(l_{d})}>$ $ $ if $ $ $d$ $>$ $1$ and $ $ $h$ $<$ $l_{d-1}$.
\end{itemize}
This implies the following results: \\  $ $ \\
\begin{itemize}
\item If $ $ $d$ = $1$, then each $x_{h}^{(j)} \mbox{ } \in  \mbox{ } k<x_{h}^{(l_{d})}, z_{l_{d}}^{-1}>$ $ $  ($i+1$ $\leq$ $h$ $\leq$ $l_{d}-1$).
\item If $ $ $d$ $>$ $1$ $ $ and $ $ $l_{d-1}$ $\leq$ $h$ $\leq$ $l_{d}-1$, then $ $ $x_{h}^{(j)} \mbox{ } \in  \mbox{ } k<x_{h}^{(l_{d})}, z_{l_{d}}^{-1}>$.
\item If $ $ $d$ $>$ $1$ $ $ and $ $ $i + 1$ $\leq$ $h$ $<$ $l_{d-1}$, $ $ then $x_{h}^{(j)} \mbox{ } \in  \mbox{ } k<x_{h}^{(l_{d})}, x_{h+1}^{(l_{d})}, \mbox{ } ... \mbox{ }, \mbox{ } x_{l_{d-1}}^{(l_{d})}, z_{l_{d}}^{-1}>$.
\end{itemize}

So, it turns out that each $ $ $P_{l}$ $ $ $\in$ $ $ $k<x_{i+1}^{(l_{d})}, \mbox{ } ... \mbox{ } , \mbox{ } x_{l_{d}-1}^{(l_{d})}, z_{l_{d}}^{-1}>$. \\ $ $ \\
\begin{itemize}
\item If $ $ $d$ = $1$ $ $ and $ $ $h$ $<$ $l_{d}$, $ $ we have (by the first point of the proposition) $ $ $x_{h}^{(l_{d})}$ =  $z_{h}$ = $0$. $ $ So, in this case, each $ $ $P_{l}$ $ $ $\in$ $ $ $k<z_{l_{d}}^{-1}>$ $ $ and
$$x_{i}^{(j)} \mbox{ } = \mbox{ } x_{i}^{(l_{d})} \mbox{ } + \mbox{ } Q_{1}z_{l_{d}}^{-1} + \mbox{ } ... \mbox{ }+ Q_{K}z_{l_{d}}^{-K}$$
with $ $ $K$ $\geq$ $1$ $ $ and each $ $ $Q_{l}$ $ $ in $ $ $k$.
\item Now, assume that $ $ $d$ $>$ $1$.\\
If $ $ $l_{d-1}$ $<$ $h$ $\leq$ $l_{d}-1$, $ $ we have (by 2.a.) $ $ $x_{h}^{(l_{d})}$ = $z_{h}$ = $0$ $ $ (since $ $ $h$ $\in$ $\Delta$). So,
\begin{itemize}
\item If $ $ $i$ $\geq$ $l_{d-1}$, each $ $ $P_{l}$ $ $ $\in$ $ $ $k<z_{l_{d}}^{-1}>$ $ $ and
$$x_{i}^{(j)} \mbox{ } = \mbox{ } x_{i}^{(l_{d})} \mbox{ } + \mbox{ } Q_{1}z_{l_{d}}^{-1} + \mbox{ } ... \mbox{ }+ Q_{K}z_{l_{d}}^{-K}$$
with $ $ $K$ $\geq$ $1$ $ $ and each $ $ $Q_{l}$ $ $ in $ $ $k$.
\item If $ $ $i$ $<$ $l_{d-1}$,
each $ $ $P_{l}$ $ $ $\in$ $ $ $k<x_{i+1}^{(l_{d})}, \mbox{ } ... \mbox{ } , \mbox{ }  x_{l_{d-1}}^{(l_{d})}, z_{l_{d}}^{-1}>$. $ $ $ $ As, by lemma 4.3. 1, $ $ $z_{l_{d}}$ = $x_{l_{d}}^{(l_{d})}$ $ $ $q$ - commutes with $ $ $x_{i+1}^{(l_{d})}, \mbox{ } ... \mbox{ } , \mbox{ } x_{l_{d-1}}^{(l_{d})}$, $ $ each $ $ $P_{l}$ $ $ can be written as follows:
$$P_{l} \mbox{ } = \mbox{ } S_{0,l} +  S_{1,l}z_{l_{d}}^{-1} + \mbox{ } ... \mbox{ }+ S_{E,l}z_{l_{d}}^{-E}$$
with $ $ $E$ $\geq$ $1$ $ $ and each $ $ $S_{j,l}$ $ $ $\in$ $ $ $k<x_{i+1}^{(l_{d})}, \mbox{ } ... \mbox{ } , \mbox{ } x_{l_{d-1}}^{(l_{d})}>.$ $ $ So, we conclude that, in this case,
$$x_{i}^{(j)} \mbox{ } = \mbox{ } x_{i}^{(l_{d})} \mbox{ } + \mbox{ } Q_{1}z_{l_{d}}^{-1} + \mbox{ } ... \mbox{ }+ Q_{K}z_{l_{d}}^{-K}$$
with $ $ $K$ $\geq$ $1$ $ $ and each $ $ $Q_{l}$ $ $ in $ $ $k<x_{i+1}^{(l_{d})}, \mbox{ } ... \mbox{ } , \mbox{ } x_{l_{d-1}}^{(l_{d})}>$.
\end{itemize}
\end{itemize}
\end{enumerate}
\end{enumerate}
\bx

\begin{prop4.3.} $ $ \\
Consider some $ $ $j$ $ $ in $ $ $\llbracket  2,$ ... , $m \rrbracket $, $ $ assume that $ $ $i+1 \mbox{ } < \mbox{ } j$ $ $ and consider some $ $ $u$ $\in$ $ $ $k<x_{i+1}^{(j)}, \mbox{ } ... \mbox{ } , \mbox{ } x_{j-1}^{(j)}>.$ \\ $ $ \\
\begin{enumerate}
\item If $ $ $j \mbox{ } \leq \mbox{ } l_{1}$, $ $ then $ $ $u$ $\in$ $ $ $k$. \\
Assume that $ $ $l_{1} \mbox{ } < \mbox{ } j$ $ $ and denote by $ $ $d$ $ $ the greatest integer such that  $ $ $l_{d} \mbox{ } < \mbox{ } j$ $ $ $(1 \mbox{ } \leq \mbox{ } d \mbox{ } \leq \mbox{ } p)$. \\
\item
\begin{enumerate}
\item If $ $ $l_{d} \mbox{ } \leq \mbox{ } i$, $ $ then $ $ $u$ $\in$ $ $ $k$.
\item If $ $ $i \mbox{ } < \mbox{ } l_{d}$, $ $ then
$$u \mbox{ } = \mbox{ } u_{1}z_{l_{d}}^{a_{1}} + \mbox{ } ... \mbox{ }+ u_{K}z_{l_{d}}^{a_{K}}$$
with $ $ $K$ $\geq$ $1$, $ $ $(a_{1}, \mbox{ } ... \mbox{ }, \mbox{ } a_{K})$ $\in$ \bbZ$^{K}$ $ $ and: \\
\begin{itemize}
\item If $ $ $d$ = $1$, $ $  then each $ $ $u_{l}$ $ $ $\in$ $ $ $k$.
\item If $ $ $d$ $>$ $1$ and $ $ $i$ $\geq$ $l_{d-1}$, $ $ then each $ $ $u_{l}$ $ $ $\in$ $ $ $k$.
\item If $ $ $d$ $>$ $1$ and $ $ $i$ $<$ $l_{d-1}$, $ $ then each $ $ $u_{l}$ $ $ $\in$ $ $ $k<x_{i+1}^{(l_{d})}, \mbox{ } ... \mbox{ } , \mbox{ } x_{l_{d-1}}^{(l_{d})}>$.
\end{itemize}  
\end{enumerate}
\end{enumerate}
\end{prop4.3.}

\bf Proof \rm \\ $ $ \\
\begin{enumerate}
\item By (proposition 4.3. 1, 1.) we have $ $ $x_{h}^{(j)}$ = $0$ $ $ for any $ $ $h$ $<$ $j$. $ $ So, $ $ $u$ $\in$ $ $ $k$.
\item
\begin{enumerate}
\item Consider some integer $ $ $h$ $ $ such that $ $ $l_{d} \mbox{ } < \mbox{ } h \mbox{ } < \mbox{ } j$. $ $ By (proposition 4.3. 1, 2.a.) we have $ $ $x_{h}^{(j)}$ = $z_{h}$ = $0$ $ $ (since $ $ $h$ $\in$ $\Delta$). $ $ So, $ $ $u$ $\in$ $ $ $k$.
\item Consider some integer $ $ $h$ $ $ with $ $ $i$ $<$ $h$ $<$ $j$. \\
If $ $ $l_{d} \mbox{ } < \mbox{ } h$, $ $ then, as in 2.a., we have $ $ $x_{h}^{(j)}$ = $0$. $ $ If $ $ $h$ = $l_{d}$, $ $ then, again by (proposition 4.3. 1, 2.a.), we have $ $ $x_{h}^{(j)}$ = $z_{h}$ = $z_{l_{d}}$. \\ $ $ \\
\begin{itemize}
\item Assume that $ $ $d$ = $1$. $ $ If $ $ $h$ $<$ $l_{d}$, $ $ then, by (proposition 4.3. 1, 2.b.), we have
$$x_{h}^{(j)} \mbox{ } = \mbox{ } x_{h}^{(l_{d})} \mbox{ } + \mbox{ } Q_{1}z_{l_{d}}^{-1} + \mbox{ } ... \mbox{ }+ Q_{M}z_{l_{d}}^{-M}$$
with $ $ $M$ $\geq$ $1$ $ $ and each $ $ $Q_{l}$ $ $ in $ $ $k$. $ $ By (proposition 4.3. 1, 1.), we also have $ $ $x_{h}^{(l_{d})}$ = $0$. $ $ So, $ $ $x_{h}^{(j)} \mbox{ } \in \mbox{ } k<z_{l_{d}}^{-1}>$ $ $ and it turns out that $ $ $u \mbox{ } \in \mbox{ } k<z_{l_{d}}^{\pm 1}>$. $ $ This implies that
$$u \mbox{ } = \mbox{ } u_{1}z_{l_{d}}^{a_{1}} + \mbox{ } ... \mbox{ }+ u_{K}z_{l_{d}}^{a_{K}}$$
with $ $ $N$ $\geq$ $1$, $ $ ($a_{1}, \mbox{ } ... \mbox{ }, \mbox{ } a_{N}$) $\in$ \bbZ$^{N}$ $ $ and each $ $ $u_{l}$ $ $ $\in$ $ $ $k$.
\item Assume that $ $ $d$ $>$ $1$ $ $ and $ $ $i$ $\geq$ $l_{d-1}$. $ $ If $ $ $i$ $<$ $h$ $<$ $l_{d}$, $ $ we have $ $ $l_{d-1}$ $<$ $h$ $<$ $l_{d}$ $ $ and, by (proposition 4.3. 1, 2.b.),
$$x_{h}^{(j)} \mbox{ } = \mbox{ } x_{h}^{(l_{d})} \mbox{ } + \mbox{ } Q_{1}z_{l_{d}}^{-1} + \mbox{ } ... \mbox{ }+ Q_{M}z_{l_{d}}^{-M}$$
with $ $ $M$ $\geq$ $1$ $ $ and each $ $ $Q_{l}$ $ $ $\in$ $ $ $k$. $ $ By (proposition 4.3. 1, 1.) we also have $ $ $x_{h}^{(l_{d})}$ = $z_{h}$ $ $ and, since $ $ $h$ $\in$ $\Delta$, $ $ $x_{h}^{(l_{d})}$ = $0$. $ $ As over, we conclude that
$$u \mbox{ } = \mbox{ } u_{1}z_{l_{d}}^{a_{1}} + \mbox{ } ... \mbox{ }+ u_{K}z_{l_{d}}^{a_{K}}$$
with $ $ $K$ $\geq$ $1$, $ $ ($a_{1}, \mbox{ } ... \mbox{ }, \mbox{ } a_{K}$) $\in$ \bbZ$^{K}$ $ $ and each $ $ $u_{l}$ $ $ $\in$ $ $ $k$.
\item Assume that $ $ $d$ $>$ $1$ and $ $ $i$ $<$ $l_{d-1}$. $ $ If $ $ $i$ $<$ $h$ $<$ $l_{d}$, $ $ we have, always by (proposition 4.3. 1, 2.b.),
$$x_{h}^{(j)} \mbox{ } = \mbox{ } x_{h}^{(l_{d})} \mbox{ } + \mbox{ } Q_{1}z_{l_{d}}^{-1} + \mbox{ } ... \mbox{ }+ Q_{M}z_{l_{d}}^{-M}$$
with $ $ $M$ $\geq$ $1$ $ $ and: \\
\begin{itemize}
\item If $ $ $l_{d-1}$ $<$ $h$, $ $ then each $ $ $Q_{l}$ $ $ $\in$ $ $ $k$ $ $ and, as in the previous case, $ $ $x_{h}^{(l_{d})}$ = $0$. $ $ So, $ $ $x_{h}^{(j)} \mbox{ } \in \mbox{ } k<z_{l_{d}}^{-1}>$.
\item If $ $ $h$ $=$ $l_{d-1}$, $ $ then each $ $ $Q_{l}$ $ $ $\in$ $ $ $k$ $ $ and, consequently, $ $ $x_{h}^{(j)} \mbox{ } \in \mbox{ } k<x_{h}^{(l_{d})}, z_{l_{d}}^{-1}> \mbox{ } = \mbox{ } k<x_{l_{d-1}}^{(l_{d})}, z_{l_{d}}^{-1}>$.
\item If $ $ $h$ $<$ $l_{d-1}$, $ $ then each $ $ $Q_{l}$ $ $ $\in$ $ $ $<x_{h+1}^{(l_{d})}, \mbox{ } ... \mbox{ } , \mbox{ } x_{l_{d-1}}^{(l_{d})}>$ $ $ and, consequently, $ $ $x_{h}^{(j)} \mbox{ } \in \mbox{ } k<x_{h}^{(l_{d})}, x_{h+1}^{(l_{d})}, \mbox{ } ... \mbox{ } , \mbox{ } x_{l_{d-1}}^{(l_{d})}, z_{l_{d}}^{-1}>$. \\
So, it turns out that $ $ $u \mbox{ } \in \mbox{ } k<x_{i+1}^{(l_{d})}, \mbox{ } ... \mbox{ } , \mbox{ } x_{l_{d-1}}^{(l_{d})}, z_{l_{d}}^{\pm 1}>$. $ $ As, by lemma 4.3. 1, $ $ $z_{l_{d}}$ = $x_{l_{d}}^{(l_{d})}$ $ $ $q$ - commutes with $ $ $x_{i+1}^{(l_{d})}, \mbox{ } ... \mbox{ } , \mbox{ } x_{l_{d-1}}^{(l_{d})}$, $ $ we conclude that
$$u \mbox{ } = \mbox{ } u_{1}z_{l_{d}}^{a_{1}} + \mbox{ } ... \mbox{ }+ u_{K}z_{l_{d}}^{a_{K}}$$
with $ $ $K$ $\geq$ $1$, $ $ ($a_{1}, \mbox{ } ... \mbox{ }, \mbox{ } a_{K}$) $\in$ \bbZ$^{K}$ $ $ and each $ $ $u_{l}$ $ $ in $ $ $k<x_{i+1}^{(l_{d})}, \mbox{ } ... \mbox{ } , \mbox{ } x_{l_{d-1}}^{(l_{d})}>$.
\end{itemize}
\end{itemize}
\end{enumerate}
\end{enumerate}
\bx

\begin{prop4.3.} $ $ \\
Consider some $ $ $j$ $ $ in $ $ $\llbracket  2,$ ... , $m \rrbracket $ $ $ with $ $ $i$ $<$ $j$. \\
\begin{enumerate}
\item If $ $ $\{l_{1}, \mbox{ } ... \mbox{ } ,\mbox{ } l_{p}\} \mbox{ } \cap \mbox{ } \llbracket  i+1,\mbox{ } ... \mbox{ } , \mbox{ } j-1 \rrbracket  \mbox{ } = \mbox{ } \emptyset$, $ $ then  $ $ $x_{i}^{(j)} \mbox{ } = \mbox{ } z_{i}$.
\item Assume that $ $ $\{l_{1}, \mbox{ } ... \mbox{ } ,\mbox{ } l_{p}\}  \cap  \llbracket  i+1,\mbox{ } ... \mbox{ } , \mbox{ } j-1 \rrbracket  \mbox{ } \neq \mbox{ } \emptyset$ $ $ and set $ $ $\{l_{1}, \mbox{ } ... \mbox{ } , \mbox{ } l_{p}\}  \cap  \llbracket  i+1,$ ... , $j-1 \rrbracket  \mbox{ } = \mbox{ } \{l_{c} \mbox{ } <  ... < \mbox{ } l_{d}\}$. $ $ Then, we have
$$x_{i}^{(j)} \mbox{ } = \mbox{ } z_{i} \mbox{ } + \mbox{ } \displaystyle \sum_{\underline{a} \mbox{ } = \mbox{ } (a_{c}, \mbox{ }  ... \mbox{ } , \mbox{ } a_{d}) \mbox{ } \in \mbox{ } F}{} \mbox{ } \eta(\underline{a}) z_{l_{c}}^{a_{c}} \mbox{ }  ... \mbox{ } z_{l_{d}}^{a_{d}},$$
where \\
\begin{itemize}
\item $F$ $ $ is a finite (possibly empty) subset of $ $ \bbZ$^{d-c+1}$.
\item If $ $ $\preceq$ $ $ denotes the inverse lexicographic order in $ $ \bbZ$^{d-c+1}$, $ $ then, for any $ $ $\underline{a}$ $ $ in $ $ $F$, $ $ we have $ $ $\underline{a}$ $ $ $\prec$ $ $ $0$ $ $ $($ie. $ $ $\underline{a} = (a_{c}, \mbox{ }  ... \mbox{ } , a_{r}, 0, \mbox{ }  ... \mbox{ } , 0)$ $ $ with $ $ $a_{r} < 0$.$)$
\item For each $ $ $\underline{a}$ $ $ in $ $ $F$, $ $ we have $ $ $\eta(\underline{a}) \mbox{ } \in \mbox{ } k^{*}$.
\end{itemize}
\end{enumerate}
\end{prop4.3.}

\bf Proof \rm \\
We prove the proposition by induction on $ $ $j - i$. $ $ In order to do this, we first prove the assertion 1. \\ $ $ \\

So, assume that $ $ $\{l_{1}, \mbox{ } ... \mbox{ } ,\mbox{ } l_{p}\} \mbox{ } \cap \mbox{ } \llbracket  i+1,\mbox{ } ... \mbox{ } , \mbox{ } j-1 \rrbracket  \mbox{ } = \mbox{ } \emptyset$, $ $ and observe that there are two possibilities: \\
\begin{enumerate}
\item \begin{enumerate}
\item $j$ $\leq$ $l_{1}$. $ $ By proposition 4.3. 1, we know that $ $ $x_{i}^{(j)} \mbox{ } = \mbox{ } z_{i}$.
\item $l_{1}$ $<$ $j$ $ $ and, if $ $ $d$ $ $ is the greatest integer such that $ $ $l_{d}$ $<$ $j$,
$ $ $l_{d} \leq  i < j$. $ $ Again, by proposition 4.3. 1, we know that $ $ $x_{i}^{(j)} \mbox{ } = \mbox{ } z_{i}$. \\ $ $ \\

Assume $ $ $j - i$ $ $ = $ $ $1$. $ $ In this case, we have $ $ $\{l_{1}, \mbox{ } ... \mbox{ } ,\mbox{ } l_{p}\} \mbox{ } \cap \mbox{ } \llbracket  i+1,\mbox{ } ... \mbox{ } , \mbox{ } j-1 \rrbracket  \mbox{ } = \mbox{ } \emptyset$. $ $ So, we are in the case of assertion 1., and the proof is over.
\end{enumerate}
\end{enumerate}

Now assume $ $ $j - i$ $ $ $>$ $ $ $1$. $ $ As the assertion 1. is already proved, we may assume that $ $ $\{l_{1}, \mbox{ } ... \mbox{ } , \mbox{ } l_{p}\} \mbox{ } \cap \mbox{ } \llbracket  i+1, ... , j-1 \rrbracket  \mbox{ } = \mbox{ } \{l_{c} \mbox{ } <  ... < \mbox{ } l_{d}\}$ $ $ is non empty. So, we have $ $ $l_{1}$ $<$ $j$, $ $ $d$ $ $ is the greatest integer such that $ $ $l_{d}$ $<$ $j$ $ $ and $ $ $i$ $<$ $l_{d}$ $<$ $j$. $ $ Now, by proposition 4.3. 1, we can write
$$x_{i}^{(j)} \mbox{ } = \mbox{ } x_{i}^{(l_{d})} + Q_{1}z_{l_{d}}^{-1} + \mbox{ } ... \mbox{ }+ Q_{K}z_{l_{d}}^{-K}$$
with $ $ $K$ $\geq$ $1$ and three possible cases: \\
\begin{enumerate}
\item $d$ = $1$. Then each $ $ $Q_{l}$ $\in$ $k$ and, since  $ $ $x_{i}^{(l_{d})}$ = $z_{i}$ $ $ (see proposition 4.3. 1, 1.), the proof is over.
\item  $d$ $>$ $1$ $ $ and $ $ $i$ $\geq$ $l_{d-1}$. Then each $ $ $Q_{l}$ $\in$ $k$ $ $ and, since  $ $ $x_{i}^{(l_{d})}$ = $z_{i}$ $ $ (see proposition 4.3. 1, 2.a.), the proof is over.
\item $d$ $>$ $1$ $ $ and $ $ $i$ $<$ $l_{d-1}$ $ $ (so that $ $ $c$ $\leq$ $d-1$). $ $ In this case, each $ $ $Q_{l}$ $\in$ $k<x_{i+1}^{(l_{d})}, \mbox{ } ... \mbox{ } , \mbox{ } x_{l_{d-1}}^{(l_{d})}>$. $ $ Since $ $ $l_{d}$ $<$ $j$, $ $ it results from the induction assumption that the proposition is true for each $ $ $x_{h}^{(l_{d})}$ $ $ with $ $ $i$ $\leq$ $h$ $<$ $l_{d}$. $ $ This implies that \\
\begin{itemize}
\item If $ $ $i+1$ $\leq$ $h$ $\leq$ $l_{d-1}$, $ $ we have $ $ $x_{h}^{(l_{d})}$ $\in$ $k<z_{l_{c}}^{\pm 1}, \mbox{ } ... \mbox{ } , \mbox{ } z_{l_{d-1}}^{\pm 1}>$ $ $ (since $ $ $\{l_{1}, \mbox{ } ... \mbox{ } , \mbox{ } l_{p}\} \mbox{ } \cap \mbox{ } \llbracket  h+1, ... , l_{d}-1 \rrbracket$  $\subseteq \mbox{ } \{l_{c} \mbox{ } <  ... < \mbox{ } l_{d-1}\}$ $ $ and $ $ $z_{h}$ $ $ is either $ $ $0$ $ $ or equal to $ $ $z_{l_{r}}$ $ $  with $ $ $c$ $\leq$ $r$ $\leq$ $d-1$). $ $ So, each $ $ $Q_{l}$ $\in$ $k<z_{l_{c}}^{\pm 1}, \mbox{ } ... \mbox{ } , \mbox{ } z_{l_{d-1}}^{\pm 1}>$.
\item $$x_{i}^{(l_{d})} \mbox{ } = \mbox{ } z_{i} \mbox{ } + \mbox{ } \displaystyle \sum_{\underline{a} \mbox{ } = \mbox{ } (a_{c}, \mbox{ }  ... \mbox{ } , \mbox{ } a_{d-1}) \mbox{ } \in \mbox{ } F_{1}}{} \mbox{ } \eta(\underline{a}) z_{l_{c}}^{a_{c}} \mbox{ }  ... \mbox{ } z_{l_{d-1}}^{a_{d-1}},$$
where \\
\begin{itemize}
\item[$\diamond$] $F_{1}$ $ $ is a finite (possibly empty) subset of $ $ \bbZ$^{d-c}$.
\item[$\diamond$] For each $ $ $\underline{a}$ $ $ in $ $ $F_{1}$, $ $ we have $ $ $\underline{a}$ $ $ $\prec$ $ $ $0$ $ $ and $ $ $\eta(\underline{a}) \mbox{ } \in \mbox{ } k^{*}$.
\end{itemize}
\end{itemize}
As each $ $ $Q_{l}$ is a summand of (ordered) Laurent monomials in $ $ $z_{l_{c}}, \mbox{ } ... \mbox{ } , \mbox{ } z_{l_{d-1}}$, it turns out that
$$x_{i}^{(j)} \mbox{ } = \mbox{ } z_{i} \mbox{ } + \mbox{ } \displaystyle \sum_{\underline{a} \mbox{ } = \mbox{ } (a_{c}, \mbox{ }  ... \mbox{ } , \mbox{ } a_{d}) \mbox{ } \in \mbox{ } F}{} \mbox{ } \eta(\underline{a}) z_{l_{c}}^{a_{c}} \mbox{ }  ... \mbox{ } z_{l_{d}}^{a_{d}},$$
where \\
\begin{itemize}
\item[$\diamond$] $F$ $ $ is a finite (possibly empty) subset of $ $ \bbZ$^{d-c+1}$.
\item[$\diamond$] For each $ $ $\underline{a}$ $ $ in $ $ $F$, $ $ we have $ $ $\underline{a}$ $ $ $\prec$ $ $ $0$ $ $ and $ $ $\eta(\underline{a}) \mbox{ } \in \mbox{ } k^{*}$. 
\end{itemize}
\end{enumerate}
\bx

This proposition implies immediately

\newtheorem{cor4.3.}{Corollary 4.3.}

\begin{cor4.3.} Consider some integer $ $ $d$ $ $ with $1$ $<$ $d$ $\leq$ $p$. $ $  
If $ $ $1$ $\leq$ $i$ $<$ $l_{d}$, $ $ then $ $
$x_{i}^{(l_{d})} \mbox{ } \in \mbox{ } k<z_{l_{1}}^{\pm 1}, \mbox{ } ... \mbox{ } , \mbox{ } z_{l_{d-1}}^{\pm 1}>.$
\end{cor4.3.}

\bf Proof \rm \\
By proposition 4.3. 3, $ $ $x_{i}^{(l_{d})} \mbox{ } = \mbox{ } z_{i}$ + $Q$ $ $ with $ $ $Q$ $ $ a (possibly $ $ $0$) Laurent polynomial in $ $ $z_{l_{1}}, \mbox{ } ... \mbox{ } , \mbox{ } z_{l_{d-1}}$. $ $ As $ $ $z_{i}$ $ $ is either $ $ $0$ $ $ or $ $ $z_{l_{r}}$ $ $ with $ $ $1$ $\leq$ $r$ $\leq$ $d-1$, $ $ the proof is over. \bx

\begin{cor4.3.} Consider some integer $ $ $d$ $ $ with $1$ $<$ $d$ $\leq$ $p$. $ $ 
In $ $ $D_{m}$ $ $ considered as a left module over $ $ $k<x_{1}^{(l_{d})}, \mbox{ } ... \mbox{ }, \mbox{ } x_{l_{d}-1}^{(l_{d})}>$, $ $ the ordered Laurent monomials $ $ $z_{l_{d}}^{a_{d}} \mbox{ } ... \mbox{ } z_{l_{p}}^{a_{p}}$ $ $ \rm(\it$a_{d}, \mbox{ } ... \mbox{ } , \mbox{ } a_{p}$ $ $ in \bbZ\rm)\it $ $ are linearly independent.
\end{cor4.3.}

\bf Proof \rm \\
By corollary 4.3. 1, we have $ $ $k<x_{1}^{(l_{d})}, \mbox{ } ... \mbox{ }, \mbox{ } x_{l_{d}-1}^{(l_{d})}> \mbox{ } \subseteq \mbox{ } k<z_{l_{1}}^{\pm 1}, \mbox{ } ... \mbox{ } , \mbox{ } z_{l_{d-1}}^{\pm 1}>$. $ $ So, the corollary results immediately from the $ $ $k$ - linear independence of the ordered Laurent monomials in $ $ $z_{l_{1}}, \mbox{ } ... \mbox{ } , \mbox{ } z_{l_{e}}$ $ $ (see proposition 4.2. 2). \bx

\subsection{A new sufficient condition for $ $ $\mathcal{P}\it^{(m)}$ $ $ to be in $ $ $Im(\phi_{m})$.}

The notations are the same as in the preceding sections, but we do not assume that $ $ $\overline{\Delta}$ $ $ is nonempty a priori.

Assume that $ $ $\mathcal{P}\it^{(m)}$ $ $ is not in $ $ $Im(\phi_{m})$ so that, by proposition 4.1. 1, $ $ $X_{m}^{(m)}$ $ $ $\in$ $ $ $\mathcal{P}\it^{(m)}$ $ $ and there is some $ $ $j$ $ $ in \\
$ $ $\llbracket  1,$ ... , $m - 1 \rrbracket $ $ $ such that
$$U \mbox{ } = \mbox{ } \Theta^{(m)}\rm(\it\delta_{m}^{(m+1)}\rm(\it X_{j}^{(m+1)}\rm))\it \mbox{ } \notin \mbox{ } \mathcal{P}\it^{(m)}.$$
Let us choose such a  $ $ $j$ $ $ maximal. So, we have
$$\Theta^{(m)}\rm(\it\delta_{m}^{(m+1)}\rm(\it X_{i}^{(m+1)}\rm))\it \mbox{ } \in \mbox{ } \mathcal{P}\it^{(m)}$$
for any $ $ $i$ $ $ in $ $ $\llbracket  j + 1,$ ... , $m - 1 \rrbracket $ $ $  and we observe that, since $ $ $U$ $ $ $\notin$ $ $ $\mathcal{P}\it^{(m)}$, $ $ we have $ $ $u$ := $f_{m} (U)$ $\neq$ $0$ $ $ in $ $ $A^{(m)}$. \\ $ $   \\

Recall that
$$\Delta  \mbox{ } = \mbox{ } \{i \in \mbox{ } \rm [\it\!|1,  \mbox{ } ...  \mbox{ } ,  \mbox{ } t|\!\rm] \it \mbox{ } | \mbox{ } Z_{i} \in \mbox{ } \mathcal{P}^{(2)} \} \mbox{ } = \mbox{ } \{j_{1} \mbox{ } < \mbox{ } ...\mbox{ } <\mbox{ } j_{s}\} \rm \mbox{ } (unless \mbox{ } this \mbox{ } set \mbox{ } is \mbox{ } empty),$$
$$\overline{\Delta} \mbox{ }  = \mbox{ } \rm [\it\!|1,  \mbox{ } ...  \mbox{ } ,  \mbox{ } t|\!\rm] \it \mbox{ } \setminus \mbox{ } \Delta \mbox{ } = \mbox{ } \{l_{1} \mbox{ } < \mbox{ } ... \mbox{ } < \mbox{ } l_{e}\} \rm \mbox{ } (unless \mbox{ } this \mbox{ } set \mbox{ } is \mbox{ } empty).$$

\newtheorem{lem4.4.}{Lemma 4.4.}

\begin{lem4.4.} $ $ \\
$\overline{\Delta}$ $ $ is nonempty, $ $ $l_{1}$ $<$ $m$, $ $ $m$ $\notin$ $\overline{\Delta}$ $ $ and, if $ $ $p$ $ $ denotes the greatest integer such that $ $ $\l_{p}$ $<$ $m$, $ $ we have $ $ $j$ $<$ $\l_{p}$ $<$ $m$.
\end{lem4.4.}

\bf Proof \rm \\
First, recall that: \\ $ $ \\
\begin{itemize}
\item  If $ $ $\overline{\Delta}$ $ $ is empty, then (observation 4.2. 1) each $ $ $x_{l}^{(m)}$ ($l$ $\in$  $[\it\!|1,  \mbox{ } ...  \mbox{ } ,  \mbox{ } t|\!\rm])$ $ $ is zero.
\item $\delta_{m}^{(m+1)}$($X_{j}^{(m+1)}$) = $P_{m,j}^{(m+1)}$. $ $ This is $ $ $0$ $ $ if $ $ $m$ = $j+1$ $ $ or a (finite) summand
$$\sum_{ \underline{a} \mbox{ } = \mbox{ } (a_{j+1}, ... ,a_{m-1})} c_{\underline{a}}(X_{j+1}^{(m+1)})^{a_{j+1}} ... (X_{m-1}^{(m+1)})^{a_{m-1}}$$
with each $ $ $\underline{a}$ $ $ in \bbN$^{ m-1-j}$ $ $ if $ $ $m$ $>$ $j+1$ $ $ and, moreover, $ $ $\underline{a}$ $ $ is nonzero when $ $ $c_{\underline{a}}$ $\neq$ $0$.
\item $\Theta^{(m)}$: $ $ $k<X_{1}^{(m+1)}, \mbox{ } ... \mbox{ } , X_{m-1}^{(m+1)}>$ $ $ $\rightarrow$ $ $ $k<X_{1}^{(m)}, \mbox{ } ... \mbox{ } , X_{m-1}^{(m)}>$ $ $ is an algebra homomorphism which transforms each $ $ $X_{l}^{(m+1)}$ $ $ in $ $ $X_{l}^{(m)}$.
\end{itemize}

Assume that $ $ $m$ = $j+1$. $ $ In this case, we have
$$ \delta_{m}^{(m+1)}(X_{j}^{(m+1)}) \mbox{ } = 0 \mbox{ } \mbox{ }  \Rightarrow \mbox{ } \mbox{ }  U \mbox{ } = \mbox{ } \Theta^{(m)} ( \delta_{m}^{(m+1)} ( X_{j}^{(m+1)} )) \mbox{ } = \mbox{ } 0 \mbox{ } \mbox{ }  \Rightarrow \mbox{ } \mbox{ }  u \mbox{ } = \mbox{ } f_{m} (U) \mbox{ } = 0$$
and we obtain a contradiction. \\ $ $ \\

So, we have $ $ $m$ $>$ $j+1$, $ $ and
$$ \delta_{m}^{(m+1)}(X_{j}^{(m+1)}) \mbox{ } = \mbox{ } \sum_{ \underline{a} \mbox{ } = \mbox{ } (a_{j+1}, ... ,a_{m-1})} c_{\underline{a}}(X_{j+1}^{(m+1)})^{a_{j+1}} ... (X_{m-1}^{(m+1)})^{a_{m-1}} \mbox{ } \mbox{ }  \Rightarrow$$
$$U \mbox{ } = \mbox{ } \Theta^{(m)}\rm (\it \delta_{m}^{(m+1)}\rm (\it X_{j}^{(m+1)}\rm ))\it \mbox{ } = \mbox{ } \sum_{ \underline{a} \mbox{ } = \mbox{ } (a_{j+1}, ... ,a_{m-1})} c_{\underline{a}}(X_{j+1}^{(m)})^{a_{j+1}} ... (X_{m-1}^{(m)})^{a_{m-1}} \mbox{ } \mbox{ }  \Rightarrow \mbox{ }$$
$$u \mbox{ } = \mbox{ } f_{m} (U) \mbox{ } = \mbox{ } \sum_{ \underline{a} \mbox{ } = \mbox{ } (a_{j+1}, ... ,a_{m-1})} c_{\underline{a}}(x_{j+1}^{(m)})^{a_{j+1}} ... (x_{m-1}^{(m)})^{a_{m-1}}.$$

Assume that $ $ $\overline{\Delta}$ $ $ is empty. Then, each $ $ $x_{l}^{(m)}$ $ $ is zero and, since $ $ $\underline{a}$ $ $ is nonzero when $ $ $c_{\underline{a}}$ $\neq$ $0$, $ $ we get that $ $ $u$ = $0$ $ $ and we still have a contradiction. So, $ $ $\overline{\Delta}$ $ $ is nonempty. \\ $ $ \\

If $ $ $l_{1}$ $\geq$ $m$, $ $ then, for each $ $ $l$ $\in$  $[\it\!|1,  \mbox{ } ...  \mbox{ } ,  \mbox{ } m - 1|\!\rm]$, $ $ we have (by lemma 4.3. 2 and, since $ $ $l$ $<$ $l_{1}$) $ $ $x_{l}^{(m)}$ = $z_{l}$ = $0$. $ $ As over, this implies that $ $ $u$ = $0$ $ $ and we still get a contradiction . So, $ $ $l_{1}$ $<$ $m$. \\ $ $ \\

Since $ $ $X_{m}^{(m)}$ $ $ $\in$ $ $ $\mathcal{P}\it^{(m)}$, $ $ we have $ $ $x_{m}^{(m)}$ = $0$. $ $ As $ $ $x_{m}^{(m)}$ = $z_{m}$ $ $ (lemma 4.3. 1), we have $ $ $Z_{m}$ $ $ $\in$ $ $ $\mathcal{P}^{(2)}$ $ $ and so,  $ $ $m$ $ $ $\in$ $ $ $\Delta$. $ $ Moreover, if $ $ $l_{p}$ $\leq$ $j$, then, for each $ $ $l$ $\in$  $[\it\!|j+1,  \mbox{ } ...  \mbox{ } ,  \mbox{ } m - 1|\!\rm]$, $ $ we have (by propostion 4.3. 1, 2.a. and, since $ $ $l_{p}$ $<$ $l$ $<$ $m$) $ $ $x_{l}^{(m)}$ = $z_{l}$ = $0$. $ $ As over, this implies that $ $ $u$ = $0$ $ $ and we still get a contradiction . So, $ $ $j$ $<$ $l_{p}$. \bx

Now, set $ $ $V$ = $\delta_{m}^{(m+1)}$($X_{j}^{(m+1)}$) = $X_{m}^{(m+1)}X_{j}^{(m+1)}$ $-$ $q^{-(\beta_{m}, \beta_{j} )}X_{j}^{(m+1)}X_{m}^{(m+1)}$ $ $ (see section 3.2), so that $ $ $U$ = $ \Theta^{(m)} (V)$.

\begin{lem4.4.} $ $ \\
Consider some integer $ $ $l$ $ $ with $ $ $j$ $<$ $l$ $<$ $m$. $ $ Then we have, in $ $ $R^{(m+1)}$,
$$X_{l}^{(m+1)}V \mbox{ } = \mbox{ } q^{(\beta_{l}, \beta_{m}-\beta_{j})} V X_{l}^{(m+1)}$$
$$ + \mbox{ }  q^{(\beta_{l}, \beta_{m})} \mbox{ } [ \mbox{ } q^{-(\beta_{j}, \beta_{l}+\beta_{m})} X_{j}^{(m+1)} \delta_{m}^{(m+1)} ( X_{l}^{(m+1)} )  \mbox{ } - \mbox{ } \delta_{m}^{(m+1)} (X_{l}^{(m+1)} ) X_{j}^{(m+1)} \mbox{ } ]$$
$$ + \mbox{ } q^{(\beta_{l}, \beta_{m})} \delta_{m}^{(m+1)} \circ \delta_{l}^{(m+1)}  ( X_{j}^{(m+1)} ).$$
\end{lem4.4.}

\bf Proof \rm \\
We observe that, as each $ $ $X_{i}^{(m+1)}$ $ $ is homogeneous of degree $ $ $\beta_{i}$ $ $ (see section 3.2), we have

$$ \delta_{l}^{(m+1)}(V) = \delta_{l}^{(m+1)}(X_{m}^{(m+1)}) X_{j}^{(m+1)} + q^{-(\beta_{l}, \beta_{m})} X_{m}^{(m+1)} \delta_{l}^{(m+1)}(X_{j}^{(m+1)})$$
$$- q^{-(\beta_{m}, \beta_{j})} \delta_{l}^{(m+1)}(X_{j}^{(m+1)}) X_{m}^{(m+1)} - q^{-[ (\beta_{m}, \beta_{j}) + (\beta_{l}, \beta_{j}) ]} X_{j}^{(m+1)} \delta_{l}^{(m+1)}(X_{m}^{(m+1)}).$$
So, we can write
$$ \delta_{l}^{(m+1)}(V) = A + B - C - D$$
with
$$A = \delta_{l}^{(m+1)}(X_{m}^{(m+1)}) X_{j}^{(m+1)},$$
$$B = q^{-(\beta_{l}, \beta_{m})} X_{m}^{(m+1)} \delta_{l}^{(m+1)}(X_{j}^{(m+1)}),$$
$$C = q^{-(\beta_{m}, \beta_{j})} \delta_{l}^{(m+1)}(X_{j}^{(m+1)}) X_{m}^{(m+1)},$$
$$D = q^{-[ (\beta_{m}, \beta_{j}) + (\beta_{l}, \beta_{j}) ]} X_{j}^{(m+1)} \delta_{l}^{(m+1)}(X_{m}^{(m+1)}).$$

Now, we compute separately the different pieces of the right member: \\ $ $ \\
\begin{itemize}
\item $\delta_{l}^{(m+1)}(X_{m}^{(m+1)})$ = $X_{l}^{(m+1)}X_{m}^{(m+1)}$ $-$ $q^{-(\beta_{l}, \beta_{m} )}X_{m}^{(m+1)}X_{l}^{(m+1)}$. $ $ As we also have $ $ $\delta_{m}^{(m+1)}(X_{l}^{(m+1)})$ = $X_{m}^{(m+1)}X_{l}^{(m+1)}$ $-$ $q^{-(\beta_{l}, \beta_{m} )}X_{l}^{(m+1)}X_{m}^{(m+1)}$, $ $ we can write $ $ $X_{l}^{(m+1)}X_{m}^{(m+1)}$ = $q^{(\beta_{l}, \beta_{m} )}$ [ $X_{m}^{(m+1)}X_{l}^{(m+1)}$ $-$ $\delta_{m}^{(m+1)}(X_{l}^{(m+1)})$ ]. \\ So, we get:

$$\delta_{l}^{(m+1)}(X_{m}^{(m+1)}) \mbox{ } = (q^{(\beta_{l}, \beta_{m} )} - q^{-(\beta_{l}, \beta_{m} )}) X_{m}^{(m+1)}X_{l}^{(m+1)} -  q^{(\beta_{l}, \beta_{m} )} \delta_{m}^{(m+1)}(X_{l}^{(m+1)}).$$
\item From this, we deduce that \\ $ $ \\
$ $ $ $ $ $ $\delta_{l}^{(m+1)}(X_{m}^{(m+1)}) X_{j}^{(m+1)}$ = $(q^{(\beta_{l}, \beta_{m} )} - q^{-(\beta_{l}, \beta_{m} )}) X_{m}^{(m+1)}X_{l}^{(m+1)}X_{j}^{(m+1)} -  q^{(\beta_{l}, \beta_{m} )} \delta_{m}^{(m+1)}(X_{l}^{(m+1)})X_{j}^{(m+1)}$. \\ $ $ \\ But we also have \\ $ $ \\ $\delta_{l}^{(m+1)}(X_{j}^{(m+1)})$ = $X_{l}^{(m+1)}X_{j}^{(m+1)}$ $-$ $q^{-(\beta_{l}, \beta_{j} )}X_{j}^{(m+1)}X_{l}^{(m+1)}$, $ $ so that $ $ $ $ $ $ 
$ \delta_{l}^{(m+1)}(X_{m}^{(m+1)}) X_{j}^{(m+1)} = (q^{(\beta_{l}, \beta_{m})} - q^{-(\beta_{l}, \beta_{m})}) X_{m}^{(m+1)}[ q^{-(\beta_{l}, \beta_{j})}X_{j}^{(m+1)}X_{l}^{(m+1)} 
 + \delta_{l}^{(m+1)}(X_{j}^{(m+1)}) ] -  q^{(\beta_{l}, \beta_{m})} \delta_{m}^{(m+1)}(X_{l}^{(m+1)})X_{j}^{(m+1)}. $ 
\\ $ $ \\
So, we get: \\ $ $ \\

\begin{center}
$A$ = $(q^{(\beta_{l}, \beta_{m})} - q^{-(\beta_{l}, \beta_{m})}) q^{-(\beta_{l}, \beta_{j})} X_{m}^{(m+1)} X_{j}^{(m+1)}X_{l}^{(m+1)}$ \\ + $(q^{(\beta_{l}, \beta_{m})} - q^{-(\beta_{l}, \beta_{m})}) X_{m}^{(m+1)} \delta_{l}^{(m+1)}(X_{j}^{(m+1)})$ \\ $-$  $q^{(\beta_{l}, \beta_{m})} \delta_{m}^{(m+1)}(X_{l}^{(m+1)})X_{j}^{(m+1)}.$
\end{center} $ $
\item Now, we have $ $ $A$ + $B$ $-$ $C$ = \\
$(q^{(\beta_{l}, \beta_{m})} - q^{-(\beta_{l}, \beta_{m})}) q^{-(\beta_{l}, \beta_{j})} X_{m}^{(m+1)} X_{j}^{(m+1)}X_{l}^{(m+1)}$ \\ + $q^{(\beta_{l}, \beta_{m})} X_{m}^{(m+1)} \delta_{l}^{(m+1)}(X_{j}^{(m+1)})$ $-$
$q^{-(\beta_{m}, \beta_{j})} \delta_{l}^{(m+1)}(X_{j}^{(m+1)}) X_{m}^{(m+1)}$ \\ $-$
$q^{(\beta_{l}, \beta_{m})} \delta_{m}^{(m+1)}(X_{l}^{(m+1)})X_{j}^{(m+1)}$ = \\

$(q^{(\beta_{l}, \beta_{m})} - q^{-(\beta_{l}, \beta_{m})}) q^{-(\beta_{l}, \beta_{j})} X_{m}^{(m+1)} X_{j}^{(m+1)}X_{l}^{(m+1)}$ \\
+ $q^{(\beta_{l}, \beta_{m})} [ X_{m}^{(m+1)} \delta_{l}^{(m+1)}(X_{j}^{(m+1)})$ $-$
$q^{-(\beta_{m}, \beta_{l} + \beta_{j})} \delta_{l}^{(m+1)}(X_{j}^{(m+1)}) X_{m}^{(m+1)} ]$ \\ $-$
$q^{(\beta_{l}, \beta_{m})} \delta_{m}^{(m+1)}(X_{l}^{(m+1)})X_{j}^{(m+1)}.$ \\

Since $ $ $\delta_{l}^{(m+1)}(X_{j}^{(m+1)})$ $ $ is homogeneous of degree $ $ $\beta_{l}$ + $\beta_{j},$ $ $ we have \\

$X_{m}^{(m+1)} \delta_{l}^{(m+1)}(X_{j}^{(m+1)})$ $-$
$q^{-(\beta_{m}, \beta_{l} + \beta_{j})} \delta_{l}^{(m+1)}(X_{j}^{(m+1)}) X_{m}^{(m+1)}$ = $\delta_{m}^{(m+1)} \circ \delta_{l}^{(m+1)}(X_{j}^{(m+1)}).$ \\ So, we get: \\ $ $ \\

\begin{center}
$A$ + $B$ $-$ $C$ = \\
$(q^{(\beta_{l}, \beta_{m})} - q^{-(\beta_{l}, \beta_{m})}) q^{-(\beta_{l}, \beta_{j})} X_{m}^{(m+1)} X_{j}^{(m+1)}X_{l}^{(m+1)}$ \\
+ $q^{(\beta_{l}, \beta_{m})} \delta_{m}^{(m+1)} \circ \delta_{l}^{(m+1)}(X_{j}^{(m+1)})$ $-$
$q^{(\beta_{l}, \beta_{m})} \delta_{m}^{(m+1)}(X_{l}^{(m+1)})X_{j}^{(m+1)}.$
\end{center} $ $ 
\item   Using the computation of $ $ $\delta_{l}^{(m+1)}(X_{m}^{(m+1)})$ $ $ made in the first point over, we can write: \\

$D$ = $q^{-[ (\beta_{m}, \beta_{j}) + (\beta_{l}, \beta_{j}) ]} X_{j}^{(m+1)} \delta_{l}^{(m+1)}(X_{m}^{(m+1)})$ =
$q^{-(\beta_{j}, \beta_{l} + \beta_{m})} X_{j}^{(m+1)} [ (q^{(\beta_{l}, \beta_{m} )} - q^{-(\beta_{l}, \beta_{m} )}) X_{m}^{(m+1)} X_{l}^{(m+1)} -  q^{(\beta_{l}, \beta_{m} )} \delta_{m}^{(m+1)}(X_{l}^{(m+1)}) ].$ \\ So, we get: \\ $ $ \\

\begin{center}
$D$ = $q^{- (\beta_{j}, \beta_{l} + \beta_{m})} (q^{(\beta_{l}, \beta_{m} )} - q^{-(\beta_{l}, \beta_{m} )}) X_{j}^{(m+1)} X_{m}^{(m+1)} X_{l}^{(m+1)}$ \\ $-$ $q^{-(\beta_{j}, \beta_{l} + \beta_{m})}  q^{(\beta_{l}, \beta_{m} )} X_{j}^{(m+1)} \delta_{m}^{(m+1)}(X_{l}^{(m+1)}).$
\end{center} $ $
\item Now, we have \\
$\delta_{l}^{(m+1)}(V)$ = $A + B - C - D$ = $q^{- (\beta_{j}, \beta_{l}}) (q^{(\beta_{l}, \beta_{m})} - q^{-(\beta_{l}, \beta_{m})}) [ X_{m}^{(m+1)} X_{j}^{(m+1)} - q^{- (\beta_{j}, \beta_{m})} X_{j}^{(m+1)} X_{m}^{(m+1)} ] X_{l}^{(m+1)}$ \\ + $q^{(\beta_{l}, \beta_{m} )} [ q^{-(\beta_{j}, \beta_{l} + \beta_{m})}  X_{j}^{(m+1)} \delta_{m}^{(m+1)}(X_{l}^{(m+1)}) - \delta_{m}^{(m+1)}(X_{l}^{(m+1)})X_{j}^{(m+1)} ]$ \\ + $q^{(\beta_{l}, \beta_{m})} \delta_{m}^{(m+1)} \circ \delta_{l}^{(m+1)}(X_{j}^{(m+1)}).$ \\ So, we get: \\ $ $ \\

$$\delta_{l}^{(m+1)}(V) = q^{- (\beta_{j}, \beta_{l}}) (q^{(\beta_{l}, \beta_{m})} - q^{-(\beta_{l}, \beta_{m})}) V X_{l}^{(m+1)}$$
$$+ \mbox{ } q^{(\beta_{l}, \beta_{m} )} [ q^{-(\beta_{j}, \beta_{l} + \beta_{m})}  X_{j}^{(m+1)} \delta_{m}^{(m+1)}(X_{l}^{(m+1)}) - \delta_{m}^{(m+1)}(X_{l}^{(m+1)})X_{j}^{(m+1)} ]$$
$$+  \mbox{ } q^{(\beta_{l}, \beta_{m})} \delta_{m}^{(m+1)} \circ \delta_{l}^{(m+1)}(X_{j}^{(m+1)})$$

As $ $ $V$ $ $ is homogeneous of degree $ $ $\beta_{j} + \beta_{m}$, $ $ we have also

$$ \delta_{l}^{(m+1)}(V) = X_{l}^{(m+1)} V - q^{-(\beta_{l}, \beta_{j} + \beta_{m})} V X_{l}^{(m+1)},$$ which gives the required formula.
\end{itemize}  \bx

Now, let us observe that, if $ $ $j$ $<$ $l$ $<$ $m$, we have the following results: \\ $ $ \\
\begin{itemize}
\item $\delta_{m}^{(m+1)}(X_{l}^{(m+1)})$ $\in$ $k< X_{l+1}^{(m+1)}, \mbox{ } ... \mbox{ } , \mbox{ } X_{m-1}^{(m+1)}>$ $ $ (or is zero if $ $ $l$ = $m-1$) $ $ (see section 3.2).
\item Similarly, $ $ $V$ = $\delta_{m}^{(m+1)}(X_{j}^{(m+1)})$ $\in$ $k< X_{j+1}^{(m+1)}, \mbox{ } ... \mbox{ } , \mbox{ } X_{m-1}^{(m+1)}>$ $ $ (observe that $ $ $j < m-1$ $ $ by lemma 4.4. 1).
\item $\delta_{l}^{(m+1)}(X_{j}^{(m+1)})$ $\in$ $k< X_{j+1}^{(m+1)}, \mbox{ } ... \mbox{ } , \mbox{ } X_{l-1}^{(m+1)}>$ $ $ (or is zero if $ $ $l$ = $j+1$), $ $ which implies that $ $ $\delta_{m}^{(m+1)} \circ \delta_{l}^{(m+1)}(X_{j}^{(m+1)})$ $\in$ $k< X_{j+1}^{(m+1)}, \mbox{ } ... \mbox{ } , \mbox{ } X_{m-1}^{(m+1)}>$.
\end{itemize}

So, $ $ $V$, $ $ $\delta_{m}^{(m+1)}(X_{l}^{(m+1)})$, $ $ $\delta_{l}^{(m+1)} \circ \delta_{l}^{(m+1)}(X_{j}^{(m+1)})$, $ $ $X_{l}^{(m+1)}$, $ $ $X_{j}^{(m+1)}$ $ $ are all in $ $ $k< X_{1}^{(m+1)}, \mbox{ } ... \mbox{ } , \mbox{ } X_{m-1}^{(m+1)}>$ $ $ and, if we transform the equality of lemma 4.4. 2 by the algebra homomorphism $ $ $\Theta^{(m)}$, $ $ we obtain:

\begin{lem4.4.} $ $ \\
$$X_{l}^{(m)} U \mbox{ } = \mbox{ } q^{(\beta_{l}, \beta_{m}-\beta_{j})} U X_{l}^{(m)}$$
$$ + \mbox{ }  q^{(\beta_{l}, \beta_{m})} \mbox{ } [ \mbox{ } q^{-(\beta_{j}, \beta_{l}+\beta_{m})} X_{j}^{(m)} \Theta^{(m)} (\delta_{m}^{(m+1)}( X_{l}^{(m+1)}))  \mbox{ } - \mbox{ } \Theta^{(m)} (\delta_{m}^{(m+1)} ( X_{l}^{(m+1)} )) X_{j}^{(m)} \mbox{ } ]$$
$$ + \mbox{ } q^{(\beta_{l}, \beta_{m})} \Theta^{(m)} (\delta_{m}^{(m+1)} \circ \delta_{l}^{(m+1)}  ( X_{j}^{(m+1)}\rm )).$$
\end{lem4.4.}

$ $ \\

$\bullet$ $ $ $ $ It results of the choice of $ $ $j$ $ $ that $ $ $\Theta^{(m)} (\delta_{m}^{(m+1)} ( X_{l}^{(m+1)}))$ $\in$ $\mathcal{P}\it^{(m)}$. $ $ \\ $ $ \\

Assume that $ $ $l$ $>$ $j+1$, $ $ so that (by section 3.2)
$$\delta_{l}^{(m+1)}(X_{j}^{(m+1)})  \mbox{ } = \mbox{ } \sum_{ \underline{a} \mbox{ } = \mbox{ } (a_{j+1}, ... ,a_{l-1})} c_{\underline{a}}(X_{j+1}^{(m+1)})^{a_{j+1}} ... (X_{l-1}^{(m+1)})^{a_{l-1}}$$
with each $ $ $c_{\underline{a}}$ $\in$ $k$. \\ $ $ \\
As each $ $ $X_{i}^{(m+1)}$ $ $ is a $ $ $h_{m}$ - eigenvector, each $ $ $\delta_{m}^{(m+1)} (X_{j+1}^{(m+1)})^{a_{j+1}} ... (X_{l-1}^{(m+1)})^{a_{l-1}})$ $ $ is a linear combination of products $ $ $M_{1} \delta_{m}^{(m+1)}(X_{i}^{(m+1)}) M_{2}$ $ $ with $ $ $j$ $<$ $i$ $<$ $l$ $ $ and $ $ $M_{1}$ $ $ (resp. $ $ $M_{2}$) $ $ an ordered monomial in $ $ $X_{j+1}^{(m+1)}$, $...$ , $X_{i}^{(m+1)}$ $ $ (resp. in $ $ $X_{i}^{(m+1)}$, $...$ , $X_{l-1}^{(m+1)}$). $ $ This implies that $ $ $\Theta^{(m)} (\delta_{m}^{(m+1)} \circ \delta_{l}^{(m+1)} \rm (\it X_{j}^{(m+1)}\rm ))$ $ $ is a linear combination of products $ $ $N_{1} \Theta^{(m)} (\delta_{m}^{(m+1)}(X_{i}^{(m+1)})) N_{2}$ $ $ with $ $ $j$ $<$ $i$ $<$ $l$, $ $ $N_{1}$ $ $ (resp. $ $ $N_{2}$) $ $ an ordered monomial in $ $ $X_{j+1}^{(m)}$, $...$ , $X_{i}^{(m)}$ $ $ (resp. in $ $ $X_{i}^{(m)}$, $...$ , $X_{l-1}^{(m)}$). $ $ \\
As $ $ $\Theta^{(m)} (\delta_{m}^{(m+1)}(X_{i}^{(m+1)}))$ $\in$ $\mathcal{P}\it ^{(m)}$ $ $ for $ $ $j$ $<$ $i$ $<$ $m$, $ $ we conclude that \\ $ $ \\

$\bullet$ $ $ $ $ $\Theta^{(m)} (\delta_{m}^{(m+1)} \circ \delta_{l}^{(m+1)} \rm (\it X_{j}^{(m+1)}\rm ))$ $\in$ $\mathcal{P}\it^{(m)}$. \\ $ $ \\

We observe that this result is also true when $ $ $l$ = $j+1$ $ $ since, in this case, $ $ $\delta_{l}^{(m+1)}(X_{j}^{(m+1)})$ $ $ is zero (see section 3.2). So, if we transform the equality of lemma 4.4. 3  by the algebra homomorphism $ $ $f_{m}$, $ $ we obtain:

\begin{lem4.4.} $ $ \\
\begin{equation}\label{xlu}
x_{l}^{(m)} u \mbox{ } = \mbox{ } q^{(\beta_{l}, \beta_{m}-\beta_{j})} u x_{l}^{(m)}.
\end{equation}
\end{lem4.4.}

$ $ \\ Recall (see section 3.2) that $ $ $V$ = $\delta_{m}^{(m+1)}(X_{j}^{(m+1)}))$ $ $ is an homogeneous polynomial in $ $ $X_{j+1}^{(m+1)}$, $...$ , $X_{m-1}^{(m+1)}$ $ $ of degree $ $ $\beta_{m}$ + $\beta_{j}$. $ $ So, since $ $ $\Theta^{(m)}$ $ $ transforms each $ $ $X_{i}^{(m+1)}$ $ $ ($1$ $\leq$ $i$ $<$ $m$) $ $ in $ $ $X_{i}^{(m)}$, $ $ $U$ = $\Theta^{(m)} (V)$ $ $ is an homogeneous polynomial in $ $ $X_{j+1}^{(m)}$, $...$ , $X_{m-1}^{(m)}$ $ $ of degree $ $ $\beta_{m}$ + $\beta_{j}$. $ $ This implies:

\begin{lem4.4.} $ $ \\
$u$ $ $ is a nonzero polynomial in $ $ $x_{j+1}^{(m)}$, $...$ , $x_{m-1}^{(m)}$. $ $ Moreover, $ $ $u$ $ $ is homogeneous of degree $ $ $\beta_{m}$ + $\beta_{j}$.
\end{lem4.4.}

\bf Proof \rm \\
$u$ = $f_{m}$($U$) $ $ is nonzero by construction. Since $ $ $f_{m}$ $ $ transforms each $ $ $X_{i}^{(m)}$ $ $ in $ $ $x_{i}^{(m)}$, $ $ it turns out that $u$ $ $ is a polynomial in $ $ $x_{j+1}^{(m)}$, $...$ , $x_{m-1}^{(m)}$. Moreover,  $u = f_{m}$($U$)$ $  is homogeneous of degree $ $ $\beta_{m}$ + $\beta_{j}$ by lemma 4.3. 1 \bx

Recall (lemma 4.4. 1) that $ $ $j$ $<$ $l_{p}$ $ $ and that $ $ $p$ $ $ is the greatest integer such that $ $ $l_{p}$ $<$ $m$. $ $ Denote by $ $ $c$ $ $ the smallest integer such that $ $ $j$ $<$ $l_{c}$ $ $ ($1$ $\leq$ $c$ $\leq$ $p$). $ $ As $ $ $u$ $\in$ $ $ $k<x_{j+1}^{(m)}, \mbox{ } ... \mbox{ } , \mbox{ } x_{m-1}^{(m)}>$, $ $ it results from (proposition 4.3. 2, 2.b.) that

$$u \mbox{ } = \mbox{ } u_{1}z_{l_{p}}^{a_{1}} + \mbox{ } ... \mbox{ }+ u_{M}z_{l_{p}}^{a_{M}}$$
with $ $ $M$ $\geq$ $1$, $ $ ($a_{1}, \mbox{ } ... \mbox{ }, \mbox{ } a_{M}$) $\in$ \bbZ$^{M}$ $ $ and: \\
$\bullet$ $ $ If $ $ $p$ = $1$, $ $  then each $ $ $u_{i}$ $ $ $\in$ $ $ $k$. \\
$\bullet$ $ $ If $ $ $p$ $>$ $1$ and $ $ $j$ $\geq$ $l_{p-1}$, $ $ then each $ $ $u_{i}$ $ $ $\in$ $ $ $k$. \\
$\bullet$ $ $ If $ $ $p$ $>$ $1$ and $ $ $j$ $<$ $l_{p-1}$, $ $ then each $ $ $u_{i}$ $ $ $\in$ $ $ $k<x_{j+1}^{(l_{p})}, \mbox{ } ... \mbox{ } , \mbox{ } x_{l_{p-1}}^{(l_{p})}>$, \\ $ $ \\

so that we can write \\ $ $ \\

$\bullet$ $ $ If $ $ $c$ = $p$, $ $  then each $ $ $u_{i}$ $ $ $\in$ $ $ $k$. \\
$\bullet$ $ $ If $ $ $c$ $<$ $p$, $ $ then each $ $ $u_{i}$ $ $ $\in$ $ $ $k<x_{j+1}^{(l_{p})}, \mbox{ } ... \mbox{ } , \mbox{ } x_{l_{p-1}}^{(l_{p})}>$. \\ $ $ \\

In both cases, we may assume that $ $ $u_{1}$, $...$ , $u_{M}$ $ $ are all nonzero, that $ $ $a_{1}$ $<$ $...$ $<$ $a_{M}$, $ $ and we observe that \\ $ $ \\

$\bullet$ $ $ Each $ $ $u_{i}$ $\in$ $k<x_{1}^{(l_{p})}, \mbox{ } ... \mbox{ } , \mbox{ } x_{l_{p-1}}^{(l_{p})}>$ $ $ and is homogeneous of degree $ $ $\beta_{m} + \beta_{j} - a_{i} \beta_{l_{p}}$. \\
In fact, since $ $ $u$ $ $ is homogeneous of degree $ $ $\beta_{m} + \beta_{j}$ $ $ and $ $ $z_{l_{p}}$ = $x_{l_{p}}^{(l_{p})}$ $ $ is homogeneous of degree $ $ $\beta_{l_{p}}$, $ $ we have, for any $ $ $\rho$ $\in$ $ $ \bbZ $\Pi$,

$$\overline{h_{\rho}}(u) \mbox{ } = \mbox{ } q^{-(\rho , \mbox{ } a_{1} \beta_{l_{p}})} \overline{h_{\rho}}(u_{1})z_{l_{p}}^{a_{1}} + \mbox{ } ... \mbox{ }+ q^{-(\rho , \mbox{ } a_{M} \beta_{l_{p}})} \overline{h_{\rho}}(u_{M})z_{l_{p}}^{a_{M}}$$

$$= \mbox{ } q^{-(\rho , \mbox{ } \beta_{m} + \beta_{j})} \mbox{ } [ u_{1}z_{l_{p}}^{a_{1}} + \mbox{ } ... \mbox{ }+ u_{M}z_{l_{p}}^{a_{M}} ]$$

and, by corollary 4.3. 2, we can identify the coefficients of $ $ $z_{l_{p}}^{a_{i}}$, $ $ so that $ $ $q^{-(\rho , \mbox{ } a_{i} \beta_{l_{p}})} \overline{h_{\rho}}(u_{i})$ = $q^{-(\rho , \mbox{ } \beta_{m} + \beta_{j})}u_{i}$ $ $ and $ $ $\overline{h_{\rho}}(u_{i})$ = $q^{-(\rho , \mbox{ } \beta_{m} + \beta_{j} - a_{i} \beta_{l_{p}})}u_{i}.$ \\ $ $ \\ $ $ \\

By (proposition 4.3. 3, 1.) we have $ $ $x_{l_{p}}^{(m)}$ = $z_{l_{p}}$ = $x_{l_{p}}^{(l_{p})}$. $ $ So, using lemma 4.4. 4 with $ $ $l$ = $l_{p}$, $ $ we obtain:

$$z_{l_{p}} u \mbox{ } = \mbox{ } q^{(\beta_{l_{p}}, \mbox{ } \beta_{m} - \beta_{j})}) u z_{l_{p}}$$

and, by (lemma 4.3. 1, 4.),

$$z_{l_{p}} u \mbox{ } = \mbox{ } \overline{h_{l_{p}}}(u_{1})z_{l_{p}}^{a_{1}+1} + \mbox{ } ... \mbox{ }+ \overline{h_{l_{p}}}(u_{M})z_{l_{p}}^{a_{M}+1}.$$

This implies that

$$q^{(\beta_{l_{p}}, \mbox{ } \beta_{m} - \beta_{j})} \mbox{ } [  u_{1}z_{l_{p}}^{a_{1}} + \mbox{ } ... \mbox{ }+ u_{M}z_{l_{p}}^{a_{M}} ] z_{l_{p}}$$

$$= \mbox{ } [ q^{-(\beta_{l_{p}} , \mbox{ } \beta_{m} + \beta_{j} - a_{1} \beta_{l_{p}})} u_{1}z_{l_{p}}^{a_{1}} + \mbox{ } ... \mbox{ }+ q^{-(\beta_{l_{p}} , \mbox{ } \beta_{m} + \beta_{j} - a_{N} \beta_{l_{p}})} u_{M}z_{l_{p}}^{a_{M}} ] z_{l_{p}}.$$

As over, if $ $ $1$ $\leq$ $i$ $\leq$ $M$, $ $ we can identify the coefficients of $ $ $z_{l_{p}}^{a_{i}+1}$, $ $ so that, since $ $ $u_{i}$ $\neq$ $0$ and q is not a root of unity, we have $- (\beta_{l_{p}} , \beta_{m} + \beta_{j} - a_{i} \beta_{l_{p}}) \mbox{ } = \mbox{ } (\beta_{l_{p}}, \beta_{m} - \beta_{j}).$ \\ $ $ \\

So, $ $ $a_{i} \| \beta_{l_{p}}  \|^{2}$ = $2(\beta_{l_{p}}, \beta_{m})$ $ $ and then $ $
$a_{i} \mbox{ } = \mbox{ } (\beta_{l_{p}}^{\vee}, \beta_{m}).$ \\ $ $ \\

This implies that $ $ $M$ = $1$ $ $ and we conclude:

\begin{lem4.4.} $ $ \\
$$u \mbox{ } = \mbox{ } u_{p}z_{l_{p}}^{a_{p}}$$
with \\
\begin{itemize}
\item If $ $ $c$ = $p$, $ $  then $ $ $u_{p}$ $ $ $\in$ $ $ $k^{*}$.
\item If $ $ $c$ $<$ $p$, $ $ then $ $ $u_{p}$ $ $ $\in$ $ $ $k<x_{j+1}^{(l_{p})}, \mbox{ } ... \mbox{ } , \mbox{ } x_{l_{p-1}}^{(l_{p})}>$ $\setminus$ $\{ 0 \}$.
\item $a_{p} \mbox{ } = \mbox{ } (\beta_{l_{p}}^{\vee}, \beta_{m}).$
\item $u_{p}$ $ $ is homogeneous of degree $ $ $\beta_{m} + \beta_{j} - a_{p} \beta_{l_{p}}$ = $\gamma_{p} + \beta_{j}$ $ $ with $ $ $\gamma_{p}$ = $s_{\beta_{l_{p}}} (\beta_{m}).$
\end{itemize}
\end{lem4.4.}

Now, we extend this result as follows:

\begin{lem4.4.} $ $ \\
Consider some integer $ $ $d$ $ $ with $ $ $c$ $\leq$ $d$ $\leq$ $p.$ $ $ Then
$$u \mbox{ } = \mbox{ } u_{d} z_{l_{d}}^{a_{d}} \mbox{ } ... \mbox{ } z_{l_{p}}^{a_{p}}$$
with \\
\begin{itemize}
\item If $ $ $c$ = $d$, $ $  then $ $ $u_{d}$ $ $ $\in$ $ $ $k^{*}$.
\item If $ $ $c$ $<$ $d$, $ $ then $ $ $u_{d}$ $ $ $\in$ $ $ $k<x_{j+1}^{(l_{d})}, \mbox{ } ... \mbox{ } , \mbox{ } x_{l_{d-1}}^{(l_{d})}>$ $\setminus$ $\{ 0 \}$.
\item Denote by $(\gamma_{d}$, $...$ , $\gamma_{p}$, $\gamma_{p+1})$ the sequence of $($non necessarily positive$)$ roots recursively defined by $ $ $\gamma_{p+1}$ = $\beta_{m}$ $ $ and, for $ $ $d$ $\leq$ $i$ $\leq$ $p$, $ $  $\gamma_{i} = s_{\beta_{l_{i}}}(\gamma_{i+1})$. Then, each $ $ $a_{i} = (\beta_{l_{i}}^{\vee} , \gamma_{i+1})$.
\item $u_{d}$ $ $ is homogeneous of degree $ $ $\gamma_{d} + \beta_{j}.$
\end{itemize}
\end{lem4.4.}

\bf Proof \rm \\

We proceed by decreasing induction on $ $ $d$. \\
If $ $ $d$ = $p$, $ $ everything has been done in lemma 4.4. 6. \\
Assume that $ $ $c$ $\leq$ $d$ $<$ $p$ $ $ and that the lemma is true at the rank $ $ $d + 1$, $ $ namely: \\ $ $ \\

$$u \mbox{ } = \mbox{ } u_{d+1} z_{l_{d+1}}^{a_{d+1}} \mbox{ } ... \mbox{ } z_{l_{p}}^{a_{p}}$$
with \\
\begin{itemize}
\item $u_{d+1}$ $ $ $\in$ $ $ $k<x_{j+1}^{(l_{d+1})}, \mbox{ } ... \mbox{ } , \mbox{ } x_{l_{d}}^{(l_{d+1})}>$ $\setminus$ $\{ 0 \}$.
\item If $ $ $(\gamma_{d+1}$, $...$ , $\gamma_{p}$, $\gamma_{p+1})$ is the sequence of $($non necessarily positive$)$ roots recursively defined by $ $ $\gamma_{p+1}$ = $\beta_{m}$ $ $ and, for $ $ $d+1$ $\leq$ $i$ $\leq$ $p$, $ $  $\gamma_{i} = s_{\beta_{l_{i}}}(\gamma_{i+1})$, then, each $ $ $a_{i} = (\beta_{l_{i}}^{\vee} , \gamma_{i+1})$.
\item $u_{d+1}$ $ $ is homogeneous of degree $ $ $\gamma_{d+1} + \beta_{j}.$
\end{itemize}

By (proposition 4.3. 2, 2.b.), we have:

$$u_{d+1} \mbox{ } = \mbox{ } v_{1} z_{l_{d}}^{b_{1}} \mbox{ } +  \mbox{ } ... \mbox{ } + \mbox{ }  v_{M} z_{l_{d}}^{b_{M}}$$
with $ $ $M$ $\geq$ $1$, $ $ ($b_{1}, \mbox{ } ... \mbox{ }, \mbox{ } b_{M}$) $\in$ \bbZ$^{M}$ $ $ and: \\
\begin{itemize}
\item If $ $ $d$ = $c$, $ $  then each $ $ $v_{i}$ $ $ $\in$ $ $ $k$.
\item If $ $ $d$ $>$ $c$, $ $ then each $ $ $v_{i}$ $ $ $\in$ $ $ $k<x_{j+1}^{(l_{d})}, \mbox{ } ... \mbox{ } , \mbox{ } x_{l_{d-1}}^{(l_{d})}>$.
\end{itemize}
 $ $ \\

As over, we may assume that $ $ $v_{1}$, $...$ , $v_{M}$ $ $ are all nonzero, that $ $ $b_{1}$ $<$ $...$ $<$ $b_{M}$, $ $ and observe that \\ $ $ \\

$\bullet$ $ $ Each $ $ $v_{i}$ $ $ is homogeneous of degree $ $ $\gamma_{d+1} + \beta_{j} - b_{i} \beta_{l_{d}}$. \\ $ $ \\

By (proposition 4.3. 3, 2.), we have

$$x_{l_{d}}^{(m)} \mbox{ } = \mbox{ } z_{l_{d}} \mbox{ } + \mbox{ } \displaystyle \sum_{\underline{s} \mbox{ } = \mbox{ } (s_{d+1}, \mbox{ }  ... \mbox{ } , \mbox{ } s_{p}) \mbox{ } \in \mbox{ } F}{} \mbox{ } \eta(\underline{s}) z_{l_{d+1}}^{s_{d+1}} \mbox{ }  ... \mbox{ } z_{l_{p}}^{s_{p}},$$
where \\
\begin{itemize}
\item $F$ $ $ is a finite (possibly empty) subset of $ $ \bbZ $^{p-d}$.
\item If $ $ $\preceq$ $ $ denotes the inverse lexicographic order in $ $ \bbZ$^{p-d}$, $ $ then, for any $ $ $\underline{s}$ $ $ in $ $ $F$, $ $ we have $ $ $\underline{s}$ $ $ $\prec$ $ $ $0$.
\item For each $ $ $\underline{s}$ $ $ in $ $ $F$, $ $ we have $ $ $\eta(\underline{s}) \mbox{ } \in \mbox{ } k^{*}$.
\end{itemize}
$ $ \\

By (lemma 4.3. 1, 4.), we have

$$z_{l_{d}} u_{d+1} \mbox{ } = \mbox{ } \overline{h_{l_{d}}}(v_{1})z_{l_{d}}^{b_{1}+1} + \mbox{ } ... \mbox{ }+ \overline{h_{l_{d}}}(v_{M})z_{l_{d}}^{b_{M}+1}.$$

For any $ $ $\underline{e} \mbox{ } = \mbox{ } (e_{d+1}, \mbox{ }  ... \mbox{ } , \mbox{ } e_{p})$ $\in$ $ $ \bbZ $ $ $^{p-d}$, $ $ let us denote $ $ $z^{\underline{e}}$ := $z_{l_{d+1}}^{e_{d+1}} \mbox{ }  ... \mbox{ } z_{l_{p}}^{e_{p}}$. \\ Again by (lemma 4.3. 1, 4.), there exists $ $ ($\lambda_{\underline{e},1}, \mbox{ }  ... \mbox{ } , \mbox{ } \lambda_{\underline{e},M})$ $\in$ $(k^{*})^{M}$ such that

$$z^{\underline{e}} u_{d+1} \mbox{ } = \mbox{ } \lambda_{\underline{e},1}v_{1}z_{l_{d}}^{b_{1}}z^{\underline{e}} + \mbox{ } ... \mbox{ }+ \lambda_{\underline{e},M}v_{M}z_{l_{d}}^{b_{M}}z^{\underline{e}}.$$

Now, if we set $ $ $\underline{a} \mbox{ } = \mbox{ } (a_{d+1}, \mbox{ }  ... \mbox{ } , \mbox{ } a_{p})$, $ $ we have

$$x_{l_{d}}^{(m)} u \mbox{ } = \mbox{ } (1) \mbox{ } + \mbox{ } (2)$$

with

$$(1) \mbox{ } = \mbox{ } z_{l_{d}} u_{d+1} z^{\underline{a}} \mbox{ } = \mbox{ } \overline{h_{l_{d}}}(v_{1})z_{l_{d}}^{b_{1}+1}z^{\underline{a}} + \mbox{ } ... \mbox{ }+ \overline{h_{l_{d}}}(v_{M})z_{l_{d}}^{b_{M}+1}z^{\underline{a}}$$

$$= \mbox{ } q^{-(\beta_{l_{d}},\mbox{ } \gamma_{d+1} + \beta_{j} - b_{1} \beta_{l_{d}})}v_{1}z_{l_{d}}^{b_{1}+1}z^{\underline{a}} + \mbox{ } ... \mbox{ }+ q^{-(\beta_{l_{d}},\mbox{ } \gamma_{d+1} + \beta_{j} - b_{M} \beta_{l_{d}})}v_{M}z_{l_{d}}^{b_{M}+1}z^{\underline{a}}$$

and

$$(2) \mbox{ } = \mbox{ } ( \displaystyle \sum_{\underline{s} \mbox{ } \in \mbox{ } F}{} \mbox{ } \eta(\underline{s}) z^{\underline{s}} \mbox{ } ) \mbox{ } u_{d+1} z^{\underline{a}} \mbox{ } = \mbox{ } \displaystyle \sum_{\underline{s} \mbox{ } \in \mbox{ } F}{} \mbox{ } \eta(\underline{s}) \mbox{ }  (\lambda_{\underline{s},1}v_{1}z_{l_{d}}^{b_{1}}z^{\underline{s}} + \mbox{ } ... \mbox{ }+ \lambda_{\underline{s},M}v_{M}z_{l_{d}}^{b_{M}}z^{\underline{s}}\mbox{ } )\mbox{ } z^{\underline{a}}.$$

Besides, we also have

$$u x_{l_{d}}^{(m)}  \mbox{ } = \mbox{ } (1^{\prime}) \mbox{ } + \mbox{ } (2^{\prime})$$

with

$$(1^{\prime}) \mbox{ } = \mbox{ }  u_{d+1} z^{\underline{a}}z_{l_{d}} \mbox{ } = \mbox{ } v_{1}z_{l_{d}}^{b_{1}}z^{\underline{a}}z_{l_{d}} + \mbox{ } ... \mbox{ }+ v_{M}z_{l_{d}}^{b_{M}}z^{\underline{a}}z_{l_{d}}$$

$$= \mbox{ } q^{-(a_{d+1}\beta_{l_{d+1}} \mbox{ } +  \mbox{ } ... \mbox{ } + \mbox{ } a_{p}\beta_{l_{p}},\mbox{ } \beta_{l_{d}})} ( v_1z_{l_{d}}^{b_{1}+1} + \mbox{ } ... \mbox{ } + v_{M}z_{l_{d}}^{b_{M}+1} ) z^{\underline{a}}$$

and

$$(2^{\prime}) \mbox{ } = \mbox{ } u_{d+1} z^{\underline{a}} \mbox{ } [ \mbox{ } \displaystyle \sum_{\underline{s} \mbox{ } \in \mbox{ } F}{} \mbox{ } \eta(\underline{s}) z^{\underline{s}} \mbox{ } ]  \mbox{ } = \mbox{ } \displaystyle \sum_{\underline{s} \mbox{ } \in \mbox{ } F}{} \mbox{ } \eta(\underline{s}) (v_{1}z_{l_{d}}^{b_{1}}z^{\underline{a}}z^{\underline{s}} + \mbox{ } ... \mbox{ }+ v_{M}z_{l_{d}}^{b_{M}}z^{\underline{a}}z^{\underline{s}}).$$

Using lemma 4.4. 4 and  with $ $ $l$ = $l_{d}$, $ $ we obtain:

$$(1) \mbox{ } + \mbox{ } (2)  \mbox{ } = \mbox{ } q^{(\beta_{l_{d}}, \mbox{ } \beta_{m} - \beta_{j})} \mbox{ } [ \mbox{ } (1^{\prime}) \mbox{ } + \mbox{ } (2^{\prime}) \mbox{ } ].$$

We observe that $ $ $(1)$ $ $ and $ $ $(1^{\prime})$ $ $ are left linear combinations of monomials of type $ $ $z_{l_{d}}^{h}z^{\underline{a}}$ $ $ ($h$ $\in$ $ $ \bbZ) $ $ with coefficients in $ $ $A$ = $ $ $k<x_{j+1}^{(l_{d})}, \mbox{ } ... \mbox{ } , \mbox{ } x_{l_{d-1}}^{(l_{d})}>$, while $ $ $(2)$ $ $ and $ $ $(2^{\prime})$ $ $ are left linear combinations of monomials of type $ $ $z_{l_{d}}^{h}z^{\underline{e}}$ $ $ ($h$ $\in$ $ $ \bbZ , $ $ $\underline{e}$ = $\underline{a}$ + $\underline{s}$, $ $ $s$ $\in$ $F$) $ $ with coefficients in $ $ $A$. $ $ Recall that, for any $ $ $\underline{s}$ $\in$ $F$, $ $ we have $ $ $\underline{s}$ $\prec$ $ $ $0$ $ $ $\Rightarrow$ $ $  $\underline{e}$ = $\underline{a}$ + $\underline{s}$ $\prec$ $ $ $\underline{a}$. $ $ So, by corollary 4.3. 2, we have
$ $ $(1)$ = $q^{(\beta_{l_{d}}, \mbox{ } \beta_{m} - \beta_{j})} (1^{\prime})$, $ $ which implies that, if $ $ $1$ $\leq$ $i$ $\leq$ $M$,

$$q^{-(\beta_{l_{d}},\mbox{ } \gamma_{d+1} + \beta_{j} - b_{i} \beta_{l_{d}})} \mbox{ } = \mbox{ } q^{(\beta_{l_{d}}, \mbox{ } \beta_{m} - \beta_{j})} q^{-(a_{d+1}\beta_{l_{d+1}} \mbox{ } +  \mbox{ } ... \mbox{ } + \mbox{ } a_{p}\beta_{l_{p}},\mbox{ } \beta_{l_{d}})}.$$

Recall that, for $ $ $d+1$ $\leq$ $h$ $\leq$ $p$, $ $  $\gamma_{h} = s_{\beta_{l_{h}}}(\gamma_{h+1})$ = $\gamma_{h+1}$ $-$ $(\beta_{l_{h}}^{\vee} , \gamma_{h+1}) \beta_{l_{h}}$ =  $\gamma_{h+1}$ $-$ $a_{h} \beta_{l_{h}}$. $ $ If we add up all these equalities, we obtain $ $ $\gamma_{d+1}$ = $\gamma_{p+1}$ $-$ $(a_{d+1}\beta_{l_{d+1}} \mbox{ } +  \mbox{ } ... \mbox{ } + \mbox{ } a_{p}\beta_{l_{p}})$. $ $ As $ $ $\gamma_{p+1}$ = $\beta_{m}$, $ $ we get the formula

$$-(\beta_{l_{d}},\mbox{ } \gamma_{d+1} + \beta_{j} - b_{i} \beta_{l_{d}}) \mbox{ } = \mbox{ } (\beta_{l_{d}}, \mbox{ } \beta_{m} - \beta_{j}) \mbox{ } - \mbox{ } (\beta_{m}  \mbox{ } - \mbox{ } \gamma_{d+1}, \mbox{ } \beta_{l_{d}})$$

for any $ $ $i$ $\in$ $\llbracket  1$, ... , $M \rrbracket $. \\ $ $ \\

This implies that $ $ $b_{i} \| \beta_{l_{d}} \|^{2}$ $-$ $(\beta_{l_{d}},\mbox{ } \gamma_{d+1})$ = $(\beta_{l_{d}},\mbox{ } \gamma_{d+1})$. $ $ So, as the integers $ $ $b_{i}$ $ $ are distincts, we necessarily have $ $ $M$ = $1$, $ $ $b_{1}$ = $(\beta_{l_{d}}^{\vee},\mbox{ } \gamma_{d+1})$ $ $ and  $ $ $v_{1}$ $ $ is homogeneous of degree $ $ $\gamma_{d+1}$ $-$ $b_{1}\beta_{l_{d}}$ + $\beta_{j}$ = $ $ $\gamma_{d}$ + $\beta_{j}$ $ $ if we set $ $ $\gamma_{d}$ = $\gamma_{d+1}$ $-$ $b_{1}\beta_{l_{d}}$ = $s_{\beta_{l_{d}}}$($\gamma_{d+1}$). \bx

Now, in the case $ $ $d$ = $c$, $ $ lemma 4.4. 7 becomes: \\ $ $ \\

$$u \mbox{ } = \mbox{ } u_{c} z_{l_{c}}^{a_{c}} \mbox{ } ... \mbox{ } z_{l_{p}}^{a_{p}}$$
with \\
\begin{itemize}
\item $u_{c}$ $ $ $\in$ $ $ $k^{*}$.
\item If $ $ $(\gamma_{c}$, $...$ , $\gamma_{p}$, $\gamma_{p+1})$ $ $ is the sequence of $($non necessarily positive$)$ roots recursively defined by $ $ $\gamma_{p+1}$ = $\beta_{m}$ $ $ and, for $ $ $c$ $\leq$ $i$ $\leq$ $p$, $ $  $\gamma_{i} = s_{\beta_{l_{i}}}(\gamma_{i+1})$, $ $ then each $ $ $a_{i} = (\beta_{l_{i}}^{\vee} , \gamma_{i+1})$.
\end{itemize}

In this case, $ $ $u_{c}$ $ $ is homogeneous of degree $ $ $0$ $ $ and, since each $ $ $z_{l_{d}}$ = $x_{l_{d}}^{(l_{d})}$ $ $ is homogeneous of degree $ $ $\beta_{l_{d}}$ $ $ (see lemma 4.3. 1), we obtain that $ $ $u$ $ $ is also homogeneous of degree $ $ $a_{c}\beta_{l_{c}} \mbox{ } +  \mbox{ } ... \mbox{ } + \mbox{ } a_{p}\beta_{l_{p}}.$ $ $ So, as the degree of $ $ $u$ $ $ is uniquely defined, we have

$$\beta_{m} \mbox{ } + \mbox{ } \beta_{j} \mbox{ } = \mbox{ } a_{c}\beta_{l_{c}} \mbox{ } +  \mbox{ } ... \mbox{ } + \mbox{ } a_{p}\beta_{l_{p}}$$

and, by proposition 2.3. 3, we conclude:

\newtheorem{prop4.4.}{Proposition 4.4.}

\begin{prop4.4.} $ $ \\
If $ $ $\mathcal{P}\it^{(m)}$ $ $ is not in $ $ $Im(\phi_{m})$, $ $ then
$$\Delta  \mbox{ } = \mbox{ } \{i \in \mbox{ } \rm [\it\!|1,  \mbox{ } ...  \mbox{ } ,  \mbox{ } t|\!\rm] \it \mbox{ } | \mbox{ } Z_{i} \in \mbox{ } \mathcal{P}^{(2)} \}$$
is not a positive diagram with respect to $ $ $(\ref{eqn:expw})$ $ $ in the sense of definition $2.2.$ $2.$
\end{prop4.4.}

\section{Connections between admissible and positive diagrams.}

$ $ $w$ $ $ and the algebra $ $ $R$ = $U^{+} [ w ]$ $ $ are defined as in section 2. The admissible diagrams (with respect to the reduced expression (\ref{eqn:expw})) are defined as in section 3.3. \\ $ $ \\

We know that the longest element $ $ $w_{0}$ $ $ has a reduced expression of the form

\begin{equation}\label{eqn:expw0}
w_{0} = s_{\alpha_{1}} \mbox{ } ... \mbox{ } s_{\alpha_{N}} \mbox{ } \mbox{ } \mbox{ } \mbox{ } \mbox{ } (\alpha_{i} \in \Pi \mbox{ }  \mbox{ }  for \mbox{ }  \mbox{ }  1 \leq i \leq N = |\Phi^{+}|)
\end{equation}

such that $ $ $s_{\alpha_{1}} \mbox{ } ... \mbox{ } s_{\alpha_{t}}$ $ $ is the reduced expression (\ref{eqn:expw}) of $ $ $w$. \\ $ $ \\

In this section, we shall use the following conventions: \\ $ $ \\
\begin{itemize}
\item Denote by  $ $ $R_{0}$ = $U^{+} [ w_{0} ]$ $ $ the algebra associated to the reduced expression (\ref{eqn:expw0}) as explained in section 2.1. and set $ $ $F_{0}$ = $Fract(R_{0})$ $ $ it's division ring of fractions.
\item As in section 3.1, we set $ $ $X_{i}$ = $X_{\beta_{i}}$ $ $ for $ $ $1$ $\leq$ $i$ $\leq$ $N$ $ $ and we observe that $ $ $X_{1}$, $ $ ... $ $, $ $ $X_{t}$ $ $ are also the generators of $ $ $R$ = $U^{+} [ w ]$ $ $ defined in section 2.1, so that $ $ $R$ = $k<X_{1}$, $ $ ... $ $, $ $ $X_{t}>$ $ $ is a subalgebra of $ $ $R_{0}$ = $k<X_{1}$, $ $ ... $ $, $ $ $X_{N}>$.
\item For any $ $ $\rho$ $\in$ $ $ \bbZ$\Pi$ $ $ we still denote by $ $ $h_{\rho}$ $ $ the automorphism of $ $ $R_{0}$ $ $ defined as in section 3.1.
\item Each $ $ $h_{\rho}$ $ $ can uniquely be extended in an automorphism, still denoted $ $ $h_{\rho}$, $ $ of $ $ $F_{0}$.
\item We still denote by $ $ $H$ $ $ the set of all the automorphisms $ $ $h_{\rho}$ $ $ of $ $ $F_{0}$ $ $ ($\rho$ $\in$ $ $ \bbZ$\Pi$). $ $ It is a subgroup of $ $ $Aut(F_{0})$.
\item For any $ $ $m$ $ $ in $ $ $\llbracket  1,  \mbox{ } ...  \mbox{ } ,  \mbox{ } N+1 \rrbracket $, $ $ we define the algebra $ $ $R_{0}^{(m)}$ $ $  and it's canonical generators $ $ $X_{1}^{(m)},  \mbox{ } ...  \mbox{ } ,  \mbox{ } X_{N}^{(m)}$ $ $ as in section 3.2.
\end{itemize}

\subsection{A property of admissible diagrams.}

Consider $ $ $\Delta$ $ $ a diagram with respect to (\ref{eqn:expw}) (ie. a subset of   $\llbracket  1,  \mbox{ } ...  \mbox{ } ,  \mbox{ } t \rrbracket $). It is also a diagram with respect to (\ref{eqn:expw0}) but it is not quite clear whether the properties $ $ "$\Delta$ $ $ is admissible with respect to (\ref{eqn:expw})"  and $ $ "$\Delta$ $ $ is admissible with respect to (\ref{eqn:expw0})"  are equivalent. In this section, we clarify this point.

\newtheorem{lem5.1.}{Lemma 5.1.}

\begin{lem5.1.} $ $ \\
\begin{enumerate}
\item There exists a unique algebra isomorphism $g$: $ $ $R$ = $k<X_{1}$, $ $ ... $ $, $ $ $X_{t}>$ $ $ $\rightarrow$ $ $ $k<X_{1}^{(t+1)}$, $ $ ... $ $, $ $ $X_{t}^{(t+1)}>$ $ $ $\subset$ $ $ $R_{0}^{(t+1)}$ $ $ which transforms each $ $ $X_{i}$ $ $ $(1$ $\leq$ $i$ $\leq$ $t)$ $ $ in  $ $ $X_{i}^{(t+1)}$.
\item For any $ $ $\rho$ $\in$ $ $ \bbZ$\Pi$ $ $ and $B$ $\in$ $R$, $ $ we have $ $ $g$ $\circ$ $h_{\rho} (B)$ = $h_{\rho}$ $\circ$ $g (B)$.
\item The map $ $ $\mathcal{P}$ $ $ $\mapsto$ $g(\mathcal{P})$ $ $ is a bijection from the set of $ $ $H$ - invariant prime ideals of  $ $ $R$ $ $ onto the set of $ $ $H$ - invariant prime ideals of  $ $ $k<X_{1}^{(t+1)}$, $ $ ... $ $, $ $ $X_{t}^{(t+1)}>$.
\end{enumerate}

\end{lem5.1.}

\bf Proof \rm \\
\begin{enumerate}
\item We know (see section 3.1) that $ $ $R$ $ $ is the $ $ $k$ - algebra generated by the variables $ $ $X_{1}$, $ $ ... $ $, $ $ $X_{t}$ $ $ submitted to the Levendorskii - Soibelman relations (\ref{LSR}) and (\ref{LSP}). As $ $ $X_{1}^{(t+1)}$, $ $ ... $ $, $ $ $X_{t}^{(t+1)}$ $ $ satisfy the "same" relations (\ref{LSRm}) and (\ref{LSPm2}) (see section 3.2), there exists a unique homomorphism $g$: $ $ $R$ = $k<X_{1}$, $ $ ... $ $, $ $ $X_{t}>$ $ $ $\rightarrow$ $ $ $k<X_{1}^{(t+1)}$, $ $ ... $ $, $ $ $X_{t}^{(t+1)}>$ $ $ $\subset$ $ $ $R_{0}^{(t+1)}$ $ $ which transforms each $ $ $X_{i}$ $ $ ($1$ $\leq$ $i$ $\leq$ $t$) $ $ in  $ $ $X_{i}^{(t+1)}$. $ $ As $ $ $g$ transforms the base of the ordered monomials in $ $ $X_{1}$, $ $ ... $ $, $ $ $X_{t}$ $ $ in the base of the ordered monomials in $ $ $X_{1}^{(t+1)}$, $ $ ... $ $, $ $ $X_{t}^{(t+1)}$, $ $ it is an isomorphism.
\item For $ $ $1$ $\leq$ $i$ $\leq$ $t$, $ $ we have $ $ $g$ $\circ$ $h_{\rho}$ ($X_{i}$) = $q^{-(\rho, \beta_{i})} g$ ($X_{i}$) = $q^{-(\rho, \beta_{i})}$ $X_{i}^{(t+1)}$ = $h_{\rho}$ ($X_{i}^{(t+1)}$) = $h_{\rho}$ $\circ$ $g$ ($X_{i}$). This implies that the homomorphisms $ $ $g$ $\circ$ $h_{\rho}$ $ $ and $ $ $h_{\rho}$ $\circ$ $g$ $ $ coincide on $ $ $R$.
\item results thoroughly of 1. and 2.
\end{enumerate}
\bx

Until the end of this section, we identify the canonical generators $ $ $X_{1}$, $ $ ... $ $, $ $ $X_{t}$ $ $ of $ $ $R$ $ $ with $ $ $X_{1}^{(t+1)}$, $ $ ... $ $, $ $ $X_{t}^{(t+1)}$ $ $ respectively, so that $ $ $R$ $ $ is identified with the subalgebra $ $ $k<X_{1}^{(t+1)}$, $ $ ... $ $, $ $ $X_{t}^{(t+1)}>$ $ $ of $ $ $R_{0}^{(t+1)}$, $ $ $H$ - $Spec(R)$ is identified with the set of $ $ $H$ - invariant prime ideals of  $ $ $k<X_{1}^{(t+1)}$, $ $ ... $ $, $ $ $X_{t}^{(t+1)}>$ $ $ and the division ring $ $ $F$ = $Fract(R)$ $ $ is identified with the subalgebra of $ $ $F_{0}$ $ $ defined as follows: $ $ $F$ = $\{as^{-1} \mbox{ } | a  \in R, \mbox{ } s  \in R \setminus \{0\} \}$. $ $ ($R$ $ $ and $ $ $F$ $ $ are both $ $ $H$ - invariant.)

\begin{lem5.1.} $ $ \\
Consider an integer $ $ $m$ $ $ with $ $ $2$ $\leq$ $m$ $\leq$ $t + 1$. \\
If $ $ $X_{1}^{(m)}$, $ $ ... $ $, $ $ $X_{N}^{(m)}$ $ $ are the canonical generators of $ $ $R_{0}^{(m)}$, $ $ then $ $ $X_{1}^{(m)}$, $ $ ... $ $, $ $ $X_{t}^{(m)}$ $ $ are the canonical generators of $ $ $R^{(m)}$.
\end{lem5.1.}

\bf Proof \rm \\
We proceed by decreasing induction on $ $ $m$. \\ $ $ \\
\begin{itemize}
\item If  $ $ $m$ = $t + 1$, $ $ $X_{1}^{(t+1)}$, $ $ ... $ $, $ $ $X_{t}^{(t+1)}$ $ $ are the canonical generators of $ $ $R$ $ $ by the previous identification. So, they also are the canonical generators of $ $ $R^{(t+1)}$ $ $ by the definition of this algebra.
\item Assume that $ $ $2$ $\leq$ $m$ $\leq$ $t$ $ $ and that $ $ $X_{1}^{(m+1)}$, $ $ ... $ $, $ $ $X_{t}^{(m+1)}$ $ $ are the canonical generators of $ $ $R^{(m+1)}$. $ $ Denote (temporarily) by $ $ $X_{1}^{\prime (m)}$, $ $ ... $ $, $ $ $X_{t}^{\prime (m)}$ $ $ the canonical generators of $ $ $R^{(m)}$. $ $ Recall (see section 3.2) that they are defined as follows:
\item For each $i$ $\in$ $\llbracket  1,$ ... , $t \rrbracket $, we have
\begin{enumerate}
\item $ $ $ $ $ $ $ $ $m \leq i \mbox{ } \Rightarrow \mbox{ } X_{i}^{\prime (m)}  = \mbox{ } X_{i}^{(m+1)}$. \\
\item $$i < m  \mbox{ } \Rightarrow \mbox{ } X_{i}^{\prime (m)} = \mbox{ } X_{i}^{(m+1)} \mbox{ } \mbox{ } +  \mbox{ } \mbox{ } \displaystyle \sum_{l = 1}^{+\infty} C_{l}^{(m+1)}(X_{m}^{(m+1)})^{-l}$$
with
$$C_{l}^{(m+1)} = \mbox{ }  \displaystyle \frac{(1-q_{m})^{-l}}{[l]!_{q_{m}}} \mbox{ } \lambda_{m,i}^{-l}  \mbox{ } (\delta_{m}^{(m+1)})^{l} \mbox{ } (X_{i}^{(m+1)}).$$
\end{enumerate}
\end{itemize}

Moreover, if $ $ $h_{m}$ $ $ is the (unique) automorphism of $ $ $k<X_{1}^{(m+1)}$, $ $ ... $ $, $ $ $X_{m-1}^{(m+1)}>$ $ $ which satisfies $ $ $h_{m}(X_{i}^{(m+1)})$ = $\lambda_{m,i}X_{i}^{(m+1)}$ $ $ ($\lambda_{m,i}$ = $q^{-(\beta_{m}, \beta_{i})}$), $ $ we have $ $ $\delta_{m}^{(m+1)}(a)$ = $X_{m}^{(m+1)}a$ $-$ $h_{m}(a)X_{m}^{(m+1)}$ $ $ for any $ $ $a$ $\in$ $k<X_{1}^{(m+1)}$, $ $ ... $ $, $ $ $X_{m-1}^{(m+1)}>$. $ $ As $ $ $X_{1}^{(m)}$, $ $ ... $ $, $ $ $X_{t}^{(m)}$ $ $ are defined by the same formulas, the proof is over. \bx

Recall (see sections 3.1 and 3.2) that each prime ideal of $ $ $R_{0}^{(t+1)}$ $ $ (resp. $ $ $R$) $ $ is completely prime. So, if $ $ $\cal Q$ $ $ is any prime ideal of $ $ $R_{0}^{(t+1)}$, $ $ then $ $ $\mathcal{P}$ = $ $ $\cal Q$ $\cap$ $R$ $ $ is a prime ideal of $ $ $R$. $ $ Moreover, since $ $ $R$ = $k<X_{1}^{(t+1)}$, $ $ ... $ $, $ $ $X_{t}^{(t+1)}>$ $ $ is $H$ - invariant, we have:  $ $ ($\cal Q$ $ $ is $H$ - invariant) $\Rightarrow$ ($\mathcal{P}$ $ $ is $H$ - invariant).

\begin{lem5.1.} $ $ \\
If $ $ $\mathcal{P}$ $ $ is any $H$ - invariant prime ideal of $ $ $R$, $ $ there exists a $H$ - invariant prime ideal $ $ $\cal Q$ $ $ of $ $ $R_{0}^{(t+1)}$, $ $ such that \\
$\bullet$ $ $ $\mathcal{P}$ = $ $ $\cal Q$ $\cap$ $R$. \\
$\bullet$ $ $ $\cal Q$ $\cap$ $\{X_{t+1}^{(t+1)}$, $ $ ... $ $, $ $ $X_{N}^{(t+1)}\}$ = $\emptyset$.
\end{lem5.1.}

\bf Proof \rm \\

Set

$$\cal Q \it \mbox{ } = \mbox{ } \displaystyle \sum_{\underline{a} \mbox{ } = \mbox{ } (a_{t+1}, \mbox{ } ... \mbox{ } , \mbox{ } a_{N})}{} \mbox{ } \mathcal{P} \it (X^{(t+1)})^{\underline{a}}$$

with $ $ $\underline{a}$ $\in$ $ $ \bbN $^{N - t}$ $ $ and $ $ $(X^{(t+1)})^{\underline{a}}$ := $(X_{t+1}^{(t+1)})^{a_{t+1}}$ $ $ ... $ $ $(X_{N}^{(t+1)})^{a_{N}}$. \\ $ $ \\

Observe that, by the results recalled in 3.2, the family $((X^{(t+1)})^{\underline{a}})$ $ $ ($\underline{a} \mbox{ } = \mbox{ } (a_{t+1}, \mbox{ } ... \mbox{ } , \mbox{ } a_{N})$ $\in$ $ $ \bbN $^{N - t})$ $ $ is a basis of $ $ $R_{0}^{(t+1)}$ $ $ as a right $ $ $R$ - module. This implies that $ $ $\cal Q$ $\cap$ $R$ = $ $ $\mathcal{P}$ $ $ and, since $ $ $1$ $\notin$ $\mathcal{P}$, $ $ that each $ $ $(X^{(t+1)})^{\underline{a}}$ $ $ $\notin$ $\cal Q$ $ $ ($\underline{a}$ $\in$ $ $ \bbN $^{N - t}).$ $ $ In particular, $ $ $\cal Q$ $\cap$ $\{X_{t+1}^{(t+1)}$, $ $ ... $ $, $ $ $X_{N}^{(t+1)}\}$ = $\emptyset$. \\

For each $ $ $\underline{a}$ = $(a_{t+1}, \mbox{ } ... \mbox{ } , \mbox{ } a_{N})$ $\in$ $ $ \bbN $^{N - t}$ $ $ and each $ $ $u$ $\in$ $R$, $ $ we have  $ $ $(X^{(t+1)})^{\underline{a}}$ $u$ = $h_{\underline{a}} (u)$ $(X^{(t+1)})^{\underline{a}}$ $ $ with $h_{\underline{a}}$ = $h_{a_{t+1}}$ $\circ$ ... $\circ$ $h_{a_{N}}$ $\in$ $ $ $H$ $ $ (see section 3.2). \\
This implies that $ $ $\cal Q$ $ $ is an $ $ $R$ - bimodule. \\ $ $ \\

Since, for any  $ $ $\underline{a}$ $ $ and $ $ $\underline{b}$ $ $ in \bbN $^{N - t}$, $ $ we have $ $ $(X^{(t+1)})^{\underline{a}}$ $(X^{(t+1)})^{\underline{b}}$ = $\lambda$ $(X^{(t+1)})^{\underline{c}}$ $ $ with $ $ $\underline{c}$ = $\underline{a}$ + $\underline{b}$ and $\lambda \in k^{*}$, $ $ $\cal Q$ $ $ is a right $ $ $R_{0}^{(t+1)}$ - module. Moreover, since $(X^{(t+1)})^{\underline{a}}$ $\mathcal{P}$ = $h_{\underline{a}} (\mathcal{P})$ $(X^{(t+1)})^{\underline{a}}$ = $\mathcal{P}$ $(X^{(t+1)})^{\underline{a}}$, $ $ we obtain that $ $ $\cal Q$ $ $ is an ideal in $ $ $R_{0}^{(t+1)}$. \\ $ $ \\

Now, consider $ $ $A$, $ $ $B$ $ $ in $ $ $R_{0}^{(t+1)}$ $ $ and assume that $ $ $AB$ $\in$ $\cal Q$. $ $ As in section 4.3, denote by $ $ $\preceq$ $ $ the inverse lexicographic order in \bbN $^{N - t}$ $ $ and assume that, neither $ $ $A$, $ $ nor $ $ $B$ is in $ $ $\cal Q$. $ $ Write

$$A \mbox{ } = \mbox{ } \displaystyle \sum_{\underline{a} \mbox{ } \in F}{} \mbox{ } A_{\underline{a}} (X^{(t+1)})^{\underline{a}}$$

with $ $ $F$ $ $ a non empty finite subset of \bbN $^{N - t}$ $ $ and each $ $ $A_{\underline{a}}$ $ $ in $ $ $R$. $ $ Since, $ $ $A$ $\notin$ $\cal Q$, $ $ at least one coefficient $ $ $A_{\underline{a}}$ $ $ is not in $ $ $\mathcal{P}$. If we choose such a coefficient with $ $ $\underline{a}$ $ $ minimal, we get:

$$A \mbox{ } = \mbox{ } A_{1} \mbox{ } + \mbox{ } A_{\underline{a}} (X^{(t+1)})^{\underline{a}} \mbox{ } + \mbox{ }\displaystyle \sum_{\underline{a^{\prime}} \mbox{ } \in F^{\prime}}{} \mbox{ } A_{\underline{a^{\prime}}} (X^{(t+1)})^{\underline{a^{\prime}}}$$

with $ $ $A_{1}$ $\in$ $\cal Q$ $ $ and $ $ $\underline{a}$ $\prec$ $\underline{a^{\prime}}$ $ $ for any $ $ $\underline{a^{\prime}}$ $ $ in $ $ $F^{\prime}$. $ $ Similarly, we have

$$B \mbox{ } = \mbox{ } B_{1} \mbox{ } + \mbox{ } B_{\underline{b}} (X^{(t+1)})^{\underline{b}} \mbox{ } + \mbox{ }\displaystyle \sum_{\underline{b^{\prime}} \mbox{ } \in G^{\prime}}{} \mbox{ } B_{\underline{b^{\prime}}} (X^{(t+1)})^{\underline{b^{\prime}}}$$

with $ $ $B_{1}$ $\in$ $\cal Q$, $ $ $B_{\underline{b}}$ $ $ is not in $ $ $\mathcal{P}$ $ $ and $ $ $\underline{b}$ $\prec$ $\underline{b^{\prime}}$ $ $ for any $ $ $\underline{b^{\prime}}$ $ $ in $ $ $G^{\prime}$. \\ $ $ \\

Now, we have
$$(A \mbox{ } - \mbox{ } A_{1}) (B \mbox{ } - \mbox{ } B_{1}) \mbox{ } = \mbox{ } \lambda A_{\underline{a}}  B_{\underline{b}}(X^{(t+1)})^{\underline{c}} \mbox{ } + \mbox{ }\displaystyle \sum_{\underline{c^{\prime}} \mbox{ } \in E^{\prime}}{} \mbox{ } C_{\underline{c^{\prime}}} (X^{(t+1)})^{\underline{c^{\prime}}}$$

with $ $ $\lambda$ $ $ in $ $ $k^{*}$, $ $ $\underline{c}$ = $\underline{a}$ + $\underline{b}$, $ $  each $ $ $C_{\underline{c^{\prime}}}$ $ $ in $ $ $R$ $ $ and, for each $ $ $\underline{c^{\prime}}$ $ $ in $ $ $E^{\prime}$, $ $ $\underline{c}$ $\prec$ $\underline{c^{\prime}}$. $ $ As $(A \mbox{ } - \mbox{ } A_{1}) (B \mbox{ } - \mbox{ } B_{1})$ $\in$ $\cal Q$, $ $ this implies that $ $ $A_{\underline{a}}  B_{\underline{b}}$ $\in$ $\mathcal{P}$, $ $ which is impossible since $ $ $\mathcal{P}$ $ $ is completely prime. This implies that  $ $ $\cal Q$ $ $ is (completely) prime. \\ $ $ \\

Since $ $ $\mathcal{P}$ $ $ is $ $ $H$ - invariant and each monomial $ $ $X^{(t+1)})^{\underline{a}}$ $ $ ($\underline{a}$ $\in$ $ $ \bbN $^{N - t}$), $ $ is a $H$ - eigenvector, $ $ $\cal Q$ $ $ is $ $ $H$ - invariant. \bx

\begin{lem5.1.} $ $ \\
Consider $ $ $\mathcal{P}\it_{0}$ $\in$ $Spec(R_{0})$ $ $ and $ $ $\mathcal{P}$ $\in$ $Spec(R)$. \\
Assume that $ $ $\mathcal{P}$ = $ $ $\cal Q$ $\cap$ $R$ $ $ with $ $ $\cal Q$ = $\mathcal{P}\it_{0}^{(t+1)}$ = the canonical image of $ $ $\mathcal{P}\it_{0}$ $ $ in $Spec(R_{0}^{(t+1)})$. \\
For any $ $ $m$ $\in$ $\llbracket  2,$ ... , $t+1 \rrbracket $, $ $ we denote by $ $ $\mathcal{P}\it^{(m)}$ $ $ $($resp. $ $ $\mathcal{P}\it_{0}\it^{(m)})$ the canonical image of $ $ $\mathcal{P}$ $ $ $($resp. $ $ $\mathcal{P}\it_{0})$ $ $ in $ $ $Spec(R^{(m)})$ $ $ $($resp. $Spec(R_{0}^{(m)}))$ $ $ and we recall that, by lemma 5.1. 2, $ $ $R^{(m)}$ $ $ is the subalgebra of $ $ $R_{0}^{(m)}$ generated by it's $ $ $t$ $ $ first canonical generators $ $ $X_{1}^{(m)}$, $ $ ... $ $, $ $ $X_{t}^{(m)}$. Then

$$\mathcal{P}\it^{(m)} \mbox{ } = \mbox{ } \mathcal{P}\it_{0}^{(m)} \mbox{ } \cap \mbox{ } R^{(m)}.$$

\end{lem5.1.}

\bf Proof \rm \\
We prove this by decreasing induction on $ $ $m$. \\ $ $ \\
\begin{itemize}
\item If $ $ $m$ = $t+1$, $ $ we have $ $ $R^{(m)}$ = $R$, $ $ $\mathcal{P}\it^{(m)}$ = $\mathcal{P}$ and $ $ $\mathcal{P}\it_{0}\it^{(m)}$ = $\cal Q$. $ $ So, we have the required equality by assumption.
\item Assume $ $ $m$ $\leq$ $t$ $ $ and $ $ $\mathcal{P}\it^{(m+1)} \mbox{ } = \mbox{ } \mathcal{P}\it_{0}\it^{(m+1)} \mbox{ } \cap \mbox{ } R^{(m+1)}.$ \\ $ $ \\
\begin{itemize}
\item Assume that $ $ $X_{m}^{(m+1)}$ $\notin$ $\mathcal{P}\it_{0}\it^{(m+1)}$, $ $ so that $ $ $X_{m}^{(m+1)}$ $\notin$ $\mathcal{P}\it^{(m+1)}$. \\
Set $ $ $S_{m}$ = $\{(X_{m}^{(m+1)})^{d} \mbox{ } | \mbox{ } d \in$ $ $ \bbN $\}$ $ $ and recall (sections 3.2 and 3.3) that \\ $ $ \\
\begin{itemize}
\item $ $ $R^{(m)}$ $S_{m}^{-1}$ = $R^{(m+1)}$ $S_{m}^{-1}$.
\item $ $ $\mathcal{P}\it^{(m)}$ = $ $ $\mathcal{P}\it^{(m+1)}$ $S_{m}^{-1}$ $ $ $\cap$ $ $ $R^{(m)}$.
\item $ $ $\mathcal{P}\it_{0}\it^{(m)}$ = $ $ $\mathcal{P}\it_{0}\it^{(m+1)}$ $S_{m}^{-1}$ $ $ $\cap$ $ $ $R_{0}^{(m)}$.
\end{itemize}
 $ $ \\
As $ $ $R^{(m)}$ $ $ $\subseteq$ $ $ $R_{0}^{(m)}$ $ $ and $ $ $\mathcal{P}\it^{(m+1)} \mbox{ } \subseteq \mbox{ } \mathcal{P}\it_{0}\it^{(m+1)}$, $ $ we have
$$\mathcal{P}\it^{(m)} \mbox{ } = \mbox{ } \mathcal{P}\it^{(m+1)} S_{m}^{-1} \mbox{ } \cap \mbox{ } R^{(m)} \mbox{ } \subseteq \mbox{ } \mathcal{P}\it_{0}\it^{(m+1)} S_{m}^{-1} \mbox{ } \cap \mbox{ } R_{0}^{(m)} \mbox{ } = \mbox{ } \mathcal{P}\it_{0}\it^{(m)}.$$

This implies that

$$\mathcal{P}\it^{(m)} \mbox{ } \subseteq \mbox{ } \mathcal{P}\it_{0}\it^{(m)} \mbox{ } \cap \mbox{ } R^{(m)}.$$

Now, consider some $ $ $u$ $\in$ $\mathcal{P}\it_{0}\it^{(m)} \mbox{ } \cap \mbox{ } R^{(m)}$. \\
We have $ $ $\mathcal{P}\it_{0}\it^{(m)}$ $ $ $\subseteq$ $\mathcal{P}\it_{0}\it^{(m+1)}$ $S_{m}^{-1}$ $ $ $\Rightarrow$ $ $ $u$ $(X_{m}^{(m+1)})^{d_{1}}$ $\in$ $\mathcal{P}\it_{0}\it^{(m+1)}$ $ $ with $ $ $d_{1}$ $\in$ $ $ \bbN. \\
We have $ $ $R^{(m)}$ $ $ $\subseteq$ $R^{(m+1)}$ $S_{m}^{-1}$ $ $ $\Rightarrow$ $ $ $u$ $(X_{m}^{(m+1)})^{d_{2}}$ $\in$ $R^{(m+1)}$ $ $ with $ $ $d_{2}$ $\in$ $ $ \bbN. \\
Since $ $ $m$ $\leq$ $t$, $ $ $X_{m}^{(m+1)}$ $ $ is in $ $ $R^{(m+1)}$ and so, with $ $ $d$ = max $(d_{1}$, $d_{2})$, we obtain that $ $ $u$ $(X_{m}^{(m+1)})^{d}$ $\in$ \\ $\mathcal{P}\it_{0}\it^{(m+1)}$ $ $ $\cap$ $ $ $R^{(m+1)}$ = $\mathcal{P}\it^{(m+1)}$ $ $ $\Rightarrow$ $ $ $u$ $\in$ $\mathcal{P}\it^{(m+1)}$ $S_{m}^{-1}$ $ $ $\cap$ $ $ $R^{(m)}$ = $\mathcal{P}\it^{(m)}$. $ $ This implies that

$$\mathcal{P}\it^{(m)} \mbox{ } = \mbox{ } \mathcal{P}\it_{0}\it^{(m)} \mbox{ } \cap \mbox{ } R^{(m)}.$$
\item Assume that $ $ $X_{m}^{(m+1)}$ $\in$ $\mathcal{P}\it_{0}\it^{(m+1)}$, $ $ so that, since $ $ $m$ $\leq$ $t$,

$$X_{m}^{(m+1)} \mbox{ } \in \mbox{ } \mathcal{P}\it_{0}^{(m+1)} \mbox{ } \cap \mbox{ } R^{(m+1)} \mbox{ } = \mbox{ } \mathcal{P}\it^{(m+1)}.$$  \\

Denote by

$$g : \mbox{ } R^{(m)} \mbox{ } \rightarrow \mbox{ } \frac{R^{(m+1)}}{\mathcal{P}\it^{(m+1)}}$$

the only algebra homomorphism which transforms each $ $ $X_{i}^{(m)}$ $ $ ($1$ $\leq$ $i$ $\leq$ $t$) $ $ in $ $ $g (X_{i}^{(m)}$) = $X_{i}^{(m+1)}$ $ + $ $\mathcal{P}\it^{(m+1)}$ $ $ (the canonical image of $ $ $X_{i}^{(m+1)}$ $ $ in $ $ $\displaystyle \frac{R^{(m+1)}}{\mathcal{P}\it^{(m+1)}}$). $ $ We recall (section 3.3) that

$$\mathcal{P}\it^{(m)} \mbox{ } = \mbox{ } ker ( g ).$$

Denote by

$$g_{0} : \mbox{ } R_{0}^{(m)} \mbox{ } \rightarrow \mbox{ } \frac{R_{0}^{(m+1)}}{\mathcal{P}\it_{0}\it^{(m+1)}}$$

the only algebra homomorphism which transforms each $ $ $X_{i}^{(m)}$ $ $ ($1$ $\leq$ $i$ $\leq$ $N$) $ $ in $ $ $g_{0} (X_{i}^{(m)}$) = $X_{i}^{(m+1)}$ $ + $ $\mathcal{P}\it_{0}\it^{(m+1)}$ $ $ (the canonical image of $ $ $X_{i}^{(m+1)}$ $ $ in $ $ $\displaystyle \frac{R_{0}^{(m+1)}}{\mathcal{P}\it_{0}\it^{(m+1)}}$), $ $ so that

$$\mathcal{P}\it_{0}\it^{(m)} \mbox{ } = \mbox{ } ker ( g_{0} ).$$

The canonical homomorphism $ $ $\displaystyle R_{0}^{(m+1)} \mbox{ } \rightarrow \mbox{ } \frac{R_{0}^{(m+1)}}{\mathcal{P}\it_{0}\it^{(m+1)}}$ $ $ restricted to $ $ $R^{(m+1)}$ $ $ has it's kernel equal to $ $ $\mathcal{P}\it_{0}\it^{(m+1)}$ $ $ $\cap$ $ $ $R^{(m+1)}$ = $\mathcal{P}\it^{(m+1)}$. $ $ So, it induces an injective homomorphism

$$\epsilon : \mbox{ } \frac{R^{(m+1)}}{\mathcal{P}\it^{(m+1)}} \mbox{ } \rightarrow \mbox{ } \frac{R_{0}^{(m+1)}}{\mathcal{P}\it_{0}\it^{(m+1)}}$$

which transforms each $ $ $X_{i}^{(m+1)}$ $ + $ $\mathcal{P}\it^{(m+1)}$ $ $ ($1$ $\leq$ $i$ $\leq$ $t$) $ $ in $ $  $X_{i}^{(m+1)}$ $ + $ $\mathcal{P}\it_{0}\it^{(m+1)}$. So,

$$\epsilon \circ g : \mbox{ } R^{(m)} \mbox{ } \rightarrow \mbox{ } \frac{R_{0}^{(m+1)}}{\mathcal{P}\it_{0}\it^{(m+1)}}$$
is an algebra homomorphism which transforms each $ $ $X_{i}^{(m)}$ $ $ ($1$ $\leq$ $i$ $\leq$ $t$) $ $ in $ $  $X_{i}^{(m+1)}$ $ + $ $\mathcal{P}\it_{0}\it^{(m+1)}$ = $g_{0}$ ($X_{i}^{(m)}$). $ $ This implies that $ $ $\epsilon \circ g$ $ $ is the restriction of $ $ $g_{0}$ $ $ to $ $ $R^{(m)}$ $ $ and so, that $ $ $Ker$ ($\epsilon \circ g$) = $\mathcal{P}\it_{0}\it^{(m)}$ $ $ $\cap$ $ $ $R^{(m)}$. $ $ Since $ $ $\epsilon$ $ $ is injective, we also have $ $ $Ker (\epsilon \circ g$) = $Ker$ ($g$) = $\mathcal{P}\it^{(m)}$ and the proof is over.
\end{itemize}
\end{itemize}
\bx

\newtheorem{prop5.1.}{Proposition 5.1.}

\begin{prop5.1.} $ $ \\
Consider a diagram $ $ $\Delta$ $ $ with respect to \rm(\ref{eqn:expw}) (\it ie. a subset of   $\llbracket  1,  \mbox{ } ...  \mbox{ } ,  \mbox{ } t \rrbracket )$, $ $ so that $ $ $\Delta$ $ $ is also a diagram with respect to \rm(\ref{eqn:expw0})\it. Then $\Delta$ $ $ is admissible with respect to \rm(\ref{eqn:expw})\it $ $ if and only if $ $ $\Delta$ $ $ is admissible with respect to \rm(\ref{eqn:expw0}).
\end{prop5.1.}

\bf Proof \rm \\ $ $ \\
\begin{itemize}
\item Assume that $ $ $\Delta$ $ $ is admissible with respect to (\ref{eqn:expw0}). This means (see section 3.3) that there exists some $ $ $\mathcal{P}\it_{0}$ $ $ in $ $ $Spec(R_{0})$ $ $ such that it's canonical image $ $ $\mathcal{P}\it_{0}^{(2)}$ $ $ in $ $ $Spec(R_{0}^{(2)})$ $ $ verifies:

$$\mathcal{P}\it_{0}^{(2)} \mbox{ } \cap \mbox{ } \{Z_{1}, \mbox{ } ... \mbox{ } , \mbox{ } Z_{N}\} \mbox{ } = \mbox{ } \{Z_{i} \mbox{ } | \mbox{ } i \mbox{ } \in  \mbox{ } \Delta\}.$$

Denote by $ $ $\cal Q$ = $\mathcal{P}\it_{0}^{(t+1)}$ $ $ the canonical image of $ $ $\mathcal{P}\it_{0}$ $ $ in $ $ $Spec(R_{0}^{(t+1)})$ $ $ and set $ $ $\mathcal{P}$ = $\cal Q$ $ $ $\cap$ $ $ $R$. $ $ We know that $ $ $\mathcal{P}$ $\in$ $Spec(R)$ and, by lemma 5.1. 4, it's canonical image $ $ $\mathcal{P}^{(2)}$ $ $ in $ $ $Spec(R^{(2)})$ $ $ verifies

$$\mathcal{P}^{(2)} \mbox{ } = \mbox{ } \mathcal{P}\it_{0}^{(2)} \mbox{ } \cap \mbox{ } R^{(2)}.$$

So,

$$\mathcal{P}^{(2)} \mbox{ } \cap \mbox{ } \{Z_{1}, \mbox{ } ... \mbox{ } , \mbox{ } Z_{t}\} \mbox{ } = \mbox{ } \mathcal{P}\it_{0}^{(2)} \mbox{ } \cap \mbox{ } R^{(2)} \mbox{ } \cap \mbox{ } \{Z_{1}, \mbox{ } ... \mbox{ } , \mbox{ } Z_{t}\} \mbox{ } = \mbox{ } \mathcal{P}\it_{0}^{(2)} \mbox{ } \cap \mbox{ } \{Z_{1}, \mbox{ } ... \mbox{ } , \mbox{ } Z_{t}\} \mbox{ }.$$

Since $ $ $\Delta$ $ $ is a subset of   $\llbracket  1,  \mbox{ } ...  \mbox{ } ,  \mbox{ } t \rrbracket $), $ $ this implies that

$$\mathcal{P}^{(2)} \mbox{ } \cap \mbox{ } \{Z_{1}, \mbox{ } ... \mbox{ } , \mbox{ } Z_{t}\} \mbox{ } = \mbox{ } \{Z_{i} \mbox{ } | \mbox{ } i \mbox{ } \in  \mbox{ } \Delta\}.$$

So $ $ $\Delta$ $ $ is admissible with respect to (\ref{eqn:expw}).
\item Assume that $ $ $\Delta$ $ $ is admissible with respect to (\ref{eqn:expw}). This implies (see section 3.3) that there exists some $ $ $\mathcal{P}$ $ $ in $ $ $H$ $-$ $Spec(R)$ $ $ whose canonical image $ $ $\mathcal{P}^{(2)}$ $ $ in $ $ $Spec(R^{(2)})$ $ $ is the ideal of $ $ $R^{(2)}$ $ $ generated by $ $ $\{Z_{i} \mbox{ } | \mbox{ } i \mbox{ } \in  \mbox{ } \Delta\}$ $ $ ($\mathcal{P}^{(2)}$ = $0$ $ $ if $ $ $\Delta$ = $\emptyset$). \\
Since $ $ $\mathcal{P}$ $ $ is $ $ $H$ $-$ invariant, it results from lemma 5.1. 3 that there exists a ($H$ $-$ invariant) prime ideal $ $ $\cal Q$ $ $ of $ $ $R_{0}^{(t+1)}$, $ $ such that \\
\begin{itemize}
\item $\mathcal{P}$ = $ $ $\cal Q$ $\cap$ $R$.
\item $\cal Q$ $\cap$ $\{X_{t+1}^{(t+1)}$, $ $ ... $ $, $ $ $X_{N}^{(t+1)}\}$ = $\emptyset$.
\end{itemize}
Now, by corollary 4.1. 1, there exists $ $ $\mathcal{P}\it_{0}$ $ $ in $ $ $Spec(R_{0})$ $ $ whose canonical image $ $ $\mathcal{P}\it_{0}^{(t+1)}$ $ $ in $ $ $Spec(R_{0}^{(t+1)})$ $ $ is equal to $ $ $\cal Q$ $ $ and whose canonical image $ $ $\mathcal{P}\it_{0}^{(2)}$ $ $ in $ $ $Spec(R_{0}^{(2)})$ $ $ does not meet the set $ $ $\{Z_{t+1}$, $ $ ... $ $, $ $ $Z_{N}\}$. \\ As over, by lemma 5.1. 4, we have

$$\mathcal{P}^{(2)} \mbox{ } = \mbox{ } \mathcal{P}\it_{0}^{(2)} \mbox{ } \cap \mbox{ } R^{(2)}.$$

Since   $ $ $\mathcal{P}\it_{0}^{(2)}$ $ $ does not meet the set $ $ $\{Z_{t+1}$, $ $ ... $ $, $ $ $Z_{N}\}$, $ $ we have

$$\mathcal{P}\it_{0}^{(2)} \mbox{ } \cap \mbox{ } \{Z_{1}, \mbox{ } ... \mbox{ } , \mbox{ } Z_{N}\} \mbox{ } = \mbox{ } \mathcal{P}\it_{0}^{(2)} \mbox{ } \cap \mbox{ } \{Z_{1}, \mbox{ } ... \mbox{ } , \mbox{ } Z_{t}\} \mbox{ }.$$

So, since $ $ $\{Z_{1}, \mbox{ } ... \mbox{ } , \mbox{ } Z_{t}\}$ $\subseteq$ $R^{(2)}$, $ $ we can write

$$\mathcal{P}\it_{0}^{(2)} \mbox{ } \cap \mbox{ } \{Z_{1}, \mbox{ } ... \mbox{ } , \mbox{ } Z_{N}\} \mbox{ } = \mbox{ } \mathcal{P}\it_{0}^{(2)}  \mbox{ } \cap \mbox{ } R^{(2)} \mbox{ } \cap \mbox{ } \{Z_{1}, \mbox{ } ... \mbox{ } , \mbox{ } Z_{t}\} \mbox{ }$$
$$= \mbox{ } \mathcal{P}^{(2)} \mbox{ } \cap \mbox{ } \{Z_{1}, \mbox{ } ... \mbox{ } , \mbox{ } Z_{N}\} \mbox{ } = \mbox{ } \{Z_{i} \mbox{ } | \mbox{ } i \mbox{ } \in  \mbox{ } \Delta\}.$$

This proves that $ $ $\Delta$ $ $ is admissible with respect to (\ref{eqn:expw0}).
\end{itemize}
\bx

\subsection{The case $w$ = $w_{0}$.}
Assume (in this subsection only) that $ $ $w$ = $w_{0}$. \\ $ $ \\

\newtheorem{lem5.2.}{Lemma 5.2.}

\begin{lem5.2.} $ $ \\
For any diagram $ $ $\Delta$ $ $ $($with respect to \rm(\ref{eqn:expw})) \it we have:
$$\Delta \mbox{ } \rm positive \it \mbox{ } \Rightarrow \mbox{ } \Delta \mbox{ } \rm admissible.$$
\end{lem5.2.}

\bf Proof \rm \\
As in section 3.3, we denote by $ $ $\mathcal{P}^{(2)}_{\Delta}$ $ $ the ideal of $ $ $\overline{R}$ = $R^{(2)}$ $ $ generated by the canonical generators $ $ $Z_{i}$ $ $ with $ $ $i$ $\in$ $\Delta$ $ $ ($\mathcal{P}^{(2)}_{\Delta}$ = $0$ $ $ if $ $ $\Delta$ = $\emptyset$). $ $ Recall (see section 3.3) that $ $ $\mathcal{P}^{(2)}_{\Delta}$ $\in$ $H$ $-$ $Spec(R^{(2)})$. \\
Consider an integer $ $ $m$ with $ $ $2$ $\leq$ $m$ $\leq$ $t+1$. $ $ We first prove, by induction on $ $ $m$, $ $ that there exists a prime ideal ideal $ $ $\mathcal{P}\it^{(m)}$ $\in$ $Spec(R^{(m)})$ $ $ such that $ $ $\mathcal{P}^{(2)}_{\Delta}$ $ $ is the canonical image of $ $ $\mathcal{P}\it^{(m)}$ $ $ in $ $ $Spec(R^{(2)})$. \\
\begin{itemize}
\item If $ $ $m$ = $2$, $ $ we just have to choose $ $ $\mathcal{P}^{(2)}$ = $\mathcal{P}^{(2)}_{\Delta}$.
\item Assume that, for some $ $ $m$ with $ $ $2$ $\leq$ $m$ $\leq$ $t$, $ $ there exists a prime ideal ideal $ $ $\mathcal{P}\it^{(m)}$ $\in$ $Spec(R^{(m)})$ such that $ $ $\mathcal{P}^{(2)}_{\Delta}$ $ $ is the canonical image of $ $ $\mathcal{P}\it^{(m)}$ $ $ in $ $ $Spec(R^{(2)})$. $ $ Since $ $ $\Delta \mbox{ }$ $ $ is positive, we deduce from proposition 4.4. 1 that $ $ $\mathcal{P}\it^{(m)}$ = $\phi_{m} (\mathcal{P}\it^{(m+1)})$ $ $ with $ $ $\mathcal{P}\it^{(m+1)}$ $ $ in $ $ $Spec(R^{(m+1)})$. \\
So, $ $ $\mathcal{P}^{(2)}_{\Delta}$ $ $ is the canonical image of $ $ $\mathcal{P}\it^{(m+1)}$ $ $ in $ $ $Spec(R^{(2)})$, $ $ and our affirmation is proved. \\
In particular, there exists $ $ $\mathcal{P}$ $\in$ $Spec(R^{(t+1)})$ = $Spec(R)$ $ $ such that $ $ $\mathcal{P}^{(2)}_{\Delta}$ $ $ is the canonical image of $ $ $\mathcal{P}$ $ $ in $ $ $Spec(R^{(2)})$. $ $ This means (see section 3.3) that $ $ $\Delta$ $ $ is admissible.
\end{itemize}
 \bx

\newtheorem{prop5.2.}{Proposition 5.2.}

\begin{prop5.2.} $ $ \\
The positive diagrams $($with respect to \rm(\ref{eqn:expw})) \it coincide with the admissible diagrams $($with respect to \rm (\ref{eqn:expw})).
\end{prop5.2.}

\bf Proof \rm \\
By lemma 5.2 1, the set of positive diagrams is contained in the set of admissible diagrams. By proposition 2.3. 1, the number of positive diagrams is equal to the number of $ $ $u$ $ $ in the Weyl group $ $ $W$ $ $ such that $ $ $u$ $\leq$ $w$ $ $ (where $ $ $\leq$ $ $ denotes the Bruhat order). Since we assume that $ $ $w$ = $w_{0}$, $ $ we get that the number of positive diagrams is equal to the cardinal $ $ $|W|$ $ $ of $ $ $W$. $ $ Since the set of admissible diagrams is in one to one correspondence with $ $ $H$ $-$ $Spec(R)$ = $H$ $-$ $Spec(U^{+})$, $ $ the number of admissible diagrams does not depend on the reduced decomposition of $ $ $w$. $ $ In \cite{M}, Antoine M\'eriaux gives a precise description of admissible diagrams, for each type of the simple complex Lie algebra $ $ $\g$, $ $ and for a particular reduced decomposition of $ $ $w_{0}$ $ $ so that he can compute their number and check that it is precisely equal to $ $ $|W|$. $ $ This implies that the positive diagrams coincide with the admissible diagrams. \bx

Let us observe that, in this proof, we use a result of Antoine M\'eriaux to prove that the number of admissible diagrams (which is also the cardinal of $ $ $H$ $-$ $Spec(U^{+})$) $ $ is equal to $ $ $|W|$. $ $ This equality can also be obtained using results of M. Gorelik, N. Andruskiewitsch and F. Dumas (\cite{G} and \cite{AnD}) or T. J. Hodges and T. Levasseur and M. Toro \cite{HLT} but this should require some (minor) restricted assumptions on the choice of the ground field $ $ $k$ $ $ (char($k$) = $0$) or on the parameter $ $ $q$ $ $ ($q$ $ $ transcendent).

\subsection{The general case.}

We can now prove, in the general case (ie. $ $ $w$ $ $ is not any longer assumed to be the longest element of $ $ $W$):

\newtheorem{th5.3.}{Theorem 5.3.}

\begin{th5.3.} $ $ \\
The positive diagrams $($with respect to \rm (\ref{eqn:expw})) \it coincide with the admissible diagrams $($with respect to \rm (\ref{eqn:expw})).
\end{th5.3.}

\bf Proof \rm \\

Consider a diagram $ $ $\Delta$ $ $ with respect to (\ref{eqn:expw}) (ie. a subset of   $\llbracket  1,  \mbox{ } ...  \mbox{ } ,  \mbox{ } t \rrbracket $), $ $ so that $ $ $\Delta$ $ $ is also a diagram with respect to (\ref{eqn:expw0}). \\
\begin{itemize}
\item Assume that $ $ $\Delta$ $ $ is positive (with respect to (\ref{eqn:expw})). By proposition 2.3. 2, $ $ $\Delta$ $ $ is positive with respect to (\ref{eqn:expw0}) and, by proposition 5.2. 1, $\Delta$ $ $ is admissible with respect to (\ref{eqn:expw0}). Now, by proposition 5.1. 1, $\Delta$ $ $ is admissible with respect to (\ref{eqn:expw}).
\item Assume that $ $ $\Delta$ $ $ is admissible (with respect to (\ref{eqn:expw})). By proposition 5.1. 1, $ $ $\Delta$ $ $ is admissible with respect to (\ref{eqn:expw0}) and, by proposition 5.2. 1, $\Delta$ $ $ is positive with respect to (\ref{eqn:expw0}). Now, by proposition 2.3. 2, $\Delta$ $ $ is positive with respect to (\ref{eqn:expw}).
\end{itemize}
\bx

\newtheorem{cor5.3.}{Corollary 5.3.} 

\begin{cor5.3.} $ $ 
$($All the diagrams are with respect to \rm (\ref{eqn:expw}).$)$ \it 
\begin{enumerate}

\item The map $ $ $\zeta$ : $ $ $ $ $\Delta$ = $\{j_{1}$, $...$ , $j_{s}\}$ $ $ $\mapsto$ $ $ $u$ = $w^{\Delta}$ = $s_{\alpha_{j_{1}}}$ $\circ$ $...$ $\circ$ $s_{\alpha_{j_{s}}}$ $ $ $ $ is a bijection from the set of admissible diagrams onto the set $ $ $ $ $\{u \in W \mbox{ } | \mbox{ } u \mbox{ } \leq \mbox{ } w \}$ $ $ $ $ $ $ 
$(\zeta (\emptyset)$ = $Id)$.

\item Consider an admissible diagram $ $ $\Delta$ = $\{j_{1}$, $...$ , $j_{s}\}$ $ $ and some integer $ $ $i$ $\in$ $\llbracket  1,  \mbox{ } ...  \mbox{ } ,  \mbox{ } t \rrbracket $. $ $ Set $ $ $\Delta$ $\cap$ $\llbracket i + 1,  \mbox{ } ...  \mbox{ } ,  \mbox{ } t \rrbracket $ = $\{j_{c}$, $...$ , $j_{s}\}$ $ $ $(1 \mbox{ } \leq \mbox{ } c \mbox{ } \leq \mbox{ } s)$. $ $ Then the expression $ $ $s_{\alpha_{i}}$ $\circ$ $s_{\alpha_{j_{c}}}$ $\circ$ $...$ $\circ$ $s_{\alpha_{j_{s}}}$ $ $ is reduced. In particular, $s_{\alpha_{j_{1}}}$ $\circ$ $...$ $\circ$ $s_{\alpha_{j_{s}}}$ $ $ is a reduced expression of $ $ $\zeta (\Delta)$. 

\item Consider some $ $ $u$ $\in$ $W$ $ $ with $ $ $u$ $\leq$ $w$. $ $  Then, the only admissible diagram $ $ $\Delta$ $ $ such that $ $ $\zeta (\Delta)$ = $u$ is recursively defined as follows: \\
$\bullet$ $ $ 
\begin{center}
$1$ $\in$ $\Delta$ $ $ $\Leftrightarrow$ $ $ $l(s_{\alpha_{1}}$ $\circ$ $u)$ = $l(u)$ $-$ $1$ $ $ $\Leftrightarrow$ $ $ $u^{-1}(\alpha_{1})$ $ $ is a negative root. 
\end{center}
$\bullet$ $ $ Consider some integer $ $ $i$ $\in$ $\llbracket  1,  \mbox{ } ...  \mbox{ } ,  \mbox{ } t - 1 \rrbracket $ $ $ and assume that $ $ $\Delta$ $\cap$ $\llbracket 1,  \mbox{ } ...  \mbox{ } ,  \mbox{ } i \rrbracket $ = $\{j_{1}$, $...$ , $j_{d}\}$. \\ Set $ $ $u_{i}$ = $s_{\alpha_{j_{d}}}$ $\circ$ $...$ $\circ$ $s_{\alpha_{j_{2}}}$ $\circ$ $s_{\alpha_{j_{1}}}$ $\circ$ $u$ $ $ $(u_{i}$ = $u$ $ $ if $ $ $\Delta$ $\cap$ $\llbracket 1,  \mbox{ } ...  \mbox{ } ,  \mbox{ } i \rrbracket $  = $\emptyset)$. $ $ Then 
\begin{center}
$ $ $i + 1$ $\in$ $\Delta$ $ $ $\Leftrightarrow$ $ $ $l(s_{\alpha_{i+1}}$ $\circ$ $u_{i})$ = $l(u_{i})$ $-$ $1$ $ $ $\Leftrightarrow$ $ $ $u_{i}^{-1}(\alpha_{i + 1})$ $ $ is a negative root. 
\end{center}

\end{enumerate}

\end{cor5.3.}

\bf Proof \rm \\

As the admissible diagrams coincide with the positive diagrams, 1. is proposition 2.3. 1 (assertion 2.) and 2. is proposition 2.2. 2. Now, let us prove 3. \\ $ $ \\

$\bullet$ $ $ 
Assume that $ $ $1$ $\in$ $\Delta$, $ $ so that $ $ $\Delta$ = $\{j_{1}$, $...$ , $j_{s}\}$ $ $ with $ $ $j_{1}$ = $1$. $ $ As $ $ $u$ = $\zeta (\Delta)$ = $s_{\alpha_{j_{1}}}$ $\circ$ $...$ $\circ$ $s_{\alpha_{j_{s}}}$, $ $ we have $ $ $l(s_{\alpha_{1}}$ $\circ$ $u)$ = $l(s_{\alpha_{j_{2}}}$ $\circ$ $...$ $\circ$ $s_{\alpha_{j_{s}}})$ = $s$ $-$ $1$ = $l(u)$ $-$ $1$. \\
Conversely, assume that  $ $ $l(s_{\alpha_{1}}$ $\circ$ $u)$ = $l(u)$ $-$ $1$. $ $ If $ $ $1$ $\notin$ $\Delta$, $ $ then $ $ $\Delta$ $\cap$ $\llbracket 2,  \mbox{ } ...  \mbox{ } ,  \mbox{ } t \rrbracket $ = $\Delta$ $ $ and, by 2., $ $ $l(s_{\alpha_{1}}$ $\circ$ $u)$ = $l(u)$ $+$ $1$. $ $ So, we obtain a contradiction. \\ Now, we have classically: $ $ $l(s_{\alpha_{1}}$ $\circ$ $u)$ = $l(u)$ $-$ $1$ $ $ $\Leftrightarrow$ $ $ $l(u^{-1}$ $\circ$ $s_{\alpha_{1}})$  = $l(u^{-1})$ $-$ $1$ $ $ $\Leftrightarrow$ $ $  $u^{-1}(\alpha_{1})$ $ $ is a negative root. 
\\ $ $ \\

$\bullet$ $ $ 
Assume that $ $ $i + 1$ $\in$ $\Delta$, $ $ so that $ $ $\Delta$ $\cap$ $\llbracket i + 1,  \mbox{ } ...  \mbox{ } ,  \mbox{ } t \rrbracket $ = $\{j_{d+1}$, $...$ , $j_{s}\}$ $ $ with $ $ $j_{d+1}$ = $i+1$. $ $ As  $ $ $u_{i}$ = $s_{\alpha_{j_{d}}}$ $\circ$ $...$ $\circ$ $s_{\alpha_{j_{2}}}$ $\circ$ $s_{\alpha_{j_{1}}}$ $\circ$ $u$ = $s_{\alpha_{j_{d+1}}}$ $\circ$ $...$ $\circ$ $s_{\alpha_{j_{s}}}$, $ $ we have $ $ $l(s_{\alpha_{i+1}}$ $\circ$ $u)$ = $l(s_{\alpha_{j_{d+2}}}$ $\circ$ $...$ $\circ$ $s_{\alpha_{j_{s}}})$ = $s$ $-$ $(d+1)$ = $l(u_{i})$ $-$ $1$. \\
Conversely, assume that  $ $ $l(s_{\alpha_{i+1}}$ $\circ$ $u_{i})$ = $l(u_{i})$ $-$ $1$. $ $ If $ $ $i + 1$ $\notin$ $\Delta$, $ $ then $ $ $\Delta$ $\cap$ $\llbracket i + 2,  \mbox{ } ...  \mbox{ } ,  \mbox{ } t \rrbracket $ = $\{j_{d + 1}$, $...$ , $j_{s}\}$ $ $ $ $ and, by 2., $ $ $l(s_{\alpha_{i + 1}}$ $\circ$ $u_{i})$ = $l(u_{i})$ $+$ $1$. $ $ So, we still have a contradiction and we conclude as over. \bx

From corollary 5.3. 1 (assertions 1. and 2.) and lemma 2.2. 1 (assertion 2.), we deduce:

\begin{cor5.3.} $ $ \\
\begin{enumerate}

\item The map $ $ $\zeta^{\prime}$ : $ $ $ $ $\Delta$ = $\{j_{1}$, $...$ , $j_{s}\}$ $ $ $\mapsto$ $ $ $u^{\prime}$ = $v^{\Delta}$ = $s_{\alpha_{j_{s}}}$ $\circ$ $...$ $\circ$ $s_{\alpha_{j_{1}}}$ $ $ $ $ is a bijection from the set of admissible diagrams onto the set $ $ $ $ $\{u \in W \mbox{ } | \mbox{ } u \mbox{ } \leq \mbox{ } v \mbox{ } = \mbox{ } w^{-1} \}$ $ $ $ $ 
$(\zeta^{\prime} (\Delta)$ = $(\zeta (\Delta))^{-1}, \mbox{ } \mbox{ }\zeta^{\prime} (\emptyset)$ = $Id)$. 

\item Consider an admissible diagram $ $ $\Delta$ = $\{j_{1}$, $...$ , $j_{s}\}$ $ $ and some integer $ $ $i$ $\in$ $\llbracket  1,  \mbox{ } ...  \mbox{ } ,  \mbox{ } t \rrbracket $. $ $ Set $ $ $\Delta$ $\cap$ $\llbracket i + 1,  \mbox{ } ...  \mbox{ } ,  \mbox{ } t \rrbracket $ = $\{j_{c}$, $...$ , $j_{s}\}$ $ $ $(1 \mbox{ } \leq \mbox{ } c \mbox{ } \leq \mbox{ } s)$. $ $ Then the expression $ $ $s_{\alpha_{j_{s}}}$ $\circ$ $...$ $\circ$ $s_{\alpha_{j_{c}}}$ $\circ$ $s_{\alpha_{i}}$ $ $ is reduced. In particular, $s_{\alpha_{j_{s}}}$ $\circ$ $...$ $\circ$ $s_{\alpha_{j_{1}}}$ $ $ is a reduced expression of $ $ $\zeta^{\prime} (\Delta)$. 

\end{enumerate}

\end{cor5.3.}

If $ $ $w$ $ $ and it's reduced decomposition are chosen as in the example of section 2.1, we know (proposition 2.1. 1) that $ $ $\upw$ $ $ is the quantum matrices algebra $ $ $O_{q}(M_{p,m}(k))$ $ $ with $ $ $m$ = $n-p+1$ and with the same canonical generators (as $ $ $\upw$). $ $ So, by \cite{CC}, the admissible diagrams are the $ $ \Le $ $ - diagrams. In this case, corollary 5.3. 2 has also been proved (with quite different methods) by A. Postnikov (\cite{P}, theorem 19.1.) and by T. Lam and L. Williams (\cite{LW}, theorem 5.3.). \\ $ $ \\

In the general case, using the bijection between the set of admissible diagrams and $ $ $H$ $-$ $Spec(\upw)$ $ $ recalled in section 3.3, the two bijections of proposition 2.3. 1 and the coincidence of admissible and positive diagrams (theorem 5.3. 1), we get

\begin{cor5.3.} $ $ \\
1. $ $ There exists a natural bijection from the set $ $ $\{u \in W \mbox{ } | \mbox{ } u \mbox{ } \leq \mbox{ } v \mbox{ } = \mbox{ } w^{-1}\}$ $ $ onto $ $ $H$ $-$ $Spec(\upw)$. $ $ It is defined by
$$v^{\Delta} \mbox{ } \mapsto \mbox{ } \mathcal{P}\it_{\Delta}$$
where $ $ $\Delta$ $ $ describes the set of admissible diagrams. \\
2. $ $ There exists a natural bijection from the set $ $ $\{u \in W \mbox{ } | \mbox{ } u \mbox{ } \leq \mbox{ } w\}$ $ $ onto $ $ $H$ $-$ $Spec(\upw)$. $ $ It is defined by
$$w^{\Delta} \mbox{ } \mapsto \mbox{ } \mathcal{P}\it_{\Delta}$$
where $ $ $\Delta$ $ $ describes the set of admissible diagrams.
\end{cor5.3.}

When $ $ $\upw$ $ $ is the quantum matrices algebra $ $ $O_{q}(M_{p,m}(k))$ $ $ as over, S. Launois \cite{L} has constructed (with different methods) a bijection from $ $ $\{u \in W \mbox{ } | \mbox{ } u \mbox{ } \leq \mbox{ } v\}$ $ $ onto $ $ $H$ $-$ $Spec(O_{q}(M_{p,m}(k)))$ $ $ which, moreover, preserves the ordering (where the Weyl group is provided with the Bruhat order and  $ $ $H$ $-$ $Spec(O_{q}(M_{p,m}(k)))$ $ $ is provided with the inclusion). So, it seems natural to ask the following questions:

\newtheorem{q5.3.}{Question 5.3.}

\begin{q5.3.} $ $ \\
Assume $ $ $w$ $ $ and it's reduced decomposition chosed as in the example of section 2.1. \\
Does the S. Launois bijection from $ $ $\{u \in W \mbox{ } | \mbox{ } u \mbox{ } \leq \mbox{ } v\}$ $ $ onto $H$ $-$ $Spec(O_{q}(M_{p,m}(k)))$ = $ $ $H$ $-$ $Spec(\upw)$ $ $ coincide with the first bijection of corollary 5.3. 3 \rm ? $ $ \it If not, is there a simple relation between those two bijections \rm ?

\end{q5.3.}

\begin{q5.3.} $ $ \\
In the general case, are the two bijections of corollary 5.3. 3 isomorphisms of ordered sets, where the set $ $ $\{u \in W \mbox{ } | \mbox{ } u \mbox{ } \leq \mbox{ } v\}$ $ $ $(resp. $ $ \{u \in W \mbox{ } | \mbox{ } u \mbox{ } \leq \mbox{ } w\})$ $ $ is provided with the Bruhat order and $ $ $H$ $-$ $Spec(\upw)$ $ $ is provided with the inclusion?
\end{q5.3.}

\end{flushleft}

\end{document}